\documentclass[trsc,nonblindrev]{informs3}

\OneAndAHalfSpacedXI

\usepackage{graphicx}
\usepackage[shortlabels]{enumitem}
\usepackage{mathtools}
\DeclarePairedDelimiter{\ceil}{\lceil}{\rceil}
\usepackage{amsfonts}
\usepackage{caption}
\usepackage{subcaption}

\usepackage{algorithm}
\usepackage[noend]{algpseudocode}
\usepackage{multirow}

\usepackage{natbib}
 \bibpunct[, ]{(}{)}{,}{a}{}{,}%
 \def\bibfont{\small}%
 \def\BIBand{and}%

\TheoremsNumberedThrough 
\ECRepeatTheorems

\EquationsNumberedThrough 

\MANUSCRIPTNO{123456}

\begin{document}


\RUNAUTHOR{Hasan and Van Hentenryck}

\RUNTITLE{Commuting with Autonomous Vehicles}


\TITLE{\vspace{-1.5cm} \\
  Commuting with Autonomous Vehicles: A Branch and Cut Algorithm with Redundant Modeling}

\ARTICLEAUTHORS{
  \AUTHOR{Mohd. Hafiz Hasan} \AFF{University of Michigan, Ann Arbor, Michigan 48105, USA, \EMAIL{hasanm@umich.edu}}
  \AUTHOR{Pascal Van Hentenryck} \AFF{Georgia Institute of Technology, Atlanta, Georgia 30332, USA, \EMAIL{pvh@isye.gatech.edu}}
} 

\ABSTRACT{
  This paper studies the benefits of autonomous vehicles in
  ride-sharing platforms dedicated to serving commuting needs. It
  considers the Commute Trip Sharing Problem with Autonomous Vehicles
  (CTSPAV), the optimization problem faced by a reservation-based
  platform that receives daily commute-trip requests and serves them
  with a fleet of autonomous vehicles.  The CTSPAV can be viewed as a
  special case of the Dial-A-Ride Problem (DARP). However, this paper
  recognizes that commuting trips exhibit special spatial and temporal
  properties that can be exploited in a branch and cut algorithm that
  leverages a redundant modeling approach. In particular, the branch
  and cut algorithm relies on a MIP formulation that schedules mini
  routes representing inbound or outbound trips. This formulation is
  effective in finding high-quality solutions quickly but its
  relaxation is relatively weak. To remedy this limitation, the
  mini-route MIP is complemented by a DARP formulation which is not as
  effective in obtaining primal solutions but has a stronger
  relaxation. A column-generation procedure to compute the DARP
  relaxation is thus executed in parallel with the core branch and cut
  algorithm and asynchronously produces a stream of increasingly
  stronger lower bounds. The benefits of the proposed approach are
  demonstrated by comparing it with another, more traditional, exact
  branch and cut procedure and a heuristic method based on mini
  routes.

  The methodological contribution is complemented by a comprehensive
  analysis of a CTSPAV platform for reducing vehicle counts, travel
  distances, and congestion. In particular, the case study for a
  medium-sized city reveals that a CTSPAV platform can reduce daily
  vehicle counts by a staggering 92\% and decrease vehicles miles by
  30\%. The platform also significantly reduces congestion, measured
  as the number of vehicles on the road per unit time, by 60\% during
  peak times. These benefits, however, come at the expense of
  introducing empty miles. Hence the paper also highlights the
  tradeoffs between future ride-sharing and car-pooling platforms.
}%


\KEYWORDS{autonomous vehicles, shared commuting, branch-and-cut, column generation}

\maketitle

%


\section{Introduction}
\label{intro}

This work is the culmination of a four-year study on the benefits of
ride-sharing and car-pooling platforms for serving commuting needs. It
was originally motivated by the desire to relieve parking pressure in
the city of Ann Arbor, Michigan. Parking structures are expensive and
are often located in prime locations for the convenience of commuters.
In Ann Arbor, the parking pressure was primarily caused by commuters
to the University of Michigan, the city's largest employer with more
than 50,000 employees.

\begin{figure}[!t]
	\centering
        \includegraphics[width=1.0\linewidth]{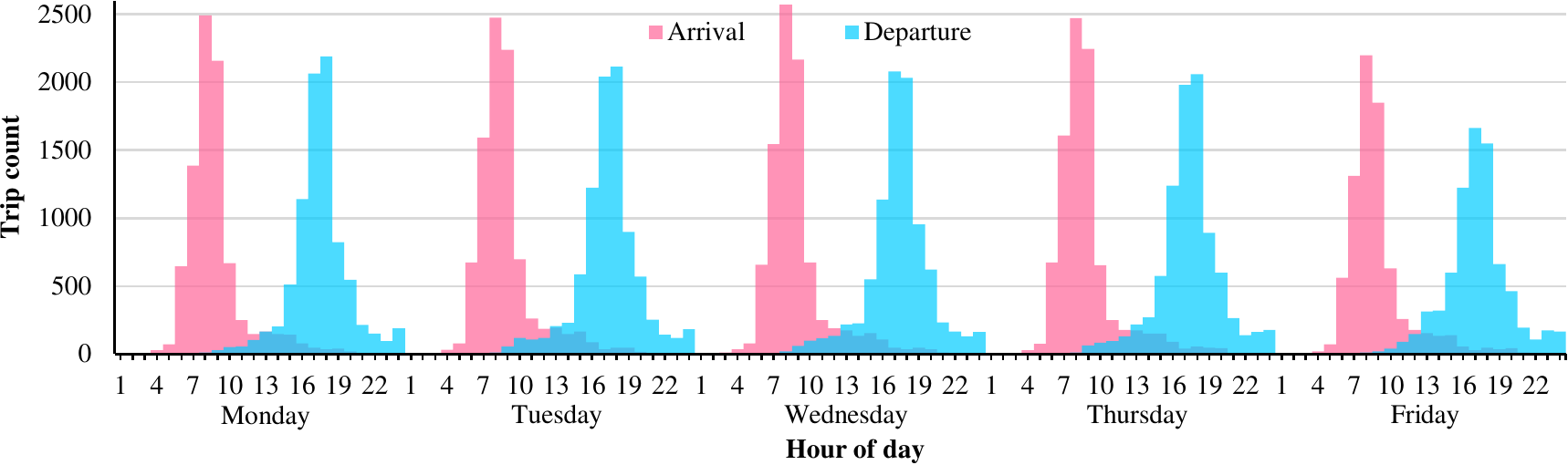}
	\caption{Distribution of Arrival and Departure Times Over Week
          2 of April 2017}
	\label{fig:arrival_departure_distribution}
\end{figure}

Detailed information about the commuting patterns of these employees
was gathered by recording trip data from approximately 15,000 drivers
who use the 15 university-operated parking structures located in the
downtown area over the month of April 2017. The data consisted of the
exact arrival and departure times of every commuter to the parking
structures, which was then joined with the precise locations of the
parking structures and the home addresses of every commuter to
reconstruct their daily trips. The dataset revealed several intriguing
temporal and spatial characteristics. First, the peak arrival and
departure times, which are depicted in Figure
\ref{fig:arrival_departure_distribution} for the weekdays of the
busiest week, coincide with the typical peak commuting
hours. Second, the strong consistency of the trip schedules was seen
as a significant opportunity for car-pooling and ride-sharing
platforms. Third, the commuting destinations (the parking structures)
are located within close vicinity of each other in the downtown area
(as they are university-owned structures), whereas the commuting
origins (the commuter homes) are located in the neighborhoods
surrounding the downtown area, as well as in Ann Arbor's neighboring
towns. This spatial structure, which is quite typical of many American
cities, was also seen as an opportunity for trip-sharing platforms.

With this in mind, \cite{hasan2020} introduced the Commute Trip
Sharing Problem (CTSP) to formalize the key optimization problem faced
by a car-pooling platform that would serve commute trips. More
precisely, the CTSP conceptualizes the platform as a reservation-based
system that receives the commute-trip requests---each consisting of a
trip request to the workplace (inbound trip) and another to return
back home (outbound trip)---ahead of time (e.g., the day ahead or the
morning of each day). Each trip request includes small time windows
for its departure and arrival, and each rider is guaranteed not to
spend more than $R\%$ of her direct trip in commuting time. The CTSP
was tailored to scenarios where: (1) The commuters travel to a
common/centralized location, e.g., the commute trips of the employees
of a corporate or university campus, or (2) The commuters live in a
common/centralized location, e.g., the commute trips originating from
a residential neighborhood or an apartment complex. These scenarios
were inspired by the spatio-temporal structure observed in the Ann
Arbor commute-trip dataset described earlier.

To implement such a platform and address the complexity of dealing
with the massive volume of the trips from the dataset,
\cite{hasan2020} applied a two-stage approach:
\begin{enumerate}
\item it first clusters commuters into artificial neighborhoods based
  on the spatial proximity of their home locations, using an
  unsupervised machine-learning algorithm;

\item it then finds optimal routes for the commuters within each cluster.
\end{enumerate}

\begin{figure}[!t]
	\centering \includegraphics[width=1.0\linewidth]{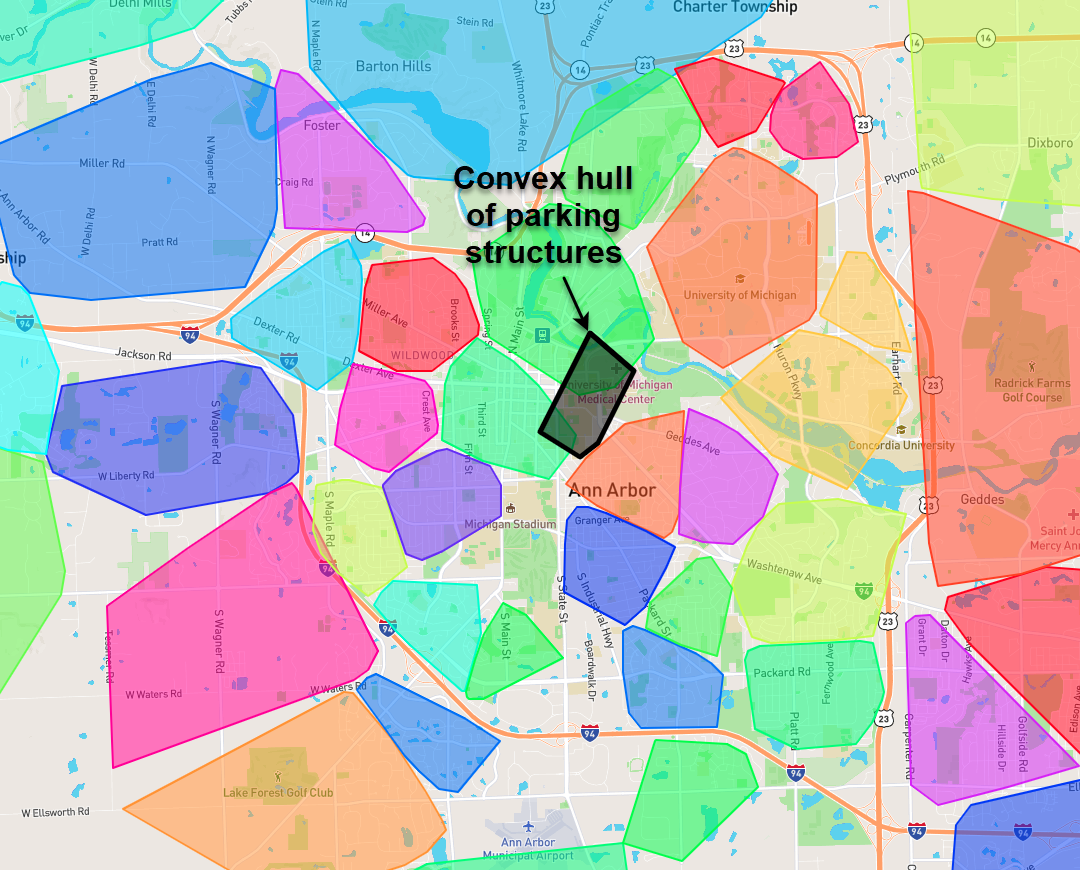}
	\caption{Convex Hulls of Artificial Neighborhoods Resulting
          from Clustering Algorithm}
	\label{fig:clusters}
\end{figure}

\noindent
Figure \ref{fig:clusters} provides an overview of the resulting
clusters within Ann Arbor's city limits: it displays the convex hulls
of the neighborhoods, as well as the convex hull of the centrally
located parking structures. The optimization problem in step 2 is the
CTSP: each day, its goal is to use private vehicles owned by the
commuters, select the set of drivers for the inbound and outbound
routes of the vehicles, and design the routes in order to minimize the
number of vehicles utilized, and hence the parking pressure. Solutions
to the CTSP were shown to reduce daily vehicle usage for the Ann Arbor
dataset by up to 57\%.

Despite this significant potential, the results also highlighted
several factors limiting further reductions in vehicle counts. They
included (1) the nature of the CTSP routes that are typically short
and (2) the necessity to synchronize the inbound and outbound routes
since they must be performed by the same set of drivers. Indeed, as the drivers
in the CTSP are selected from the set of commuters themselves, each
route must begin and end at the origin and destination of its
driver. This book-ending requirement subjects the total duration of
the route to the temporal constraints of the driver, restricting its
length and consequently its ability to serve more trips. This,
combined with the necessity of selecting an identical set of drivers
for the inbound and outbound routes, limits the flexibility of the
routes that can be generated and used in a CTSP routing plan.

The Commute Trip Sharing Problem with Autonomous Vehicles (CTSPAV)
considered in this paper was originally proposed by \cite{hasan2021}:
its goal was to overcome these shortcomings by leveraging Autonomous
Vehicle (AV) technology that is lurking in the horizon. By removing
driver-related constraints, the CTSPAV was anticipated to allow the AV
routes to be significantly longer than the CTSP routes. While these
longer routes would significantly increase the number of commute trips
that can be covered by each AV on any day, the algorithmic complexity
for finding them was also expected to increase
significantly. \cite{hasan2021} therefore proposed a column-generation
solution procedure, dubbed the CTSPAV procedure, that is a departure
from the classical column-generation approach for solving typical
Vehicle Routing Problems (VRPs). The latter typically entails solving
a set-partitioning/covering master problem that ensures that each
customer is served, and a pricing subproblem that searches for
feasible routes that depart from and return to a depot and have
negative reduced costs. The CTSPAV procedure circumvents the
anticipated complexity of searching for the long AV routes in the
pricing subproblem by shifting some of the burden to the master
problem and {\em exploiting the spatio-temporal structure of the
  dataset}. It uses a pricing subproblem that only searches for
feasible ``mini'' routes with negative reduced costs instead. The mini
routes are short by construction: each covers only inbound or outbound
trips exclusively, and each has distinct pickup, transit, and drop-off
phases during which it first visits trip origins, then travels from an
origin to a destination, and finally visits destinations. These three
phases are naturally encountered by each vehicle as it travels from a
residential neighborhood to the workplace in the morning to serve
inbound trips, and vice versa in the evening to serve outbound trips.

In order for these mini routes to be feasible, they must visit each
location within a specified time window, ensure that the time spent of
the vehicle by each rider does not exceed a specified limit, and
cannot exceed the vehicle capacity. In other words, they must satisfy
time-window, ride-duration, and vehicle-capacity
constraints. Furthermore, they must also satisfy pairing and
precedence constraints, which require a route visiting the origin of
trip to also visit its destination in the correct order. The master
problem of the CTSPAV procedure is then responsible for stitching or
chaining together the feasible mini routes to form longer AV routes
that begin and end at a depot. In addition to ensuring that each trip
is covered, the master problem must also select mini routes that are
temporally compatible with each other, i.e., it needs to ensure that
it is possible to travel from the last destination of one mini route
to the first origin of another without violating the temporal
constraints of the selected mini routes. All of this is done in
service of a lexicographic objective function that first minimizes the
number of formed AV routes (i.e., the vehicle count if each route is
assigned to an AV) and then minimizes their total travel distance.

Since the routes of the CTSPAV satisfy time-window, ride-duration,
capacity, pairing, and precedence constraints which are identical to
those for the Dial-A-Ride Problem (DARP)
\citep{cordeau2003a,cordeau2007}, the CTSPAV can be seen as a special
version of the DARP that serves inbound-outbound trip pairs using
AVs. In fact, any DARP algorithm can be used to solve the
CTSPAV. \cite{hasan2021} explored this possibility as well by
investigating a DARP procedure for solving the CTSPAV. The procedure
can be thought of as a model-driven approach that borrows heavily from
an algorithm for the DARP proposed by \cite{gschwind2015}, as it
relies on the classical column-generation approach but uses a novel,
label-setting dynamic program to solve its pricing subproblem.
\cite{hasan2021} discovered that, while the complexity of discovering
the long AV routes in its pricing subproblem severely hampered the
algorithm ability to find strong integer solutions within a
time-constrained setting, the DARP model also produced superior primal
lower bounds for the primary objective. On the other hand, the CTSPAV
procedure produces stronger integer solutions within a similar
time-constrained setting, but it does so at the expense of generating
weaker lower bounds.

This paper aims at addressing these limitations with two goals
in mind:
\begin{enumerate}
\item to propose an exact algorithm for the CTSPAV;
\item to provide a conclusive and comprehensive analysis of the
  potential of the CTSPAV for reducing vehicle counts, travel
  distances, and congestion.
\end{enumerate}
To meet the first goal, the paper presents an exact algorithm that
improves upon the CTSPAV procedure of \cite{hasan2021} by combining
the insights from both approaches in a redundant modeling framework
\citep{Liberti2004,Ruiz2011}. The proposed algorithm leverages the best
characteristics of the CTSPAV and DARP procedures, i.e., the former's
capability of producing strong integer solutions and the latter's
ability of generating strong primal lower bounds. More specifically,
the paper describes a branch-and-cut procedure which is capable of
solving medium-sized CTSPAV instances exactly, unlike the CTSPAV
procedure of \cite{hasan2021}.
This procedure is then compared against a branch-and-cut procedure
using other families of valid inequalities, as well as against the
CTSPAV procedure of \cite{hasan2021} for problem instances derived
from the Ann Arbor commute-trip dataset. With the exact CTSPAV
algorithm available, the paper can then perform a systematic analysis
of the CTSPAV potential in reducing vehicle counts, travel distance,
and congestion. Moreover, the paper can contrast the existing
situation where commuters drive mostly alone with car-pooling and
automomous ride-sharing platforms, highlighting the various trade-offs
on a real case study.

The {\em methodological contribution} of this paper is to propose a
branch-and-cut algorithm for solving the CTSPAV exploiting a novel
dual-modeling technique. The branch and cut algorithm solves a
mathematical model that exploits the spatio-temporal structure of the
data, making it conducive to finding high-quality solutions
quickly. But the branch and cut algorithm also uses another
mathematical model for the same problem to generate valid inequalities
that are separated by a column-generation procedure and produce strong
lower bounds. The paper demonstrates the benefits of this
dual-modeling approach through a comparison with a dedicated
branch-and-cut procedure based on well-established families of valid
inequalities, and with the heuristic column-generation procedure of
\cite{hasan2021}. The proposed exact branch and cut procedure is also
embedded into a end-to-end approach combining clustering and
optimization to solve large-scale, real-world instances of the CTSPAV.

The methodology ontribution is completemented by a case study that
provides unique insights on the potential benefits of ride sharing and
autonomous vehicles for serving the commuting needs of many cities
around the world. The case study demonstrates that a ride-sharing
platform based on autonomous vehicles can provide substantial
reductions in vehicle counts and congestion, as well as improvements
in travel miles. In addition, the paper contrasts, for the first time,
the potential benefits and drawbacks of car-pooling and ride-sharing
platforms along those dimensions.

The rest of this paper is organized as follows. Section
\ref{sec:relatedwork} briefly discusses related work.  Section
\ref{sec:notation} introduces the terminologies and assumptions used
throughout the work. Section \ref{sec:clustering} describes the
clustering algorithm.  Section \ref{sec:CTSPAV} specifies the CTSPAV
model and describes an algorithm for enumerating mini routes. Section
\ref{sec:validinequalities} provides an overview of the branch-and-cut
algorithm and covers the different families of valid inequalities
considered in this work together with the heuristics used to separate
them. Section \ref{sec:results} outlines how the algorithm is
evaluated and presents the computational results. Section
\ref{sec:casestudy} documents the insights obtained on the case study.
Finally, Section \ref{sec:conclusion} provides some concluding
remarks.

\section{Related Work} \label{sec:relatedwork}

The Vehicle Routing Problem with Time Windows (VRPTW) is perhaps the
most well-studied variant of VRPs; It seeks an optimal routing plan
that consists of a set of minimum cost routes, each departing and
returning to a designated depot, to service a set of customers. Each
customer has a capacity demand and a time window specifying allowable
service times, therefore the plan must ensure every customer is served
exactly once within their time windows while not exceeding the
capacity of the vehicles utilized, i.e., its routes must satisfy
time-window and vehicle-capacity constraints. The problem is
well-known to be NP-hard as finding a feasible solution to the version
of the problem with a fixed vehicle count has been shown to be
NP-complete by \cite{savelsbergh1985}. Nevertheless, numerous
approaches ranging from heuristics to exact methods have been proposed
for the problem, and they have been comprehensively reviewed by
\cite{cordeau2002}. The VRPTW was generalized to the Pickup and
Delivery Problem with Time Windows (PDPTW) by \cite{dumas1991} to
model services that first pick up and then deliver merchandise within
specified time windows. The routes of the problem therefore need to
satisfy pairing and precedence constraints in addition to time-window
and vehicle-capacity. The former two require that each route visit a
pair of locations associated with each customer in a specific order,
the first representing a pickup location and the second representing a
delivery location. The PDPTW was then generalized to the DARP which is
used to model door-to-door transportation services for the disabled or
the elderly. The ride duration becomes a critical factor for ensuring
the quality of these services as they are now transporting
humans. Therefore the DARP introduces ride-duration constraints to the
PDPTW, which limit the time elapsed between every pair of pickup and
delivery location to ensure that the customers are not spending
excessive amounts of time on the vehicle. The various algorithms and
techniques that have been proposed for the DARP have been reviewed by
\cite{cordeau2003a,cordeau2007}.

Of the many solution approaches that have been proposed for the
different variants of the VRP, column generation is perhaps to most
popular due to its ability to generate strong lower bounds to the
problem objective and due to its elegance of only considering a subset
of feasible routes that can improve the objective function. The
typical column-generation approach for solving VRPs begins with the
application of the Dantzig-Wolfe decomposition \citep{dantzig1960} on
an edge-flow formulation of the problem to produce a master problem
and a pricing subproblem. The master problem typically solves a
set-partitioning/covering problem on a set of feasible routes to
ensure every customer is served, whereas the pricing subproblem
searches for new feasible routes to be added to the set. The latter
problem uses the duals of the linear relaxation of the master problem
to identify new routes with negative reduced costs, and it is
typically cast as a Shortest Path Problem with Resource Constraints
(SPPRC), a class of problems that has been extensively reviewed by
\cite{irnich2005}. The SPPRC seeks a route with minimum cost, and the
feasibility of the discovered route is guaranteed through the
enforcement of numerous resource constraints that model the
route-feasibility constraints. Some of the approaches that have been
used to solve these SPPRCs include Lagrangian relaxation
\citep{beasley1989,borndrfer2001}, constraint programming
\citep{rousseau2004}, heuristics \citep{desaulniers2008}, and cutting
planes \citep{drexl2013}, but perhaps the most popular approach uses
dynamic programming, e.g., the generalized label-setting algorithm for
multiple resource constraints by \cite{desrochers1988b}. Examples of
successful applications of column generation on the different VRP
variants include \cite{desrosiers1984,desrochers1992} for the VRPTW,
\cite{dumas1991,ropke2009} for the PDPTW, and \cite{gschwind2015} for
the DARP.

Another common approach for solving routing problems is the polyhedral
approach which generates cutting planes to progressively ``trim'' the
convex hull defining the feasible region of the problem's linear
relaxation. Its application on VRPs traces its roots back to the
seminal work by \cite{dantzig1954} for solving the Traveling Salesman
Problem (TSP). Their procedure uses an edge-flow formulation of the
problem which is iteratively solved to identify subtours which break
the feasibility of the solution. A family of valid inequalities,
commonly referred to now as the DFJ subtour elimination constraints
(SECs), are then progressively introduced to prevent generation of the
subtours in subsequent solutions. \cite{grotshcel1975} later proved
that the DFJ SECs induce facets of the polytope of the convex hull of
the feasible solutions, which explained why they were so effective at
strengthening the linear-programming (LP) bound, while
\cite{padberg1990} proposed an exact algorithm for separating the
inequalities. In a similar vein, many other works have focused on
identifying facet-defining inequalities together with
algorithms/heuristics for separating them, e.g., $D_k^+$ and $D_k^-$
inequalities for the TSP by \cite{grotshcel1985}, predecessor and
successor inequalities for the Precedence-Constrained Asymmetric TSP
(PCATSP) by \cite{balas1995}, tournament and generalized tournament
constraints for the Asymmetric TSP with Time Windows (ATSPTW) by
\cite{ascheuer2000}, and 2-path cuts for the VRPTW by
\cite{kohl1999}. Most approaches to routing problems embed
cutting-plane generation within the classical branch-and-bound
framework for solving mixed-integer programs (MIPs) to produce a more
sophisticated branch-and-cut procedure, whereby heuristics for
separating violated valid inequalities are executed on the solution of
the LP relaxation that is obtained in the bounding phase of each tree
node. The separated inequalities are then introduced into the problem
formulation to strengthen the LP bound of the procedure. The proposed
branch-and-cut algorithms typically begin with an edge-flow
formulation and then introduce numerous existing and/or new families
of valid inequalities that are tailored specifically for the type of
routing problem being solved. Examples of these branch-and-cut
algorithms include \cite{padberg1991} for the TSP,
\cite{fischetti1997} for the Asymmetric TSP (ATSP), \cite{ruland1997}
for the Pickup and Delivery Problem (PDP), \cite{ascheuer2001} for the
ATSPTW, \cite{naddef2001} for the Capacitated VRP (CVRP),
\cite{bard2002,kallehauge2007} for the VRPTW, and \cite{cordeau2006}
for the DARP.

The prevalence of large-scale datasets of real-world trips, e.g., the
New York City (NYC) Taxi \& Limousine Commission (TLC) trip record
data \citep{nyctriprecord} which stores trip information of more than
one billion taxi rides in NYC, combined with the growing awareness and
concern for the sustainability of passenger mobility systems have
increased attention towards the optimization of car-pooling and
ridesharing services. For instance, \cite{santi2014} formalized the
notion of shareability networks as a tool to quantify the ridesharing
potential of the trips from the TLC dataset, while
\cite{alonso-mora2017} proposed an anytime optimal algorithm that
utilizes shareability graphs to optimize ridesharing for on-demand
trip requests extracted from the TLC dataset. Studies involving other
real-world datasets include \cite{baldacci2004} who proposed a
Lagrangian column-generation method to optimize the Car-Pooling
Problem (CPP) for commuting trips to a research institution in Italy
and \cite{agatz2011} who used graph matching within a rolling-horizon
framework to optimize ridesharing for real-time, non-recurring trips
from metro Atlanta. Classifications of the different variants of
shared mobility problems together with reviews of the proposed
optimization approaches for them are provided by \cite{agatz2012} and
\cite{mourad2019}. The impending arrival of fully autonomous vehicles
has also spurred a growing interest in the potential of Shared
Autonomous Vehicle (SAV) services, due to the perceived benefits that
are afforded by this new mode of transportation, be it reducing
traffic \citep{martinez2017,alazzawi2018,salazar2018}, increasing road
capacity
\citep{friedrich2015,tientrakool2011,talebpour2016,menaoreja2018,olia2018},
or decreasing parking demand
\citep{zhang2015,dia2017,zhang2017}. \cite{narayanan2020}, which
reviewed the numerous potential impacts of SAV services to society and
the environment, also suggested classifying them as either on-demand
or reservation-based systems, with the former being tailored for
dynamic trips whose requests are made in real time and the latter for
recurring trips whose requests are made way in advance. Several
optimization approaches have also been proposed for conceptual systems
of each type. For example, \cite{farhan2018} proposed a three-step
approach---which clusters trip requests from discretized time
intervals by assigning them to their nearest vehicles and then solving
the requests for each cluster as a VRPTW---to optimize a fleet of SAVs
for on-demand trips, while \cite{ma2017} proposed an LP approach to
optimize vehicle sharing of a fleet of SAVs for trip requests that are
known ahead of time.

The work on the CTSPAV traces its roots back to the authors' initial
desire to solve the parking problem in downtown Ann Arbor, Michigan,
that was partly caused by the massive infusion of trips from the
thousands of commuters driving to the University of Michigan campus
daily. Having access to a large-scale, high-fidelity dataset of these
commute trips, they wanted to investigate the vehicle reduction
potential of an optimized car-pooling or ridesharing platform.
\cite{hasan2018} began by investigating the performance of several
car-pooling and car-sharing models, each with different driver and
passenger matching constraints, and discovered that the model that
requires the commuters to adopt different roles and to ride with
different passengers and drivers daily had the best vehicle reduction
potential. In other words, the flexibility in driving and sharing
preferences is critical to maximizing trip shareablity. In
\citep{hasan2020}, the best performing car-pooling model was refined
and subsequently formalized as the CTSP, a model that maximizes trip
sharing while selecting an identical set of drivers for the inbound
and outbound routes from the set of commuters on a daily basis. Two
exact algorithms were proposed: the first exhaustively enumerates
feasible routes before their selection is optimized with a MIP, while
the second uses column generation to search for feasible routes on
demand within a branch-and-price framework. Subsequent application of
the algorithms on the commute-trip dataset revealed an ability to
reduce daily vehicle counts by more than 50\%. \cite{hasan2020b} then
proposed a method to handle potential uncertainties in the trip
schedules of the CTSP by incorporating a randomized, scenario-sampling
technique within a two-stage optimization approach. The method was
shown to be capable of producing routing plans that are robust to
changes in trip schedules, but the increase in robustness comes at the
price of an increase in vehicle utilization. A method to properly
evaluate this trade-off was then proposed. The CTSPAV was formally
conceptualized in \cite{hasan2021} to address a key shortcoming of the
CTSP---its short routes which limited the potential to further reduce
daily vehicle counts---through the utilization of a SAV platform. The
work explored two methods for optimizing its routes: (1) an approach
which uses column generation to search for mini routes which are then
assembled in a master problem, and (2) an approach which relied on a
more classical column-generation technique originally conceived for
the DARP. They discovered that each method had complementary
performance trade-offs, with the former being able to produce stronger
integer solutions and the latter being able to generate stronger lower
bounds. All of these earlier works have culminated into this study
which hopes to develop an algorithm that melds together both
approaches proposed from \cite{hasan2021} in order to leverage their
unique strengths in effectively solving the CTSPAV. Accomplishing this
goal uniquely positions this work to glean additional insights into
the strengths and weaknesses of an optimized SAV platform relative to
car-pooling platforms that uses conventional vehicles for maximizing
large-scale ridesharing of commute trips.

\section{Preliminaries} \label{sec:notation}

This section introduces the main concepts used throughout this paper:
trips, mini routes, and AV routes. It also describes the constraints
that mini routes and AV routes must satisfy. This work assumes that a
homogeneous fleet of vehicles with capacity $K$ is available to serve
all rides, and that the triangle inequality is satisfied for all
travel times.

\paragraph{Trips} A trip $t=\{o,dt,d,at\}$ is a tuple that consists of an origin $o$, 
a departure time $dt$, a destination $d$, and an arrival time $at$ of
a trip request.  Every day, a commuter $c$ makes two trips: a trip
$t_c^+$ to the workplace and a return trip $t_c^-$ back home. These
trips are called inbound and outbound trips respectively.

\paragraph{Mini Routes} A mini route $r$ is a sequence of locations that visits each origin
and destination from a set of inbound or outbound trips exactly
once. Let $\mathcal{C}_r$ denote the set of riders served in $r$. A
mini route $r$ must respect the vehicle capacity, i.e.,
$|\mathcal{C}_r| \leq K$, and consists of three phases: a pickup phase
where the passengers are picked up, a transit phase where the vehicle
travels to the destination, and a drop-off phase where all the
passengers are dropped off. During the pickup (resp., drop-off) phase,
the vehicle visits only origins (resp., destinations), whereas it
travels from an origin to a destination in the transit phase. For
instance, a possible mini route for a car with $K=4$ serving trips
$t_1=\{o_1,dt_1,d_1,at_1\}$, $t_2=\{o_2,dt_2,d_2,at_2\}$, and
$t_3=\{o_3,dt_3,d_3,at_3\}$ is $r=o_2\rightarrow o_1\rightarrow
o_3\rightarrow d_1\rightarrow d_2\rightarrow d_3$, and its pickup,
transit, and drop-off phases are given by $o_2\rightarrow
o_1\rightarrow o_3$, $o_3\rightarrow d_1$, and $d_1\rightarrow
d_2\rightarrow d_3$ respectively. An inbound mini route $r^+$ covers
only inbound trips and an outbound mini route $r^-$ covers only
outbound trips.

\begin{definition}[Valid Mini Route]
	A valid mini route $r$ serving a set of riders $\mathcal{C}_r$
        visits all of its origins, $\{o_c:c\in\mathcal{C}_r\}$, before
        its destinations, $\{d_c:c\in\mathcal{C}_r\}$, and respects
        the vehicle capacity, i.e., it has $|\mathcal{C}_r|\leq K$.
\end{definition}

\noindent
Let $T_i$ denote the time at which service begins at location $i$,
$s_i$ the service duration at $i$, $pred(i)$ the location visited just
before $i$, $\tau_{(i,j)}$ the estimated travel time for the shortest
path between locations $i$ and $j$, and $\dot{\mathcal{C}}_r$ the
first commuter served on $r$.  Commuters sharing rides are willing to
tolerate some inconvenience in terms of deviations to their desired
departure and arrival times, as well as in terms of their ride
durations compared to their individual, direct trips. Therefore, a
time window $[a_i,b_i]$ is constructed around the desired departure
times and is associated with each pickup location $i$, where $a_i$ and
$b_i$ denote the earliest and latest times at which service may begin
at $i$ respectively. Conversely, only an upper bound $b_j$ is
associated with each drop-off location $j$ as the arrival time at $j$
is implicitly bounded from below by $a_j = a_i + s_i + \tau_{(i,j)}$,
where $i$ is the corresponding pickup location for $j$. On top of
that, a duration limit $L_c$ is associated with each rider $c$ to
denote her maximum ride duration.

\begin{definition}[Feasible Mini Route]
	A feasible mini route $r$ is valid, has pickup and drop-off
        times $T_i\in [a_i,b_i]$ for each location $i\in r$, and
        ensures the ride duration of each rider $c\in \mathcal{C}_r$
        does not exceed $L_c$.
\end{definition}

\noindent
Determining if a valid mini route $r$ is feasible amounts to solving a
feasibility problem defined by the following constraints on $r$.
\begin{flalign}
	& a_{o_c}\leq T_{o_c} \leq b_{o_c} \qquad \forall c\in
  \mathcal{C}_r \label{eqn:time_window_origin} \\ & T_{d_c} \leq
  b_{d_c} \qquad \forall c\in
  \mathcal{C}_r \label{eqn:time_window_destination} \\ & T_{pred(o_c)}
  + s_{pred(o_c)} + \tau_{(pred(o_c),o_c)} \leq T_{o_c} \qquad \forall
  c\in
  \mathcal{C}_r\setminus\dot{\mathcal{C}}_r \label{eqn:travel_time_origin}
  \\ & T_{pred(d_c)} + s_{pred(d_c)} + \tau_{(pred(d_c),d_c)} =
  T_{d_c} \qquad \forall c\in
  \mathcal{C}_r \label{eqn:travel_time_destination} \\ & T_{d_c} -
  (T_{o_c} + s_{o_c}) \leq L_c \qquad \forall c\in
  \mathcal{C}_r \label{eqn:ride_duration_limit}
\end{flalign}

\noindent
Constraints \eqref{eqn:time_window_origin} and
\eqref{eqn:time_window_destination} are time-window constraints for
pickup and drop-off locations respectively, while constraints
\eqref{eqn:travel_time_origin} and \eqref{eqn:travel_time_destination}
describe compatibility requirements between pickup/drop-off times and
travel times between consecutive locations along the route. Finally,
constraints \eqref{eqn:ride_duration_limit} specify the ride-duration
limit for each rider. Note that constraints
\eqref{eqn:travel_time_origin} allow waiting at pickup locations.
Moreover, the service starting times on consecutive locations along
$r$ are strictly increasing, which ensures that the route is
elementary.  Numerous algorithms have been proposed for solving this
feasibility problem efficiently,
e.g. \cite{tang2010,haughland2010,firat2011,gschwind2015}. In the
following, the Boolean function $feasible(r)$ is used to indicate
whether mini route $r$ admits a feasible solution to constraints
\eqref{eqn:time_window_origin}--\eqref{eqn:ride_duration_limit}. This
work implements the labeling procedure proposed by \cite{gschwind2015}
for this function.

\paragraph{AV Routes} An AV route $\rho=v_s\rightarrow r_1\rightarrow \ldots\rightarrow
r_k\rightarrow v_t$ is a sequence of $k$ distinct mini routes that
starts at a source node $v_s$ and ends at a sink node $v_t$, both
representing a designated depot.

\begin{definition}[Feasible AV Route]
	A feasible AV route $\rho$ is one that consists of a sequence
        of distinct, feasible mini routes and starts and ends at a
        designated depot.
\end{definition}

\noindent
In other words, for $\rho$ to be feasible, each of its mini routes
must be valid and satisfy constraints
\eqref{eqn:time_window_origin}--\eqref{eqn:ride_duration_limit}. Let
$\dot{r}$ denote the first location visited on $r$ and $\ddot{r}$
denote the last. Each mini route $r_i \; (1 \leq i \leq k)$ must also
satisfy the following constraints:
\begin{flalign}
	& T_{v_s} + \tau_{(v_s, \dot{r}_1)} =
  T_{\dot{r}_1} \label{eqn:travel_time_begin} \\ & T_{\ddot{r}_i} +
  s_{\ddot{r}_i} + \tau_{(\ddot{r}_i, \dot{r}_{i+1})} \leq
  T_{\dot{r}_{i+1}} \qquad \forall i = 1, \ldots,
  k-1 \label{eqn:travel_time_routes} \\ & T_{\ddot{r}_k} +
  s_{\ddot{r}_k} + \tau_{(\ddot{r}_k, v_t)} =
  T_{v_t} \label{eqn:travel_time_end}
\end{flalign}

\noindent
Constraints \eqref{eqn:travel_time_begin}--\eqref{eqn:travel_time_end}
describe compatibility requirements between the beginning/ending
service times of consecutive mini routes along $\rho$ and the travel
times between them. The constraints, together with
\eqref{eqn:travel_time_origin} and
\eqref{eqn:travel_time_destination}, enforce strictly increasing
starting times for service on all consecutive locations along $\rho$,
therefore ensuring that $\rho$ is elementary.

\section{The Clustering Algorithm} \label{sec:clustering}

This section describes a clustering algorithm used to decompose the
large volume of commute trips in our case study into smaller, more
manageable problem instances. This strategy is congruent with the
conclusion of \cite{agatz2012} that acknowledges the necessity of
effective decomposition approaches for the computational feasibility
of large-scale problems. The idea behind this clustering approach is
simply to construct artificial neighborhoods within which ridesharing
is performed exclusively, and the neighborhoods are constructed by
algorithmically grouping up to $N$ commuters together based on the
spatial proximity of their residential locations. Obviously, this
approach precludes the discovery of a global optimal solution, but it
is seen as a practical necessity to ensure that the problem is
computationally tractable.

The algorithm proceeds in a fashion that is very similar to the
$k$-means clustering algorithm by \cite{lloyd1982}, with the exception
that its assignment step limits the number of elements assigned to
each cluster by a parameter $N$ to produce groups that are
approximately equal in size. It represents each commuter as a point in
$\mathbb{R}^2$ whose GPS coordinates are first obtained by geocoding
the commuter home address. In the rest of this section, $\mathcal{C}$
denotes the set of point coordinates for every commuter (i.e., a set
of 2D vectors, each storing the 2D coordinates of a commuter home),
$\mathcal{U}$ the set of coordinates of cluster centers (similarly, a
set of 2D vectors, each consisting of the 2D coordinates of a cluster
center), $S(\boldsymbol{x})$ the Euclidean distance from a point
$\boldsymbol{x}$ to the nearest cluster center, and
$S(\boldsymbol{x},\boldsymbol{y})$ the Euclidean distance between
points $\boldsymbol{x}$ and $\boldsymbol{y}$.

The algorithm begins with the identification of $k$, the number of
clusters, using $k = \ceil{|\mathcal{C}|/N}$. The $k$ cluster centers
are then initialized randomly using the {\tt $k$-means++} method by
\cite{arthur2007}. The method first selects a point uniformly at
random from $\mathcal{C}$ as the first center, $\boldsymbol{u}_1$, and
then selects the $i^\text{th}$ center, $\boldsymbol{u}_i$, from
$\mathcal{C}$ with probability $S(\boldsymbol{u}_i)^2 /
[\sum_{\boldsymbol{c}\in\mathcal{C}} S(\boldsymbol{c})^2]$ until $k$
centers are selected. Each point $\boldsymbol{c}\in\mathcal{C}$ is
then assigned to its nearest cluster center subject to the constraint
that each center is assigned at most $N$ points. This assignment step
is accomplished by solving the generalized-assignment problem
described in Figure \ref{fig:clustering}. The formulation uses a
binary variable $x_{\boldsymbol{c},\boldsymbol{u}}$ that indicates a
point $\boldsymbol{c}$ is assigned to center $\boldsymbol{u}$ when
set. Its objective function \eqref{eqn:min_total_distance} minimizes
the total distance between all points and their assigned
centers. Constraints \eqref{eqn:assignment} assign each point to a
center, while constraints \eqref{eqn:size_limit} limit the number of
points assigned to each center by $N$.

\begin{figure}[!t]
\begin{flalign}
	\min & \quad \sum_{\boldsymbol{c}\in\mathcal{C}}\sum_{\boldsymbol{u}\in\mathcal{U}} S(\boldsymbol{c},\boldsymbol{u}) x_{\boldsymbol{c},\boldsymbol{u}} & \label{eqn:min_total_distance} \\
	\text{s.t.}\nonumber & \\
	& \quad \sum_{\boldsymbol{u}\in\mathcal{U}} x_{\boldsymbol{c},\boldsymbol{u}}=1\qquad\forall \boldsymbol{c}\in\mathcal{C} \label{eqn:assignment} \\
	& \quad \sum_{\boldsymbol{c}\in\mathcal{C}} x_{\boldsymbol{c},\boldsymbol{u}}\leq N\qquad\forall \boldsymbol{u}\in\mathcal{U} \label{eqn:size_limit} \\
	& \quad x_{\boldsymbol{c},\boldsymbol{u}}\in\{0,1\}\qquad\forall \boldsymbol{c}\in\mathcal{C},\forall \boldsymbol{u}\in\mathcal{U} \label{eqn:assign_var}
\end{flalign}
\caption{The Clustering Formulation.}
\label{fig:clustering}
\end{figure}

The assignment step is followed by an update step which recalculates
the coordinates of each cluster center by averaging the coordinates of
its assigned points:
\begin{equation}
	\boldsymbol{u} = \frac{\sum_{\boldsymbol{c}\in\mathcal{C}} x_{\boldsymbol{c},\boldsymbol{u}} \boldsymbol{c}}{\sum_{\boldsymbol{c}\in\mathcal{C}}x_{\boldsymbol{c},\boldsymbol{u}}} \qquad \forall \boldsymbol{u}\in\mathcal{U}
\end{equation}
The assignment and update steps are then repeated until the
point-center assignments stabilize, i.e., until the centers every
point are assigned to remain the same in consecutive iterations.

\section{The Commute Trip Sharing Problem for Autonomous Vehicles} \label{sec:CTSPAV}

This section specifies the CTSPAV, a problem which seeks a set of
minimal cost AV routes to serve every inbound and outbound trip of a
set of commuters, $\mathcal{C}$.

\subsection{Notation}
Let $n$ denote the total number of commuters, i.e., $n =
|\mathcal{C}|$. For every commuter $i \in \mathcal{C}$, let nodes $i$,
$n+i$, $2n+i$, and $3n+i$ represent the inbound pickup, inbound
drop-off, outbound pickup, and outbound drop-off locations of the
rider trips respectively. Then let the sets of all inbound pickup, all
inbound drop-off, all outbound pickup, and all outbound drop-off nodes
be denoted by $\mathcal{P}^+ = \{1,\ldots,n\}$, $\mathcal{D}^+ =
\{n+1,\ldots,2n\}$, $\mathcal{P}^- = \{2n+1,\ldots,3n\}$, and
$\mathcal{D}^- = \{3n+1,\ldots,4n\}$ respectively. Furthermore, let
$\mathcal{P} = \mathcal{P}^+\cup\mathcal{P}^-$ and $\mathcal{D} =
\mathcal{D}^+\cup\mathcal{D}^-$. With this notation, note that $n+i$
provides the drop-off node corresponding to any pickup node
$i\in\mathcal{P}$. By definition of AV routes, the following
precedence constraints apply to the following set of nodes:
\begin{equation} 
	i \prec n + i \prec 2n + i \prec 3n + i \qquad \forall i\in\mathcal{P}^+ \label{eqn:precedence}
\end{equation}
where $i \prec j$ denotes the precedence relation between nodes $i$
and $j$, i.e., the constraint indicating that $i$ must be visited
before $j$ if both $i$ and $j$ are served by the same AV route.

The directed graph $\mathcal{G} = (\mathcal{N}, \mathcal{A})$ with the
node set $\mathcal{N} = \mathcal{P}\cup\mathcal{D}\cup\{v_s,v_t\}$
contains all pickup and drop-off nodes together with a source and a
sink node (both representing the designated depot) and its edge set
$\mathcal{A}=\{(i,j):i,j\in\mathcal{N},i\neq j\}$ consists of all
possible edges as a first approximation. A time window $[a_i,b_i]$ and
a service duration $s_i$ are then associated with each node
$i\in\mathcal{P}\cup\mathcal{D}$. No time windows are associated with
$v_s$ and $v_t$ as it is assumed that the AVs may start and end their
routes at any time of the day. Additionally, a ride-duration limit
$L_i$ is associated with each node $i\in\mathcal{P}$. Finally, a
travel time $\tau_{(i,j)}$, a distance $\varsigma_{(i,j)}$, and a cost
$c_{(i,j)}$ are associated with each edge $(i,j)\in\mathcal{A}$, and
$\delta^+(i)$ and $\delta^-(i)$ denote the sets of all outgoing and
incoming edges of node $i$ respectively.

\subsection{A MIP Model for the CTSPAV} \label{sec:MIP}

This section introduces a MIP model for the CTSPAV. The MIP is
summarized in Figure \ref{fig:MIP}: it formalizes the CTSPAV and is
defined on the graph $\mathcal{G}$ and {\em the set
$\mathrm{\Omega}$ of all feasible mini routes}. The MIP formulation uses
two sets of binary variables: variable $X_r$ indicates whether mini
route $r\in\mathrm{\Omega}$ is selected and variable $Y_{(i,j)}$
indicates whether edge $(i,j)\in\mathcal{A}$ is used, i.e., whether
node $j$ should be visited immediately after node $i$ by an AV route
in the optimal solution. Additionally, the model uses a continuous
variable $T_i$ that represents the start of service time at node
$i\in\mathcal{P}\cup\mathcal{D}$.

The objective function \eqref{eqn:ctspav_obj} minimizes the total cost
of all selected edges. Contraints \eqref{eqn:ctspav_route_cover}
ensure each trip is served by exactly one mini route, while
constraints \eqref{eqn:ctspav_edge_select} select edges belonging to
selected mini routes. Constraints \eqref{eqn:ctspav_outgoing} and
\eqref{eqn:ctspav_incoming} simultaneously ensure each pickup and
drop-off node is visited exactly once while conserving the flow
through each. Constraints \eqref{eqn:ctspav_travel_time_1} and
\eqref{eqn:ctspav_travel_time_2} ensure the start of service time at
the tail and head of every selected edge is compatible with the travel
time along the edge using large constants $M_{(i,j)}$ and
$\bar{M}_{(i,j)}$. Finally, constraints
\eqref{eqn:ctspav_ride_duration} and \eqref{eqn:ctspav_time_window_1}
describe the ride-duration limit of every trip and the time-window
constraint of every pickup and drop-off node respectively.

Note that constraints \eqref{eqn:ctspav_travel_time_1} and
\eqref{eqn:ctspav_travel_time_2} are generalizations of the popular
Miller-Tucker-Zemlin (MTZ) subtour-elimination constraints for the TSP
\citep{miller1960}. They utilize big-$M$ constants and enforce the
underlying constraints on a subset of edges:
\begin{flalign}
	& M_{(i,j)} = \max \{0, b_i + s_i + \tau_{(i,j)} - a_j\} \qquad \forall i,j\in\mathcal{P}\cup\mathcal{D} \\
	& \bar{M}_{(i,j)} = \max \{0, b_j - a_i - s_i - \tau_{(i,j)}\} \qquad \forall i\in\mathcal{P}\cup\mathcal{D}, \forall j\in\mathcal{D}
\end{flalign}

\begin{figure}[!th]
\begin{flalign}
	& \min \sum_{e\in\mathcal{A}} c_e Y_e \label{eqn:ctspav_obj} \\
	& \text{subject to} \nonumber \\
	& \sum_{r\in\mathrm{\Omega}:i\in r} X_r = 1 \qquad \forall i\in\mathcal{P} \label{eqn:ctspav_route_cover} \\
	& \sum_{r\in\mathrm{\Omega}:e\in r} X_r - Y_e \leq 0 \qquad \forall e\in\mathcal{A} \setminus \{\delta^+(v_s)\cup\delta^-(v_t)\} \label{eqn:ctspav_edge_select} \\
	& \sum_{e\in\delta^+(i)} Y_e = 1 \qquad \forall i\in\mathcal{P}\cup\mathcal{D} \label{eqn:ctspav_outgoing} \\
	& \sum_{e\in\delta^-(i)} Y_e = 1 \qquad \forall i\in\mathcal{P}\cup\mathcal{D} \label{eqn:ctspav_incoming} \\
	& T_i + s_i + \tau_{(i,j)} \leq T_j + M_{(i,j)}(1 - Y_{(i,j)}) \qquad \forall i,j\in\mathcal{P}\cup\mathcal{D} \label{eqn:ctspav_travel_time_1} \\
	& T_i + s_i + \tau_{(i,j)} \geq T_j - \bar{M}_{(i,j)}(1 - Y_{(i,j)}) \qquad \forall i\in\mathcal{P}\cup\mathcal{D}, \forall j\in\mathcal{D} \label{eqn:ctspav_travel_time_2} \\
	& T_{i+n} - (T_i + s_i) \leq L_i \qquad \forall i\in\mathcal{P} \label{eqn:ctspav_ride_duration} \\
	& a_i \leq T_i \leq b_i \qquad \forall i\in\mathcal{P}\cup\mathcal{D} \label{eqn:ctspav_time_window_1} \\
	& X_r \in\{0,1\} \qquad \forall r\in\mathrm{\Omega} \label{eqn:ctspav_route_var} \\
	& Y_e \in\{0,1\} \qquad \forall e\in\mathcal{A} \label{eqn:ctspav_edge_var} 
\end{flalign}
\caption{The MIP Model for the CTSPAV.}
\label{fig:MIP}
\end{figure}

The model adopts a lexicographic objective whose primary objective is
to minimize the number of vehicles used and whose secondary objective
is to minimize the total travel distance. This lexicographic ordering
is accomplished by weighting the sub-objectives: an identical, large
fixed cost and a variable cost that is proportional to the route total
distance are assigned to each AV route. The edge costs are defined as
follows to accomplish this goal:
\begin{equation}
	c_e=
	\begin{cases} \label{eqn:edge_costs}
		\varsigma_e + 100 \hat{\varsigma}_\text{max}&\qquad\forall e\in\delta^+(v_s)\\
		\varsigma_e&\qquad\text{otherwise}
	\end{cases}
\end{equation} 
where $\hat{\varsigma}_\text{max}$ is a constant equal to the length (total distance) of the longest AV route. Letting $\mathcal{R}$ denote the set of all feasible AV routes, $\hat{\varsigma}_\text{max}$ is given by:
\begin{equation}
	\hat{\varsigma}_\text{max} = \max_{\rho\in\mathcal{R}} \sum_{(i,j)\in\rho} \varsigma_{(i,j)}
\end{equation}

\noindent
The CTSPAV model essentially solves a scheduling problem that selects
and assembles feasible mini routes to form longer, feasible AV routes
to cover all trips and minimize the total cost. The optimal AV routes
are obtained by constructing paths beginning at $v_s$ and ending at
$v_t$ from the selected edges, and the start and end times can be
calculated using Equations \eqref{eqn:travel_time_begin} and
\eqref{eqn:travel_time_end} respectively.

\subsection{The Mini Route-Enumeration Algorithm}

Since the MIP model is defined in terms of all mini-routes, this
section describes the Mini Route-Enumeration Algorithm (MREA), a
procedure for enumerating all the feasible mini routes in
$\mathrm{\Omega}$ that is based on the algorithm proposed by
\cite{hasan2020}. The set $\mathrm{\Omega}$ can be partitioned into
$\mathrm{\Omega} = \mathrm{\Omega}^+ \cup \mathrm{\Omega}^-$, where
$\mathrm{\Omega}^+$ represents the set of feasible inbound mini routes
(which covers only inbound trips) while $\mathrm{\Omega}^-$ represents
the set of feasible outbound mini routes (which covers only outbound
trips). Without loss of generality, this section only describes the
procedure for enumerating the mini routes in $\mathrm{\Omega}^+$.

The procedure is summarized in Algorithm \ref{alg:MREA}. It requires
as inputs the set $\mathcal{T}^+$ of all inbound trips and the vehicle
capacity $K$. It begins by considering all feasible inbound mini
routes for a vehicle capacity of 1 by adding the routes for all direct
trips from $\mathcal{T}^+$ to $\mathrm{\Omega}^+$ (lines 2--3). It
then enumerates feasible routes for progressively increasing vehicle
capacities by increasing a parameter $k$ which represents the current
vehicle capacity from 2 to $K$ (line 4). For each $k$, the procedure
first enumerates all $k$-combinations of trips from $\mathcal{T}^+$
(line 5). Let $\mathcal{Q}_k$ represent the set of all $k$-trip
combinations. It then enumerates all valid mini routes for each trip
combination $q\in\mathcal{Q}_k$. Let $\mathrm{\Omega}^v_q$ be this set
of routes for a trip combination $q$. The procedure checks the
feasibility of each route in $\mathrm{\Omega}^v_q$ (using the
$feasible$ function) and adds the ones that are feasible to
$\mathrm{\Omega}^+$ (lines 8--10).

\begin{algorithm}[!t]
	\caption{Mini Route-Enumeration Algorithm for $\mathrm{\Omega}^+$}\label{alg:MREA}
	\begin{algorithmic}[1]
		\Require $\mathcal{T}^+,K$
		\State $\mathrm{\Omega}^+ \leftarrow \text{\O}$
		\For{each $t_c^+\in\mathcal{T}^+$}
		\State $\mathrm{\Omega}^+ \leftarrow \mathrm{\Omega}^+ \cup \{o_c^+\rightarrow d_c^+\}$
		\EndFor
		\For{$k=2$ to $K$}
		\State $\mathcal{Q}_k \leftarrow \{$all $k$-combinations of $\mathcal{T}^+\}$
		\For{each $q\in \mathcal{Q}_k$}
		\State $\mathrm{\Omega}^v_q \leftarrow \{$all valid mini routes of $q\}$
		\For{each $r^+\in\mathrm{\Omega}^v_q$}
		\If{$feasible(r^+)$}
		\State $\mathrm{\Omega}^+ \leftarrow \mathrm{\Omega}^+ \cup \{r^+\}$
		\EndIf
		\EndFor
		\EndFor
		\EndFor
		\State \textbf{return} $\mathrm{\Omega}^+$
	\end{algorithmic}
\end{algorithm}

The labeling procedure by \cite{gschwind2015} makes it possible to
check feasibility when extending partial mini routes and permits a more
efficient implementation of lines 7--10. The set of feasible mini
routes for any trip combination $q$ can be enumerated by performing a
depth-first search which checks the feasibility of each partial route
as it is being extended and backtracks when an extension is
infeasible. Furthermore, the independence of the search procedure for
each trip combination $q\in \mathcal{Q}_k$ allows each combination to
be performed in parallel.

In summary, the enumeration procedure considers all trip combinations
of size $k\leq K$ (of which there are $O(n^K)$ combinations). For each
$k$-combination, it enumerates $(k!)^2$ valid route permutations ($k!$
pickup node permutations followed by $k!$ drop-off node permutations
for each pickup permutation) and checks the feasibility of each. The
procedure therefore has a time and space complexity of $O([K!]^2
n^K)$. \cite{hasan2021} have shown that capacities greater than 5
bring only marginal benefits for the case study, which will also be
confirmed later in this paper.

\subsection{Filtering of Graph $\mathcal{G}$}

Graph $\mathcal{G}$ can be made more compact by only retaining edges
that satisfy a priori route-feasibility constraints. This is done by
pre-processing time-window, pairing, precedence, and ride-duration
limit constraints on $\mathcal{A}$ to identify and eliminate edges
that are infeasible, i.e., those that cannot belong to any feasible AV
route. In this work, the set of infeasible edges is identified using a
combination of rules proposed by \cite{dumas1991} and
\cite{cordeau2006}. These rules are presented in the Appendix.

\section{Valid Inequalities for the CTSPAV} \label{sec:validinequalities}

The CTSPAV MIP is solved with a traditional branch-and-cut procedure that
expoits a number of valid inequalities for the MIP formulation. The
inequalities are valid for all nodes in the branch and bound tree, and
the LP relaxation at each node incorporates all inequalities
discovered up to that point. Numerous families of valid inequalities,
that have been proposed for the TSP
\citep{dantzig1954,grotshcel1985,padberg1991}, ATSP
\citep{fischetti1997}, PCATSP \citep{balas1995}, PDP
\citep{ruland1997}, ATSPTW \citep{ascheuer2000,ascheuer2001}, VRPTW
\citep{kohl1999,bard2002,kallehauge2007}, PDPTW \citep{ropke2009}, and
DARP \citep{cordeau2006}, are also valid for the CTSPAV as the
CTSPAV is a generalization of the DARP. However, this work only
considers inequalities that specifically improve the lower bound on
the vehicle count (the primary objective). This is because extensive
computational experiments from an earlier work \citep{hasan2021}
showed that the LP relaxation already provides a sufficiently strong
lower bound for the secondary objective (total distance).  This
section describes the considered valid inequalities with their
respective separation heuristics when applicable. The following
notation is used to simplify the exposition. For any set of edges
$\mathcal{A}' \subseteq \mathcal{A}$, let $Y(\mathcal{A}') =
\sum_{e\in\mathcal{A}'} Y_e$. For a set of nodes $S \subseteq
\mathcal{N}$, let $\bar{S}$ denote its complement, i.e., $\bar{S} =
\{i\in\mathcal{N}\,|\,i\notin S\}$. For any two node sets
$S,T\subseteq\mathcal{N}$, let $(S,T) = \{(i,j)\in\mathcal{A}\,|\,i\in
S, j\in T\}$. For brevity, $Y(S,T)$ is used to represent
$Y((S,T))$. Finally, for node set
$S\subseteq\mathcal{P}\cup\mathcal{D}$, let $\pi(S) =
\{i\in\mathcal{P}\,|\,n+i\in S\}$ and $\sigma(S) =
\{n+i\in\mathcal{D}\,|\,i\in S\}$ denote the sets of predecessors and
successors of $S$ respectively.

\subsection{Rounded Vehicle-Count Inequalities}
\label{section:darplb}

Suppose that a (fractional) lower bound $\chi_\text{LB}$ is known for the
vehicle count. The inequality
\begin{equation} \label{eqn:rounded_vc_inequality}
	Y(\delta^+(v_s)) \geq \ceil{\chi_\text{LB}}
\end{equation}
is a direct consequence of the integrality of the vehicle count. Such
a lower bound can be obtained by selecting the best bound in the
branch-and-bound algorithm. Let $Y_e^*$ denote the value of $Y_e$ in
the LP-relaxation for this best bound. The lower bound
$\chi_\text{BB}$ can be obtained by
\begin{equation}
	\chi_\text{BB} = \sum_{e\in\delta^+(v_s)} Y_e^*	
\end{equation}
and used in place of $\chi_\text{LB}$ in \eqref{eqn:rounded_vc_inequality}.

\subsection{The Column-Generation Procedure for Deriving Vehicle-Count Lower Bounds}

A stronger lower bound may be obtained from a column-generation
procedure that solves the CTSPAV as a DARP.  This recognition is based
on an earlier work \citep{hasan2021} which discovered that a
column-generation procedure which resembles that used by
\cite{gschwind2015} for solving the DARP is capable of producing
strong lower bounds for the vehicle count of the CTSPAV when it is
paired with an appropriate objective function. This work leverages the
procedure to strengthen the vehicle-count lower bound of the CTSPAV
MIP.

The DARP column-generation procedure of \cite{hasan2021} features a
Pricing Subproblem (PSP) that searches for AV routes with negative
reduced costs to improve the objective function of a set-covering
master problem (MP) whose columns consist of the routes. More
specifically, it utilizes a restricted master problem (RMP) which is
the linear relaxation of the MP that is defined on a subset
$\mathcal{R}' \subseteq \mathcal{R}$ of all feasible AV routes. The
discovered routes are progressively added to $\mathcal{R}'$ as the RMP
and the PSP are solved iteratively. The column generation terminates
when the PSP cannot produce AV routes with negative reduced costs. At
this stage, the objective value $z_\text{RMP}$ of the RMP is identical
to the optimal objective $z^*$ of the linear relaxation of the
original MP. In this work, the column-generation procedure is not used
to obtain a solution to the CTSPAV per se; instead it is used to
extract (potentially strong) lower bounds to the primary objective of
the CTSPAV. The following describes the procedure for obtaining these
lower bounds.

\paragraph{The Restricted Master Problem}

The RMP is a set-covering formulation:
\begin{flalign}
	& \min\,z = \sum_{\rho\in\mathcal{R}'} X_\rho \label{eqn:darp_obj} \\
	& \text{subject to} \nonumber \\
	& \sum_{\rho\in\mathcal{R}'} a_{i,\rho} X_\rho \geq 1 \qquad \forall i\in\mathcal{P} \label{eqn:darp_route_cover} \\
	& X_\rho \geq 0 \qquad \forall \rho\in\mathcal{R}' \label{eqn:darp_route_var} 
\end{flalign}
It is defined on a subset $\mathcal{R}' \subseteq \mathcal{R}$ of all
feasible AV routes, and uses a variable $X_\rho$ to indicate whether
AV route $\rho\in\mathcal{R}'$ is used in the optimal solution. Its
objective function \eqref{eqn:darp_obj} minimizes the number $z$ of
selected AV routes and is therefore identical to the primary
objective of the CTSPAV. Constraints \eqref{eqn:darp_route_cover}
ensure each pickup node is covered in the solution, and $a_{i,\rho}$
is a constant that indicates the number of times node $i$ is visited
by route $\rho$.

\paragraph{The Pricing Subproblem}
The PSP searches for AV routes with negative reduced costs to be added
to $\mathcal{R}'$. It uses $\{\mu_i : i\in\mathcal{P}\}$, the set of
optimal duals of constraints \eqref{eqn:darp_route_cover}, to
compute the reduced costs of the undiscovered routes. The reduced
cost of a route $\rho$ is given by
\begin{equation} \label{eqn:darp_reduced_cost_1}
\bar{c}_\rho = 1 - \sum_{i\in\mathcal{P}} a_{i,\rho} \mu_i.
\end{equation}
To find these routes, a graph $\mathcal{G}$ identical to that defined
in Section \ref{sec:CTSPAV} is first constructed. A reduced cost
$\bar{c}_{(i,j)}$ is then associated with each edge
$(i,j)\in\mathcal{A}$, and it is defined as follows so that the total
cost of any path in $\mathcal{G}$ from $v_s$ to $v_t$ is equivalent to
\eqref{eqn:darp_reduced_cost_1}:
\begin{equation} \label{eqn:darp_edge_reduced_cost_1}
	\bar{c}_{(i,j)} = 
	\begin{cases}
		1&\qquad \forall (i,j)\in\delta^+(v_s)\\
		-\mu_i&\qquad \forall i\in\mathcal{P}, \forall j\in\mathcal{N}\\
		0&\qquad \forall i\in\mathcal{D}, \forall j\in\mathcal{N}.
	\end{cases}
\end{equation}
Obtaining a solution to the PSP is then a matter of finding a feasible
AV route, i.e., a path from $v_s$ to $v_t$ that satisfies the
time-window, capacity, pairing, precedence, and ride-duration limit
constraints, with negative reduced cost. The PSP can be solved by
first finding the least-cost feasible path from $v_s$ to $v_t$ and
then adding it to $\mathcal{R}'$ if the cost is negative. This
approach makes the problem an ESPPRC which can be solved by the
label-setting dynamic program proposed by \cite{gschwind2015}. The
necessity of ensuring elementarity of the path (to ensure its
feasibility), however, makes the problem especially hard to solve
\citep{dror1994}. Since we are only interested in deriving lower
bounds to the vehicle count from this procedure and not in discovering
AV routes per se, the elementarity requirement can be relaxed to admit
a pseudo-polynomial solution from the label-setting algorithm. While
the relaxation, in theory, may cause $z_\text{RMP}$ to converge to a
weaker primal bound as the PSP admits a larger set of routes
$\mathcal{R}'' \supseteq \mathcal{R}'$, other works that have adopted
a similar strategy (e.g., \cite{ropke2009} and \cite{gschwind2015})
have discovered that the lower bound is only slightly weaker in
practice.

\paragraph{Extracting a Lower Bound to the Vehicle Count from the PSP}
As mentioned earlier, $z_\text{RMP}$ converges to $z^*$ and therefore
becomes a valid lower bound to the vehicle count of the CTSPAV when
the PSP is unable to discover a new AV route with negative reduced
cost. However, reaching this point in the procedure typically requires
many column-generation iterations and thus a long computation
time. Prior to it, $z_\text{RMP}$ only represents an upper bound to
$z^*$ and therefore it cannot be used to bound the vehicle
count. Fortunately, the identical unit cost of each AV route in the
RMP allows for the derivation of a lower bound to $z^*$ using the
method proposed by \cite{farley1990}. The Farley bound after the
$k$\textsuperscript{th} column-generation iteration is given by:
\begin{equation}
	z_\text{Farley}^k = \frac{z_\text{RMP}}{1 - \bar{c}^k_\rho}
\end{equation}
where $\bar{c}^k_\rho$ represents the smallest route reduced cost discovered by the PSP after the $k$\textsuperscript{th} iteration. As the value of $z_\text{Farley}^k$ tends to fluctuate between iterations, a monotonically non-decreasing lower bound to $z^*$ can be obtained with the following equation:
\begin{equation} \label{eqn:monotone}
	z_\text{Farley}^k = \max \left\{\frac{z_\text{RMP}}{1 - \bar{c}^k_\rho}, z_\text{Farley}^{k-1}\right\}
\end{equation}

As $z_\text{Farley}^k$ is a lower bound to $z^*$, it is also a valid lower bound to the vehicle count of the CTSPAV. Therefore, $\chi_\text{LB}$ for cut \eqref{eqn:rounded_vc_inequality} may be defined as follows:
\begin{equation}
	\chi_\text{LB} = \max \left\{\chi_\text{BB}, z_\text{Farley}^k\right\}
\end{equation}
Since $z_\text{Farley}^k$ as defined in \eqref{eqn:monotone} is monotonically non-decreasing and improves with the number of column-generation iterations, it is practical to dedicate a single thread for executing this column-generation procedure and use the remaining thread(s) for solving the CTSPAV MIP in parallel. The CTSPAV MIP may then check for the most up-to-date value of $z_\text{Farley}^k$ from the column-generation thread after evaluating the LP relaxation of each tree node and introduce cut \eqref{eqn:rounded_vc_inequality} when there is an improvement to the rounded lower bound.

\subsection{Two-Path Inequalities}

The two-path inequality was originally conceived by \cite{kohl1999}
for the VRPTW. It has been shown to be particularly effective at
strengthening the lower bound for the vehicle count of the VRPTW
\citep{bard2002} and the PDPTW \citep{ropke2009} when vehicle-count
minimization is (part of) the objective function. For a set of nodes
$S \subseteq \mathcal{P}\cup\mathcal{D}$, let $\kappa(S)$ denote the
\emph{minimum} number of vehicles needed to serve $S$, i.e., the
minimum number of vehicles needed to serve all nodes in $S$ while
satisfying all route-feasibility constraints. The following two-path
inequality,
\begin{equation} \label{eqn:2pathcut}
	Y(S,\bar{S}) \geq 2 \qquad \forall S \subseteq \mathcal{P}\cup\mathcal{D}, \kappa(S) > 1
\end{equation}
is valid when it is known that a single vehicle cannot feasibly serve a set $S$, i.e., when $\kappa(S) > 1$. Inequality \eqref{eqn:2pathcut} has a form that is similar to the cutset inequality:
\begin{equation} \label{eqn:cutset}
	Y(S,\bar{S}) \geq 1 \qquad \forall S \subseteq \mathcal{P}\cup\mathcal{D}, |S| \geq 2
\end{equation}
which, in turn, is equivalent to the Dantzig-Fulkerson-Johnson (DFJ) subtour-elimination constraint (SEC) \citep{dantzig1954}:
\begin{equation} \label{eqn:dfj}
	Y(S,S) \leq |S| - 1 \qquad \forall S \subseteq \mathcal{P}\cup\mathcal{D}, |S| \geq 2
\end{equation}
Inequalities \eqref{eqn:cutset} or \eqref{eqn:dfj} are typically used to eliminate subtours in the TSP and ATSP (the subtours manifest themselves as cycles when two separate nodes are used to represent the depot, as is done in this work). For instance, \eqref{eqn:cutset} does so by requiring at least a unit of flow emanating from any set $S$ with two or more nodes. The two-path inequality \eqref{eqn:2pathcut} can therefore be seen as a strengthened SEC as it requires at least two units of flow emanating from any set $S$, however its validity also requires a stronger condition, i.e., $\kappa(S) > 1$. 

\begin{figure}[!t]
	\centering
	\includegraphics[width=0.4\linewidth]{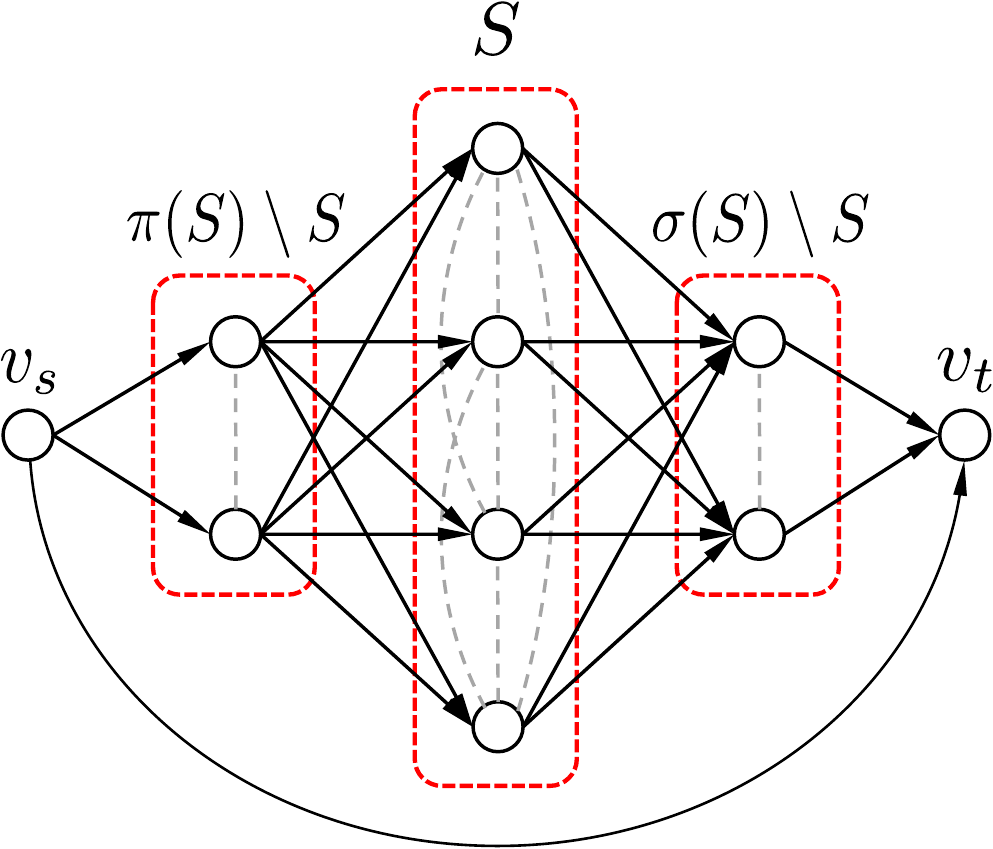}
	\caption{Graph $\mathcal{G}_S$ (Each Dotted Line Represents a Pair of Bidirectional Edges).}
	\label{fig:graph_2pathcut}
\end{figure}

For the CTSPAV, a method similar to that proposed in \cite{ropke2009}
may be used to determine if $\kappa(S) > 1$ for any given set $S$. It
essentially requires one to determine if there exists a
\emph{feasible} path that first visits all the nodes in
$\pi(S)\setminus S$, followed by all the nodes in $S$, and then all
the nodes in $\sigma(S)\setminus S$. If the path does not exist, then
$\kappa(S) > 1$. The task of determining the existence of this path
can be accomplished by first constructing a three-layered graph
$\mathcal{G}_S = (\mathcal{N}_S, \mathcal{A}_S)$ with nodes
$\mathcal{N}_S = \pi(S) \cup S \cup \sigma(S) \cup \{v_s, v_t\}$ and
an initially empty edge set $\mathcal{A}_S$. The nodes from
$\mathcal{N}_S \setminus \{v_s, v_t\}$ are grouped into three layers,
the first consisting of $\pi(S) \setminus S$, the second consisting of
$S$, and the third containing $\sigma(S) \setminus S$. The following
sets of edges are then introduced into $\mathcal{A}_S$:
\begin{itemize}
	\item $\{(v_s,v_t)\}$
	\item $(\{v_s\}, \pi(S)\setminus S) \cap \mathcal{A}$
	\item $(\pi(S)\setminus S, S) \cap \mathcal{A}$
	\item $(S, \sigma(S)\setminus S) \cap \mathcal{A}$
	\item $(\sigma(S)\setminus S, \{v_t\}) \cap \mathcal{A}$
	\item $(\pi(S)\setminus S, \pi(S)\setminus S) \cap \mathcal{A}$
	\item $(S,S) \cap \mathcal{A}$
	\item $(\sigma(S)\setminus S, \sigma(S)\setminus S) \cap \mathcal{A}$
\end{itemize}
where $\mathcal{A}$ denotes the set of feasible edges of graph $\mathcal{G}$ (after they have been filtered). Figure \ref{fig:graph_2pathcut} provides a sketch of $\mathcal{G}_S$. The following sets of edges are introduced into $\mathcal{A}_S$ should either $\pi(S)\setminus S$ or $\sigma(S)\setminus S$ be empty:
\begin{itemize}
	\item If $\pi(S)\setminus S = $ \O, introduce $(\{v_s\}, S) \cap \mathcal{A}$
	\item If $\sigma(S)\setminus S = $ \O, introduce $(S, \{v_t\}) \cap \mathcal{A}$
\end{itemize}

One now needs to determine if there exists a feasible path from $v_s$
to $v_t$ that visits every node of $\mathcal{G}_S$. This problem can
be treated as an ESPPRC, whereby an edge cost of $-1$ is first
assigned to all edges leaving the pickup nodes of $\mathcal{N}_S$
(i.e., edges in $(\pi(S), \mathcal{N}_S \setminus \pi(S))$). A
feasible path from $v_s$ to $v_t$ that visits every node of
$\mathcal{G}_S$ then exists if and only if the least-cost
\emph{elementary} path from $v_s$ to $v_t$ has a total cost of
$-|\pi(S)|$. While this ESPPRC is well-known to be NP-hard
\citep{dror1994}, it can be solved efficiently using the label-setting
algorithm by \cite{gschwind2015} for small $S$. Therefore, one just
needs to solve the ESPPRC and check the total cost of the resulting
elementary path. Should it be greater than $-|\pi(S)|$, then the nodes
of $S$ cannot be feasibly served by a single vehicle, $\kappa(S) > 1$,
and \eqref{eqn:2pathcut} becomes a valid inequality.

\subsubsection{Separation Heuristic} \label{sec:2path_separation}

The separation heuristic for the two-path inequalities first
identifies sets of nodes $S$ for which $\kappa(S) > 1$. As the
two-path inequality is essentially a strengthened SEC, the heuristic
utilized in this work first identifies sets of nodes that form
subtours (cycles) in the LP-relaxation solution at each tree node. Let
$Y_e^*$ denote the value of $Y_e$ from the solution of the LP
relaxation. For every subtour $S$ considered, the heuristic then
checks if $\sum_{e\in(S, \bar{S})} Y_e^* < 2$ and then if $\kappa(S) >
1$. Satisfaction of these two conditions indicates that the two-path
inequality is valid for $S$, and that it is violated by $S$ in the
LP-relaxation solution. The heuristic therefore adds the two-path
inequality to eliminate generation of the subtour from subsequent LP
solutions.

To identify subtours from the LP relaxation at each tree node, the
heuristic by \cite{drexl2013} is used. The heuristic was proposed as a
cheaper yet effective alternative for identifying violated SECs to the
exact method proposed by \cite{gomory1961}, as it has an $O(n^2)$
complexity compared to the $O(n^4)$ complexity of the latter. For any
LP solution, a support graph, $\mathcal{G}_\text{sp} =
(\mathcal{N}_\text{sp}, \mathcal{A}_\text{sp})$, is first constructed
with nodes $\mathcal{N}_\text{sp} = \mathcal{N}$ and edges
$\mathcal{A}_\text{sp} = \{e\in\mathcal{A}\,|\,Y_e^* > 0\}$. All
strongly-connected components (SCCs) of $\mathcal{G}_\text{sp}$ are
then identified, where an SCC of a graph is its subgraph with more
than one node whereby there exists a path between all pairs of its
nodes. The rationale behind identification of SCCs is that each forms
a subtour (the nodes of the SCC form a cycle(s) as every node is
reachable from another). In practice, all SCCs of
$\mathcal{G}_\text{sp}$ can be computed using the algorithm by
\cite{tarjan1972} which has a time complexity of
$O(|\mathcal{N}_\text{sp}| + |\mathcal{A}_\text{sp}|)$. Let
$\mathcal{S}_\text{sp}$ denote the set of all SCCs of
$\mathcal{G}_\text{sp}$, and for each SCC $c\in\mathcal{S}_\text{sp}$,
let $S_c$ denote its set of nodes. For every
$c\in\mathcal{S}_\text{sp}$, the heuristic then checks if the total
flow leaving $S_c$ is less than 2, i.e., if $\sum_{e\in(S_c,
  \bar{S}_c)} Y_e^* < 2$. If this condition is satisfied for $S_c$,
the heuristic then determines if $\kappa(S_c) > 1$ using the procedure
described earlier. Finally, the two-path cut $Y(S_c, \bar{S}_c) \geq
2$ is introduced to the MIP if $\kappa(S_c) > 1$.

Due to the expensive nature of the procedure for determining if
$\kappa(S_c) > 1$, results of the procedure for every set $S_c$ are
stored in a hash table, and the hash table is examined first before
the procedure is performed on any set $S$ to ensure that the same
calculations are not repeated. Furthermore, the part of the procedure
which solves an ESPPRC on graph $\mathcal{G}_{S}$ can also be made
more efficient. Instead of directly applying the label-setting
algorithm of \cite{gschwind2015} which proposes keeping track of all
visited pickup nodes and preventing path extensions to the already
visited nodes to ensure elementarity, the procedure proposed by
\cite{boland2006} can be used. The latter entails iteratively solving
a sequence of relaxed SPPRCs, whereby the elementarity requirement is
completely relaxed in the very beginning. A repeated node from the
solution of the relaxed problem is selected and added to a set
$\mathcal{U}$, after which the problem is solved again, this time with
an additional restriction that the nodes in $\mathcal{U}$ can only be
visited once. The procedure is repeated with $\mathcal{U}$ being
progressively enlarged until an elementary path is discovered. The
rationale behind this procedure is that solving the sequence of
relaxed SPPRCs is usually less expensive than solving a single ESPPRC
in practice, as often times the former discovers an elementary path
without having to include all pickup nodes in the set
$\mathcal{U}$. \cite{desaulniers2008} proposed adding only the first
repeated node from the solution of the relaxed problem to
$\mathcal{U}$ after each iteration, and our initial evaluations show
that this approach works very well in practice.

\subsection{Predecessor and Successor Inequalities}
Predecessor and successor inequalities were first introduced by \cite{balas1995} for the PCATSP. The predecessor inequality ($\pi$-inequality) is given by:
\begin{equation} \label{eqn:pred_inequality}
	Y(S\setminus \pi(S), \bar{S}\setminus \pi(S)) \geq 1  \qquad \forall S \subseteq \mathcal{P}\cup\mathcal{D}, |S| \geq 2
\end{equation}
and the successor inequality ($\sigma$-inequality) is given by:
\begin{equation} \label{eqn:succ_inequality}
	Y(\bar{S}\setminus \sigma(S), S\setminus \sigma(S)) \geq 1  \qquad \forall S \subseteq \mathcal{P}\cup\mathcal{D}, |S| \geq 2
\end{equation}
These inequalities are essentially lifted versions of the cutset inequality \eqref{eqn:cutset}. They are also valid for the CTSPAV as it generalizes the PCATSP.

\subsubsection{Separation Heuristic}

The heuristic utilized to separate $\pi$- and $\sigma$-inequalities is
very similar to that described in Section \ref{sec:2path_separation}
for the two-path inequality. At each tree node, values of $Y_e^*$ are
first used to construct a support graph $\mathcal{G}_\text{sp}$, after
which $\mathcal{S}_\text{sp}$ which represents the set of all SCCs of
$\mathcal{G}_\text{sp}$ are identified. For each
$c\in\mathcal{S}_\text{sp}$, the heuristic then checks if either
inequalities \eqref{eqn:pred_inequality} or
\eqref{eqn:succ_inequality} have been violated for $S_c$, i.e., if
either $Y(S_c\setminus \pi(S_c), \bar{S}_c\setminus \pi(S_c)) < 1$ or
$Y(\bar{S}_c\setminus \sigma(S_c), S_c\setminus \sigma(S_c)) <
1$. Finally, corresponding $\pi$- or $\sigma$-inequalities are
introduced to the MIP for each violation.

\subsection{Lifted MTZ Inequalities}

The lifted MTZ inequality was initially proposed by
\cite{desrochers1991} for the VRPTW. They were intended to strengthen
MTZ constraints that are similar to \eqref{eqn:ctspav_travel_time_1}
and \eqref{eqn:ctspav_travel_time_2} which are well-known to produce
weak LP relaxations \citep{langevin1990,gouveia1999}. The MTZ
constraints for an edge $(i,j)$ is strengthened by taking into
consideration the flow along the opposite edge $(j,i)$ combined with
the fact that only one of the edges may have positive flow in a
feasible integer solution. The lifted versions constraints
\eqref{eqn:ctspav_travel_time_1} and \eqref{eqn:ctspav_travel_time_2}
are given by \eqref{eqnliftedMTZ1} and \eqref{eqnliftedMTZ2}
respectively.
\begin{flalign} 
	& T_i + s_i + \tau_{(i,j)} \leq T_j + M_{(i,j)}(1 - Y_{(i,j)}) - \alpha_{(j,i)} Y_{(j,i)} \qquad \forall i,j\in\mathcal{P}\cup\mathcal{D} \label{eqnliftedMTZ1} \\ 
	& T_i + s_i + \tau_{(i,j)} \geq T_j - \bar{M}_{(i,j)}(1 - Y_{(i,j)}) - \beta_{(j,i)} Y_{(j,i)}\qquad \forall i\in\mathcal{P}\cup\mathcal{D}, \forall j\in\mathcal{D} \label{eqnliftedMTZ2}
\end{flalign}

To correctly lift the constraints using this technique, the coefficients of the flow variable of the opposite edge, $\alpha_{(j,i)}$ and $\beta_{(j,i)}$, are assigned values that are as large as possible while ensuring that inequalities \eqref{eqnliftedMTZ1} and \eqref{eqnliftedMTZ2} are still valid for any feasible integer solution. \cite{desrochers1991} proposed coefficient values for the VRPTW that ensure the earliest start of service times for every node. As serving pickup nodes as early as possible may not be desirable for the CTSPAV (as doing so lengthens the ride duration of the picked-up rider and thus increases the likelihood of exceeding her ride-duration limit), the coefficients are adjusted to \eqref{eqn:alpha} and \eqref{eqn:beta} for the CTSPAV.
\begin{flalign} 
	& \alpha_{(j,i)} = 
	\begin{cases}
		M_{(i,j)} - s_i - \tau_{(i,j)} - s_j - \tau_{(j,i)}&\qquad \text{if } i \in\mathcal{D}\\
		M_{(i,j)} - s_i - \tau_{(i,j)} - b_i + a_j&\qquad \text{otherwise}
	\end{cases} \label{eqn:alpha} \\
	& \beta_{(j,i)} = -\bar{M}_{(i,j)} - s_i - \tau_{(i,j)} - s_j - \tau_{(j,i)} \label{eqn:beta}
\end{flalign}

\noindent
The validity of the lifted constraints can be verified by first
substituting \eqref{eqn:alpha} and \eqref{eqn:beta} into
\eqref{eqnliftedMTZ1} and \eqref{eqnliftedMTZ2} respectively, and then
setting the flows along edges $(i,j)$ and $(j,i)$ to zero or
setting the flow along either edge to one. Firstly, setting both
$Y_{(i,j)}$ and $Y_{(j,i)}$ to zero just disables constraints
\eqref{eqnliftedMTZ1} and \eqref{eqnliftedMTZ2} for both edges. Next,
setting $Y_{(i,j)} = 1$ and $Y_{(j,i)} = 0$ produces the following
constraints,
\begin{flalign} 
	& T_i + s_i + \tau_{(i,j)} \leq T_j \qquad \text{if }i,j\in\mathcal{P}\cup\mathcal{D} \\ 
	& T_i + s_i + \tau_{(i,j)} \geq T_j \qquad \text{if }i\in\mathcal{P}\cup\mathcal{D},j\in\mathcal{D}
\end{flalign}
which simply enforce the increasing service time requirement along
edge $(i,j)$. Finally, setting $Y_{(i,j)} = 0$ and $Y_{(j,i)} = 1$
results in the following set of constraints:
\begin{flalign} 
	& T_j + s_j + \tau_{(j,i)} \geq T_i \qquad \text{if }j\in\mathcal{P}\cup\mathcal{D},i\in\mathcal{D} \label{eqn1} \\ 
	& T_i - T_j \leq b_i - a_j \qquad \text{if }j\in\mathcal{P}\cup\mathcal{D},i\in\mathcal{P} \label{eqn2} \\ 
	& T_j + s_j + \tau_{(j,i)} \leq T_i \qquad \text{if }j\in\mathcal{D},i\in\mathcal{P}\cup\mathcal{D} \label{eqn3}
\end{flalign}
Constraints \eqref{eqn1} and \eqref{eqn3} simply enforce increasing
service times along edge $(j,i)$, while \eqref{eqn2} is obviously a
valid inequality if edge $(j,i)$ is selected.

\subsection{Lifted Time-Bound Inequalities}

The lifted time-bound inequalities were also proposed by
\cite{desrochers1991} to strengthen the time-window constraints of the
VRPTW. Inequalities \eqref{eqn:liftedlowerbound} and
\eqref{eqn:liftedupperbound} strengthen the time-window constraints of
node $i$ by taking into consideration the temporal requirements along
the node's incoming and outgoing edges with positive flow.
\begin{flalign} 
	& T_i \geq a_i + \sum_{(j,i)\in\delta^-(i)} \max \{0,a_j-a_i+s_j+\tau_{(j,i)}\} Y_{(j,i)} \qquad \forall i\in\mathcal{P}\cup\mathcal{D} \label{eqn:liftedlowerbound} \\
	& T_i \leq b_i - \sum_{(i,j)\in\delta^+(i)} \max \{0,b_i-b_j+s_i+\tau_{(i,j)}\} Y_{(i,j)} \qquad \forall i\in\mathcal{P}\cup\mathcal{D} \label{eqn:liftedupperbound}
\end{flalign}

\section{Computational Results} \label{sec:results}

This section presents the computational results of the branch-and-cut algorithm on problem instances derived from a real-world dataset of commute trips. 

\subsection{Algorithmic Settings}

Three variants of the branch-and-cut algorithm are considered and
contrasted in the evaluations; they are named
CTSPAV\textsubscript{Base}, CTSPAV\textsubscript{SEC}, and
CTSPAV\textsubscript{Hybrid}. Each is differentiated by the types of
valid inequalities included in its implementation. They are specificied as follows:
\begin{itemize}
\item CTSPAV\textsubscript{Base} is the core algorithm and implements
the simplest valid inequalities: lifted time bounds, lifted
MTZ, and rounded vehicle count which uses $\chi_{\text{BB}}$ as its
lower bound;
\item CTSPAV\textsubscript{SEC} is CTSPAV\textsubscript{Base} with the
  two-path, predecessor, and successor inequalities;
\item CTSPAV\textsubscript{Hybrid} is CTSPAV\textsubscript{Base} with the DARP lower
  bound from Section \ref{section:darplb}.
\end{itemize}
The latter variant also uses the interior-point, dual-stabilization
method proposed by \cite{rousseau2007} to accelerate the convergence of
its column-generation procedure. Furthermore, instead of only
selecting the least-cost feasible path with negative reduced cost in
its PSP, all non-dominated paths resulting from the label-setting
algorithm with negative reduced costs are added to $\mathcal{R}'$ to
further accelerate convergence.

\subsection{Construction of Problem Instances}

Problem instances for the computational evaluations are derived from
the commute trip dataset first used by \cite{hasan2018}. It consists
of the real-world arrival and departure times to 15 parking structures
located in downtown Ann Arbor, Michigan, of approximately 15,000
commuters that were collected throughout the month of April 2017. This
information, when joined with the home addresses of every commuter,
allowed the reconstruction of their daily commute trips. The
performance evaluations utilize the trips made by commuters living
within Ann Arbor's city limits, the region bounded by highways US-23,
M-14, and I-94. More specifically, the 2,200 commute trips from this
region made on the busiest day of the month (Wednesday of week 2) were
first selected and then partitioned into smaller problem instances
using the clustering algorithm described in Section
\ref{sec:clustering}. Trip sharing is then only considered
intra-cluster with the largest parking structure arbitrarily
designated as the depot for all clusters.

In addition to this, the following assumptions are made in order to
define the time windows and ride-duration limits of each
trip. Consistent with past works on the DARP (e.g.,
\cite{jaw1986,cordeau2003,cordeau2006}), each rider $i$ specifies a
desired arrival time $at_i^+$ at the destination of her inbound trip
and a desired departure time $dt_i^-$ at the origin of her outbound
trip when requesting a trip. Riders also tolerate a maximum shift of
$\pm \mathrm{\Delta}$ to the desired times. By considering the arrival
and departure times to and from the parking structures as the desired
times, an arrival-time upper bound at node $n+i$ of $b_{n+i} = at_i^+
+ \mathrm{\Delta}$ and a time window at node $2n+i$ of $[a_{2n+i},
  b_{2n+i}] = [dt_i^- - \mathrm{\Delta}, dt_i^- + \mathrm{\Delta}]$
are defined for each $i\in\mathcal{P}^+$. Consequently, the time
window at node $i$ is given by $[a_i,b_i] = [b_{n+i} - s_i - L_i -
  2\Delta,b_{n+i} - s_i - L_i]$ and the arrival-time upper bound at
node $3n+i$ is given by $b_{3n+i} = b_{2n+i} + s_{2n+i} + L_{2n+i}$
for each $i\in\mathcal{P}^+$. Finally, consistent with
\cite{hunsaker2002}, the ride-duration limit of each trip is defined
as an $R\%$ extension to the direct trip, i.e., $L_i =
(1+R)\tau_{i,n+i}$ for each $i\in\mathcal{P}$.

A set of tight, medium, and large problem instances are constructed by
varying parameter $N$ in the clustering algorithm together with
$\mathrm{\Delta}$ and $R$. The parameter combinations are
carefully selected to highlight performance differences in the
three variants of the branch-and-cut algorithm considered. A vehicle
capacity of $K = 4$ is used in all instances to represent the use of
autonomous cars. Table \ref{tab:instance_params} shows the parameters
used together with the number of instances created when the clustering
algorithm is applied on the set of 2,200 commuters:

\begin{table}[!t]
	\centering
	\caption{Parameters for Constructing Problem Instances}
	\label{tab:instance_params}
	\begin{tabular}{cccccc}
		\hline\up
		Problem size & $N$   & $\mathrm{\Delta}$   & $R$    & $K$ & Number of instances\down\\
		\hline\up
		Large        & 100 	 & 10 mins 			   & 0.50   & 4   & 22                  \\
		Medium       & 75    & 10 mins             & 0.50   & 4   & 30                  \\
		Tight        & 100 	 & 5 mins              & 0.25   & 4   & 22\down\\
		\hline
	\end{tabular}
\end{table}

\subsection{Experimental Settings}

All algorithms are implemented in C++. Parallelization of the mini
route-enumeration algorithm is handled with OpenMP, while the parallel
execution of the column-generation procedure and the MIP of
CTSPAV\textsubscript{Hybrid} is handled with the thread class from the
C++11 standard library. All LPs and MIPs are solved with Gurobi 9.0.2,
while graph algorithms from the Boost Graph Library (version 1.70.0)
are used to calculate SCCs of a graph and to implement the
label-setting algorithm of \cite{gschwind2015}. Gurobi's callback
feature is used to implement the bespoke cutting-plane separation and
insertion, while the MIP solver is configured with its default
parameters.  For problem instance construction, Geocodio is used to
geocode GPS coordinates of every address considered, after which
GraphHopper's Directions API is used in conjunction with OpenStreetMap
data to estimate the shortest path, travel time, and travel distance
between any two nodes.  Unless stated otherwise, every problem
instance is solved on a compute cluster, each utilizing 4 cores of a
3.0 GHz Intel Xeon Gold 6154 processor and 16 GB of RAM. All four
cores are used for the MREA. For CTSPAV\textsubscript{Hybrid}, one
core is dedicated for the column-generation procedure while the
remaining three are used for solving the MIP. All four cores are used
for solving the MIPs of CTSPAV\textsubscript{SEC} and
CTSPAV\textsubscript{Base}. Finally, a 2-hour time budget is allocated
for solving all MIPs.

\subsection{Algorithm Performance Comparison}

\begin{table}[!t]
	\centering
	\caption{Average Vehicle Count and Optimality Gaps of Every CTSPAV Variant for Every Problem Size}
	\label{tab:my-table7}
	\begin{tabular}{ccccccc}
		\hline\up
		\multirow{2}{*}{\begin{tabular}[c]{@{}c@{}}CTSPAV\\ variant\end{tabular}} & \multicolumn{3}{c}{Average vehicle count gap} & \multicolumn{3}{c}{Average optimality gap\down} \\
		\cline{2-7}\rule{0pt}{12pt}
		& Large & Medium & Tight & Large & Medium & Tight\down\\
		\hline\up
		Hybrid & 1.18 & 0.50 & 0.00 & 31.8\% & 16.6\% & 0.0\% \\
		SEC & 1.73 & 0.73 & 0.09 & 45.5\% & 23.8\% & 1.7\% \\
		Base & 2.50 & 1.67 & 0.14 & 68.0\% & 59.0\% & 3.2\%\down\\
		\hline
	\end{tabular}
\end{table}

Table \ref{tab:my-table7} first summarizes the average vehicle count
gaps and average optimality gaps obtained for every problem size and
every CTSPAV variant. $\chi_\text{MIP}$, $z_\text{MIP}$, and
$z_\text{BB}$ denote the vehicle count, the objective value of the
best incumbent solution, and its best bound respectively. The vehicle
count gap is given by $\chi_\text{MIP} - \ceil{\chi_\text{LB}}$, while
the optimality gap is given by $(z_\text{MIP} -
z_\text{BB})/z_\text{MIP}$. The complete results of all the
computational experiments are listed in Tables
\ref{tab:my-table}--\ref{tab:my-table6} in the Appendix. Note that the
route enumeration times for every problem instance are consistently
less than 60 seconds, which highlights the efficiency of the MREA.

The average optimality gaps for large and medium instances appear to
be relatively large.  However, a closer examination paints a different
picture, as their values are relatively small across the board. In
fact, the average count gap for CTSPAV\textsubscript{Hybrid} is only a
little above one for the large problem instances, and is less than one for
the tighter instances. The values for CTSPAV\textsubscript{Hybrid} are
also consistently smaller across the board than those of
CTSPAV\textsubscript{SEC} which, in turn, are smaller than those of
CTSPAV\textsubscript{Base}. This observation provides the first
evidence of the capability of CTSPAV\textsubscript{Hybrid}'s
column-generation procedure at producing very strong lower bounds for
the primary objective; it also demonstrates the effectiveness of the
combination of the two-path, successor, and predecessor inequalities
at closing the vehicle count gap (compared to an implementation that
only adopts the three basic inequalities). While the latter set of
inequalities produces significant improvements in closing the primary
gap, they are nevertheless outperformed by the rounded vehicle-count
inequalities of CTSPAV\textsubscript{Hybrid}.

\begin{figure}[!t]
	\centering
	\includegraphics[width=1.0\linewidth]{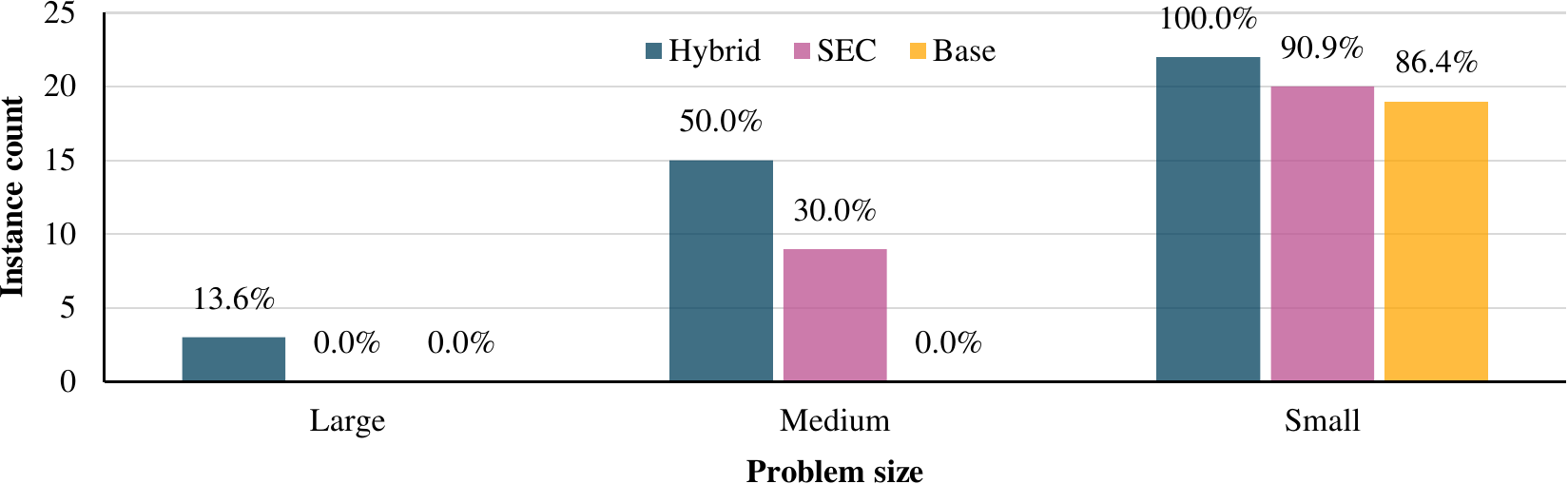}
	\caption{Number of Problem Instances Whereby Vehicle Count Gap is Closed by Every CTSPAV Variant.}
	\label{fig:optimal_count_chart}
\end{figure}

Figure \ref{fig:optimal_count_chart} provides a different perspective
by summarizing the number of problem instances whose vehicle count
gaps are successfully closed within the 2-hour time limit for every
CTSPAV variant. It also displays each count as a fraction of the total
number of instances considered. For the large instances,
CTSPAV\textsubscript{Hybrid} could only close the gap for three
instances, while the other two variants could not for any of the
problems from the set. This number improves for the medium problem
instances, where CTSPAV\textsubscript{Hybrid} could now close the gap
for 15 out of the 30 instances, while CTSPAV\textsubscript{SEC} could
do the same for 9 of the instances. However,
CTSPAV\textsubscript{Base} still cannot close the primary gap for
any. Finally, for the tight problem instances,
CTSPAV\textsubscript{Hybrid} produces the optimal solution for all of
them, while CTSPAV\textsubscript{SEC} closes the primary gap for
90.9\% of the instances and CTSPAV\textsubscript{Base} does the same
for 86.4\% of them. Regardless of the set of problem instances being
considered, the trend is clear: (1) The additional set of inequalities
adopted by CTSPAV\textsubscript{SEC} allows it to successfully close
the primary gap of more instances than CTSPAV\textsubscript{Base}, and
(2) CTSPAV\textsubscript{Hybrid} consistently outperforms the other
two CTSPAV variants at closing the optimality gap.  The latter
observation provides yet another evidence of the efficacy of the
CTSPAV\textsubscript{Hybrid}'s column-generation procedure at
generating strong lower bounds for the primary objective.

\begin{figure}[!t]
	\centering
	\includegraphics[width=1.0\linewidth]{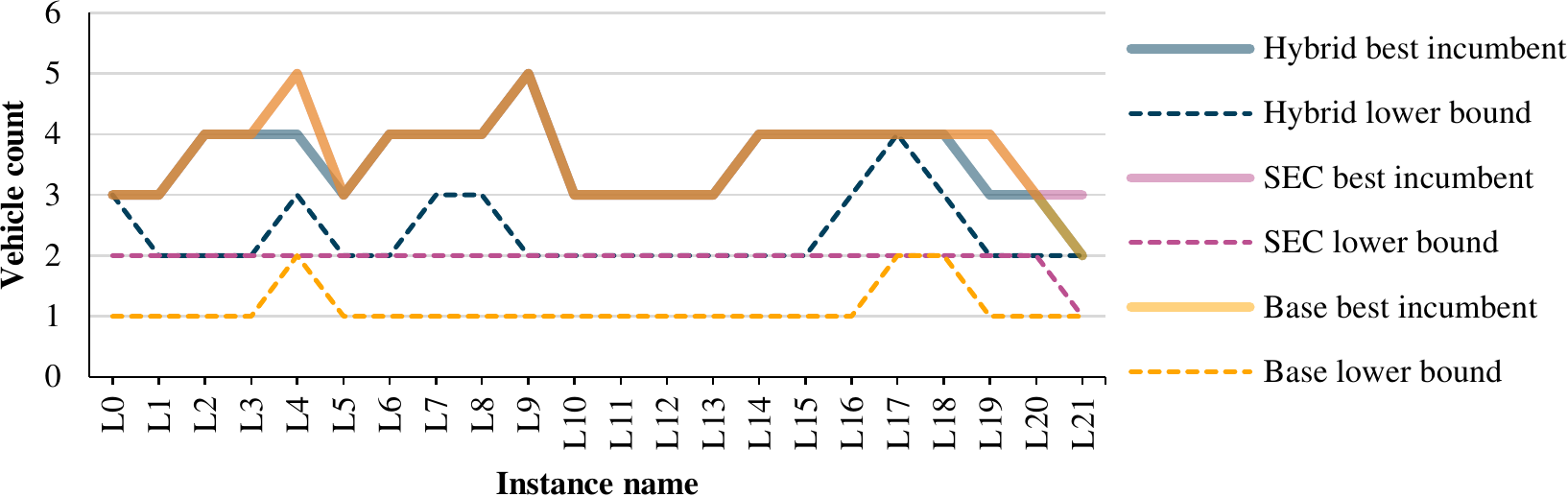}
	\caption{Best Incumbent Solution and Lower Bound for Vehicle Count of Every CTSPAV Variant for Every Large Problem Instance.}
	\label{fig:large_instance_gaps}
\end{figure}

\begin{figure}[!t]
	\centering
	\includegraphics[width=1.0\linewidth]{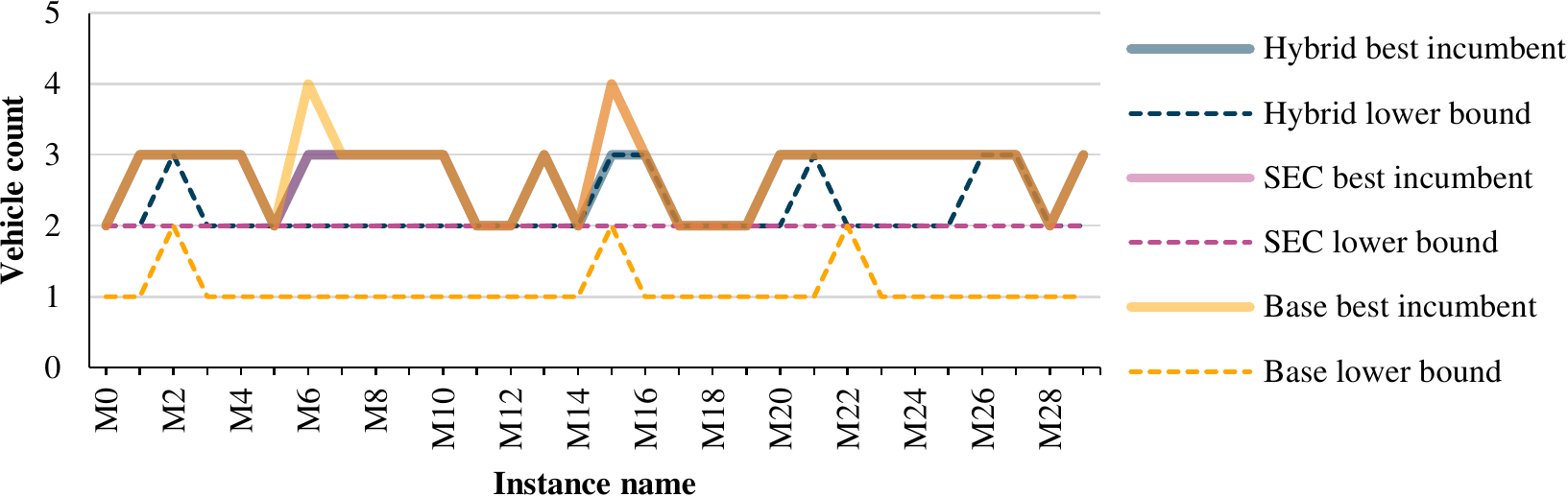}
	\caption{Best Incumbent Solution and Lower Bound for Vehicle Count of Every CTSPAV Variant for Every Medium Problem Instance.}
	\label{fig:medium_instance_gaps}
\end{figure}

\begin{figure}[!t]
	\centering
	\includegraphics[width=1.0\linewidth]{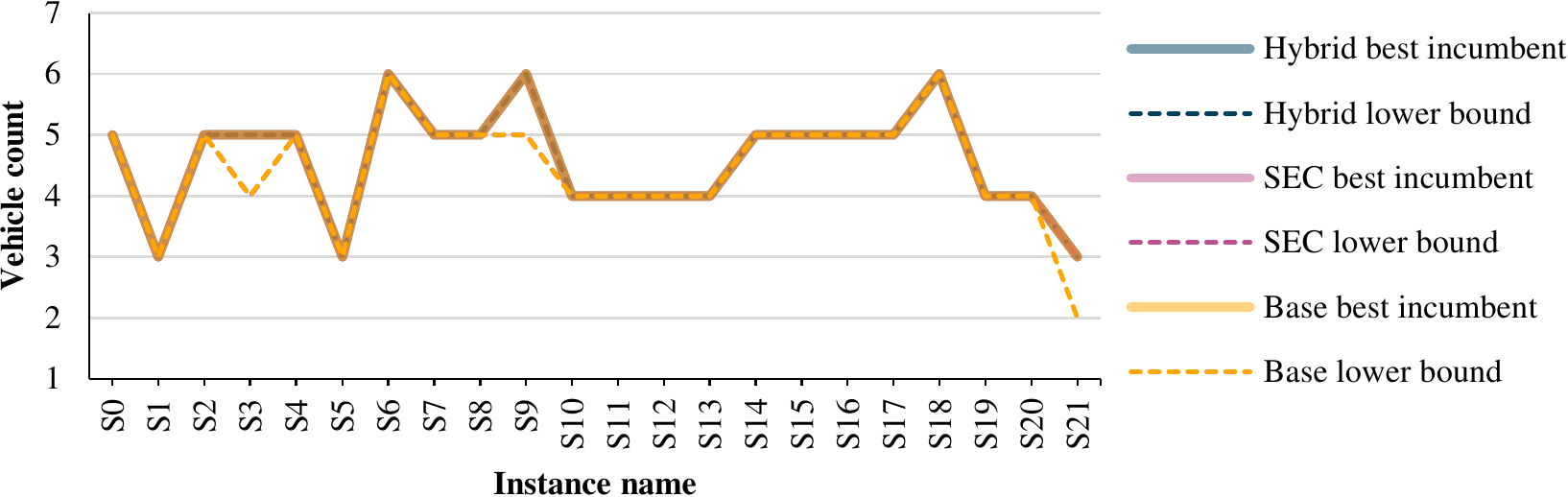}
	\caption{Best Incumbent Solution and Lower Bound for Vehicle Count of Every CTSPAV Variant for Every Tight Problem Instance.}
	\label{fig:small_instance_gaps}
\end{figure}

Instead of aggregating the results from each problem set, Figures
\ref{fig:large_instance_gaps}, \ref{fig:medium_instance_gaps}, and
\ref{fig:small_instance_gaps} provide a closer look at the primary
objective value and its corresponding lower bound for every problem
instance from the large, medium, and tight sets respectively. For
instance, Figure \ref{fig:large_instance_gaps} shows the best
incumbent solution and the lower bound for the vehicle count of every
CTSPAV variant for every large problem instance. The figure reveals
that, except for a few instances, all three variants produced
identical final vehicle counts. The difference, however, lies in their
lower bounds. The lower bounds of CTSPAV\textsubscript{Hybrid}
dominate those of CTSPAV\textsubscript{SEC} in every instance. In
turn, those of the latter dominate the lower bounds of
CTSPAV\textsubscript{Base} in every instance as well. The same
observation is carried over to Figure \ref{fig:medium_instance_gaps}
which summarizes the primary gap of every instance from the medium
set. While CTSPAV\textsubscript{Hybrid} and CTSPAV\textsubscript{SEC}
produce identical lower bounds for more instances from this set, on
the whole, lower bounds of CTSPAV\textsubscript{SEC} are still
dominated by those of CTSPAV\textsubscript{Hybrid}. Similarly, they
both dominate the lower bounds of CTSPAV\textsubscript{Base}. Finally,
Figure \ref{fig:small_instance_gaps} summarizes the results of the
tight instances, and confirms the observations from the previous two
figures. The observations from Figures \ref{fig:large_instance_gaps},
\ref{fig:medium_instance_gaps}, and \ref{fig:small_instance_gaps} lead
to the following conclusion: Regardless of the size of the problem
considered, there is a clear delineation between the strengths of the
lower bounds for the primary objective of the three CTSPAV
variants. CTSPAV\textsubscript{Hybrid} dominates
CTSPAV\textsubscript{SEC} which, in turn, dominates
CTSPAV\textsubscript{Base}. The relative strength of
CTSPAV\textsubscript{Hybrid}'s lower bound directly contributes to its
ability to close or narrow the optimality gap of more problem
instances than the other two variants.

\subsection{Analysis of the Lower Bounds}

\begin{figure}[!t]
	\centering
	\includegraphics[width=1.0\linewidth]{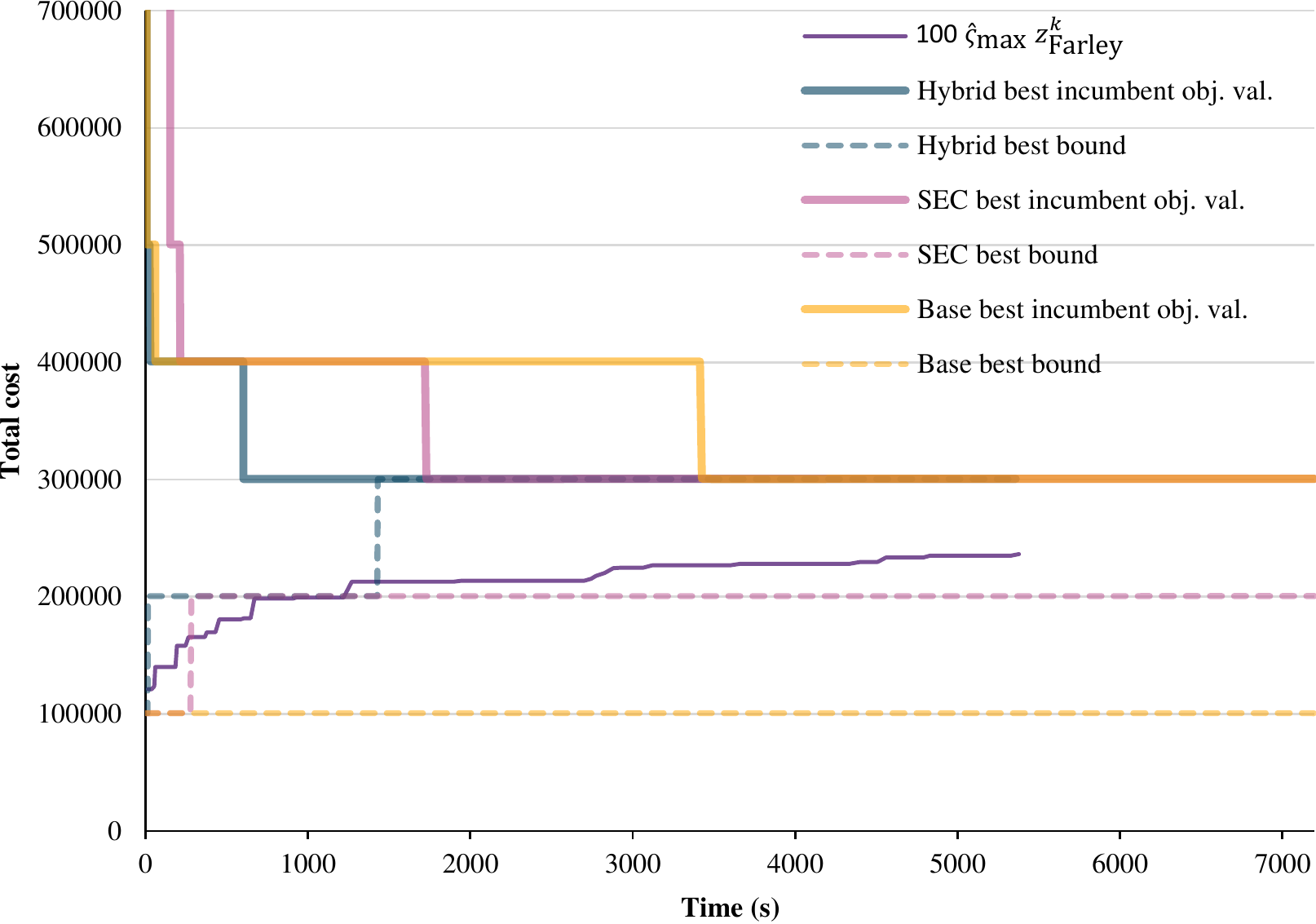}
	\caption{Evolution of Best Incumbent Objective Value and Best Bound of Every CTSP Variant for Problem Instance L0}
	\label{fig:bounds_evolution}
\end{figure}

Figure \ref{fig:bounds_evolution} presents a closer examination of the
evolution of the best bound and best incumbent objective value of
every CTSPAV variant over time for a specific problem instance
(instance L0). It also shows the progression of $z^k_\text{Farley}$
(after it has been scaled by $100 \hat{\varsigma}_\text{max}$) over
time; the lower bound is obtained by rounding it to the smallest
multiple of $100 \hat{\varsigma}_\text{max}$.  Since the MIP solver,
using its default heuristics, is able to discover strong integer
solutions fairly quickly for this formulation, the critical challenge
lies in closing the optimality gap quickly. Unfortunately, the CTSPAV
formulation uses big-$M$ constants in constraints
\eqref{eqn:ctspav_travel_time_1} and \eqref{eqn:ctspav_travel_time_2}
which produce weak LP relaxations.

The lifted MTZ and lifted time-bound inequalities only provide
marginal improvements to the LP relaxation. While the rounded
vehicle-count inequality has the capability of rectifying the issue,
$\chi_\text{BB}$ rarely becomes fractional in practice, and thus the
version of the inequality that only uses $\chi_\text{BB}$ as its lower
bound rarely improves the vehicle-count lower bound. This explains why
CTSPAV\textsubscript{Base} always produces the weakest lower
bounds. Separation heuristics of the two-path, successor, and
predecessor inequalities attempt to alleviate this situation by first
searching for subtours that result from the flow of an LP-relaxation
solution, and then introducing the respective inequalities to remove
these subtour flows from subsequent LP relaxations. The experimental
results of CTSPAV\textsubscript{SEC} demonstrate that these
inequalities are indeed effective at further strengthening the LP
bounds, however the results also show that their effect on the best
bound tends to stagnate over time.

The CTSPAV\textsubscript{Hybrid} attempts to circumvent the CTSPAV
formulation's weak LP bound by dedicating a computational thread to
solving the same problem using a DARP formulation that focuses only on
the primary objective. The Farley bound $z^k_\text{Farley}$ of the
DARP relaxation provides a lower bound, and its scaled values in
Figure \ref{fig:bounds_evolution} show that it progressively improves
over time even after the best bounds of CTSPAV\textsubscript{Base} and
CTSPAV\textsubscript{SEC} begin to stagnate. The ability of the column-generation
to produce relatively stronger lower bounds
can be attributed to a few factors:
\begin{enumerate}
\item The RMP formulation does not utilize any big-$M$ constants.
\item The RMP uses only one set of binary variables ($X_\rho$), as
  opposed to two by the CTSPAV MIP ($X_r$ and $Y_e$). Therefore, fewer
  convex combinations of its routes are allowed in its LP relaxation,
  which leads to stronger primal (and dual) lower bounds.
\item \cite{ropke2006} showed that the set-covering formulation
  actually implies several valid inequalities (precedence and
  strengthened precedence inequalities) that would otherwise need to
  be enforced explicitly in an edge flow formulation.
\end{enumerate}


\noindent
The approach of dedicating a single thread for executing the
column-generation procedure also has a side benefit: it allows the
branch-and-bound algorithm to freely explore more tree nodes without
being encumbered by expensive separation heuristics. This is evident
from a comparison of the number of explored nodes for several problem
instances, for example, those of CTSPAV\textsubscript{Hybrid} and
CTSPAV\textsubscript{SEC} for instances L1, L5, and L12 from Tables
\ref{tab:my-table} and \ref{tab:my-table2}. The results show that the
former was able to explore significantly more nodes, and this could,
in turn, lead to the discovery of better integer solutions. While
CTSPAV\textsubscript{Hybrid} had one fewer thread for solving its MIP,
it also did not have to execute any of the expensive separation
heuristics of CTSPAV\textsubscript{SEC} which consequently resulted in
a net gain in terms of the number of nodes it could explore.

\subsection{Analysis of the Column-Generation Heuristic}

It is useful to contrast these results with the column-generation
heuristic proposed by \cite{hasan2021}.  The heuristic does not
exhaustively enumerate all the mini routes in
$\mathrm{\Omega}$. Instead it uses a column-generation procedure
consisting of a restricted master problem
(RMP\textsubscript{CTSPAV})---the linear relaxation of MIP model
\eqref{eqn:ctspav_obj}--\eqref{eqn:ctspav_edge_var} defined on only a
subset $\mathrm{\Omega}' \subseteq \mathrm{\Omega}$ of the mini
routes--- and a pricing subproblem (PSP\textsubscript{CTSPAV}) that
searches for mini routes with negative reduced costs to augment
$\mathrm{\Omega}'$. The RMP\textsubscript{CTSPAV} and
PSP\textsubscript{CTSPAV} are solved repeatedly until the
PSP\textsubscript{CTSPAV} is unable to find any mini route with
negative reduced cost. Then the heuristic solves the
RMP\textsubscript{CTSPAV} as a MIP (that does not incorporate the
valid inequalities considered in this work) to obtain a feasible
integer solution. Since the heuristic only considers a subset of the
feasible mini routes, \emph{it is incapable of proving the optimality
  of its solution unless the solution of its RMP\textsubscript{CTSPAV}
  at convergence is integral} (which is never the case for the
instances considered). Nevertheless, it is still instructive to
compare its results against those of the exact
CTSPAV\textsubscript{Hybrid} method to gauge the effectiveness of its
column-generation procedure at identifying useful mini routes.

Tables \ref{tab:my-table9}, \ref{tab:my-table10}, and
\ref{tab:my-table11} (in the Appendix) give comprehensive results for
the heuristic on every large, medium, and tight instance
respectively. The results show that significantly fewer columns (mini
routes) are considered by the heuristic. On average, it considers
66\%, 62\%, and 16\% fewer columns for the large, medium, and tight
instances respectively compared to
CTSPAV\textsubscript{Hybrid}. However, the final vehicle counts and
total distances of the heuristics and CTSPAV\textsubscript{Hybrid} are
very similar. In fact, the vehicle count results of the heuristic are
identical to those of CTSPAV\textsubscript{Hybrid} in all except three
instances: L19, M15, and S7. For these three instances, the counts of
the heuristic are only greater than those of
CTSPAV\textsubscript{Hybrid} by one vehicle. Moreover, the percentage
difference in the total distance results are consistently less than
1.50\% (on average, they differ by 0.01\%). This similarity bodes very
well for the heuristic; it highlights the effectiveness of its
negative reduced cost criterion for identifying the subset of mini
routes that are critical for producing strong integer solutions. It
also indicates that the heuristic is more than sufficient for
producing high-quality solutions, especially in applications whereby
proving the optimality of the final solution is not of paramount
importance. As mentioned earlier, the heuristic is incapable of
closing the vehicle count or optimality gap for any of the instances,
so CTSPAV\textsubscript{Hybrid} remains the better candidate in
applications where closing or narrowing the optimality gap is
critical.

\section{Case Study of Shared Commuting in Ann Arbor, Michigan} \label{sec:casestudy}

This section summarizes the results of a case study that applies the
CTSPAV to optimize the commuting trips from the Ann Arbor
dataset. More specifically, it considers all trips (of commuters
living inside and outside city limits) for the first four weekdays
(Monday--Thursday) of the busiest week of April 2017 (week 2). The
parameters $N$, $\Delta$, and $R$ are set to 100, 10 minutes, and 50\%
respectively for this case study.\footnotemark Its goal is to
demonstrate the effectiveness of the CTSPAV at reducing vehicle usage
and miles traveled, as well as to examine some of the real-world
benefits and drawbacks of the AV ridesharing platform.

\footnotetext{Part of the results for this case study is obtained by
  performing further analysis on the results from an earlier work
  \citep{hasan2021} which utilized the column-generation heuristic to
  solve the CTSPAV.}

\begin{figure}[!t]
	\begin{subfigure}[b]{0.5\linewidth}
		\centering
		\includegraphics[width=0.98\linewidth]{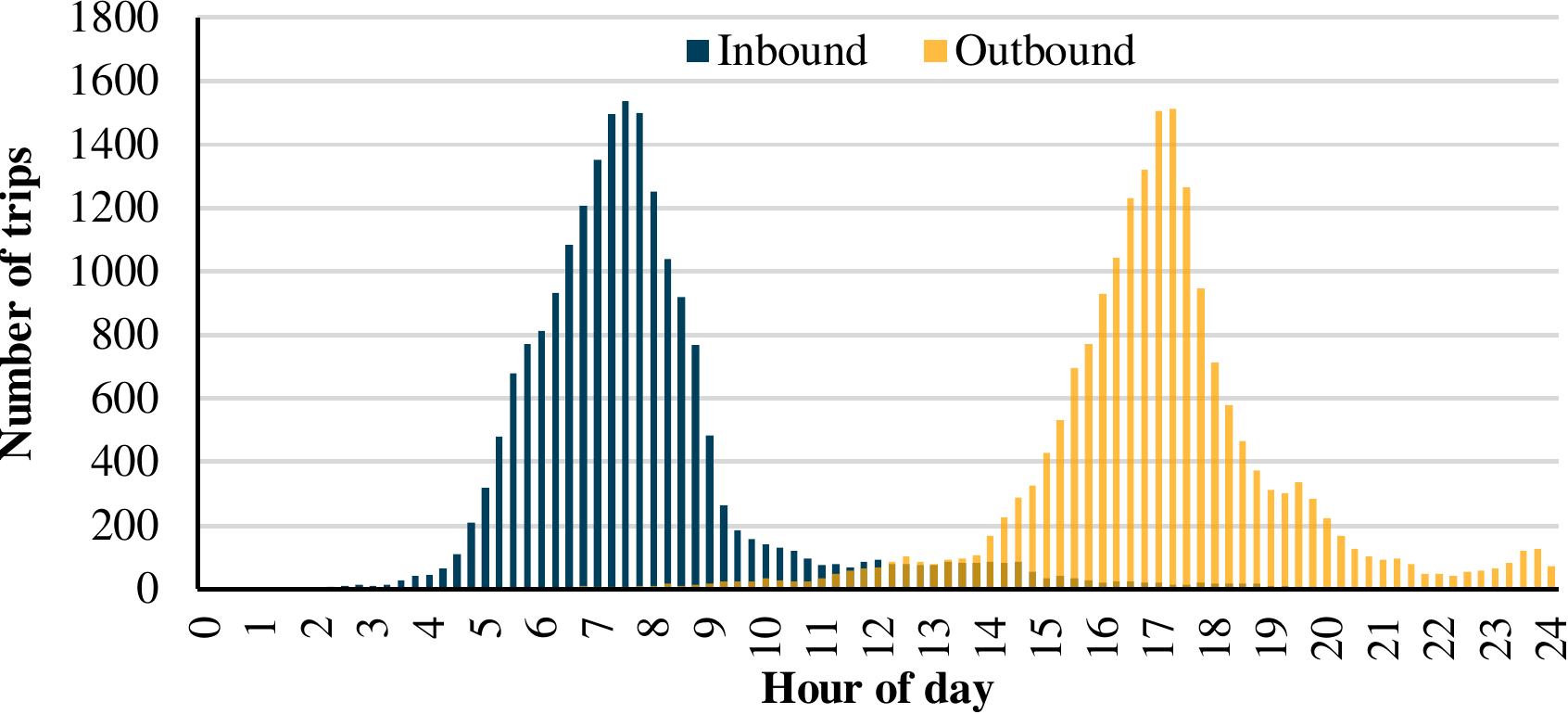} 
		\caption{Monday} 
		\label{fig:demand_a} 
		\vspace{1ex}
	\end{subfigure}
	\begin{subfigure}[b]{0.5\linewidth}
		\centering
		\includegraphics[width=0.98\linewidth]{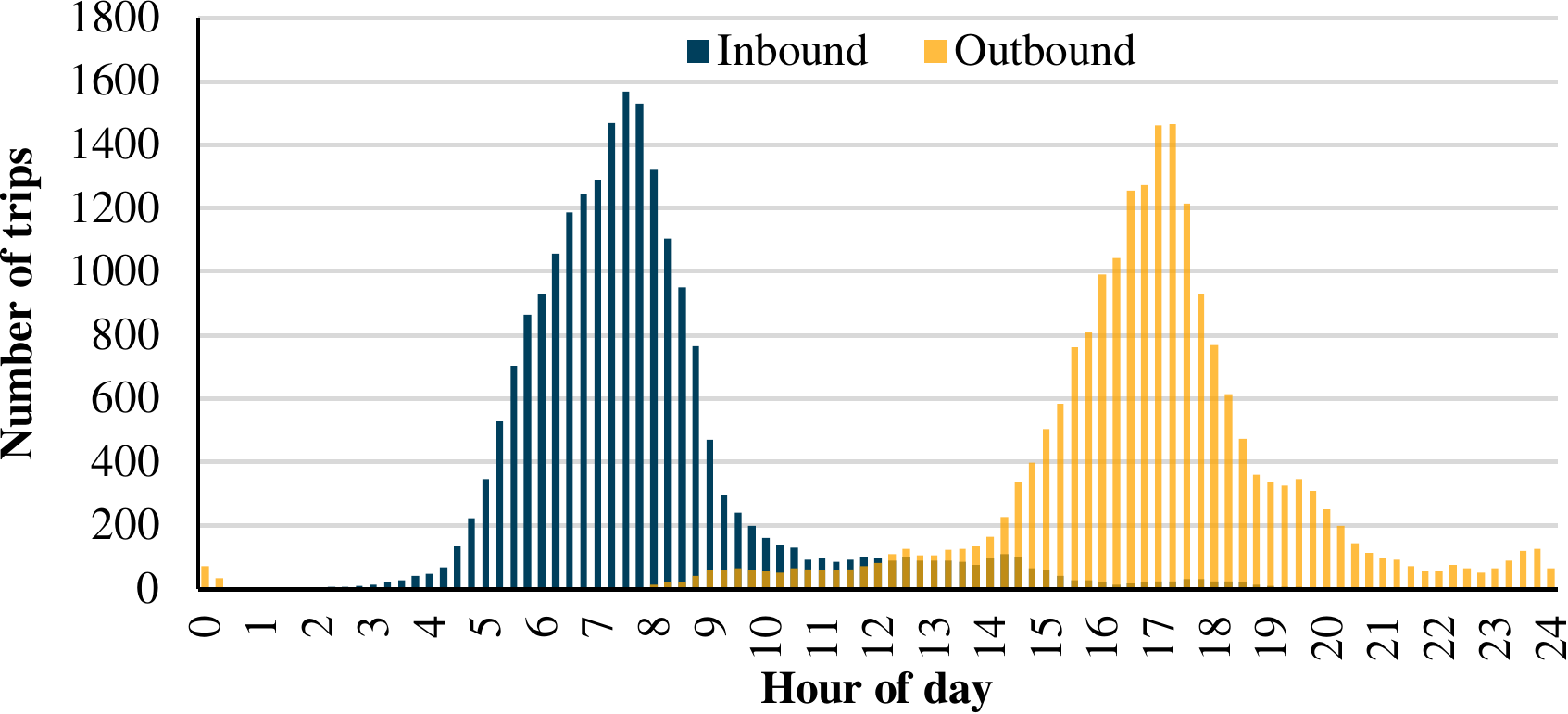} 
		\caption{Tuesday} 
		\label{fig:demand_b} 
		\vspace{1ex}
	\end{subfigure} 
	\begin{subfigure}[b]{0.5\linewidth}
		\centering
		\includegraphics[width=0.98\linewidth]{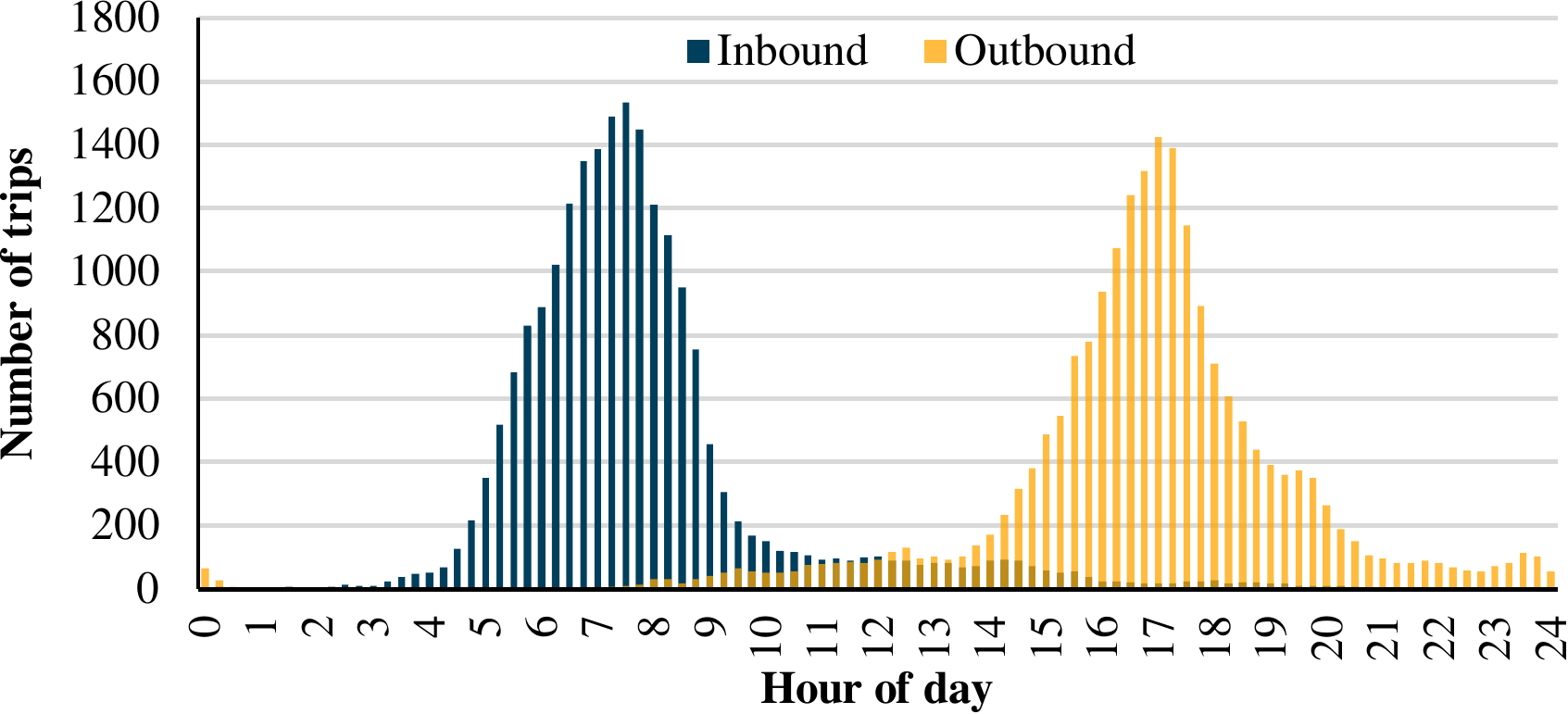} 
		\caption{Wednesday} 
		\label{fig:demand_c} 
	\end{subfigure}
	\begin{subfigure}[b]{0.5\linewidth}
		\centering
		\includegraphics[width=0.98\linewidth]{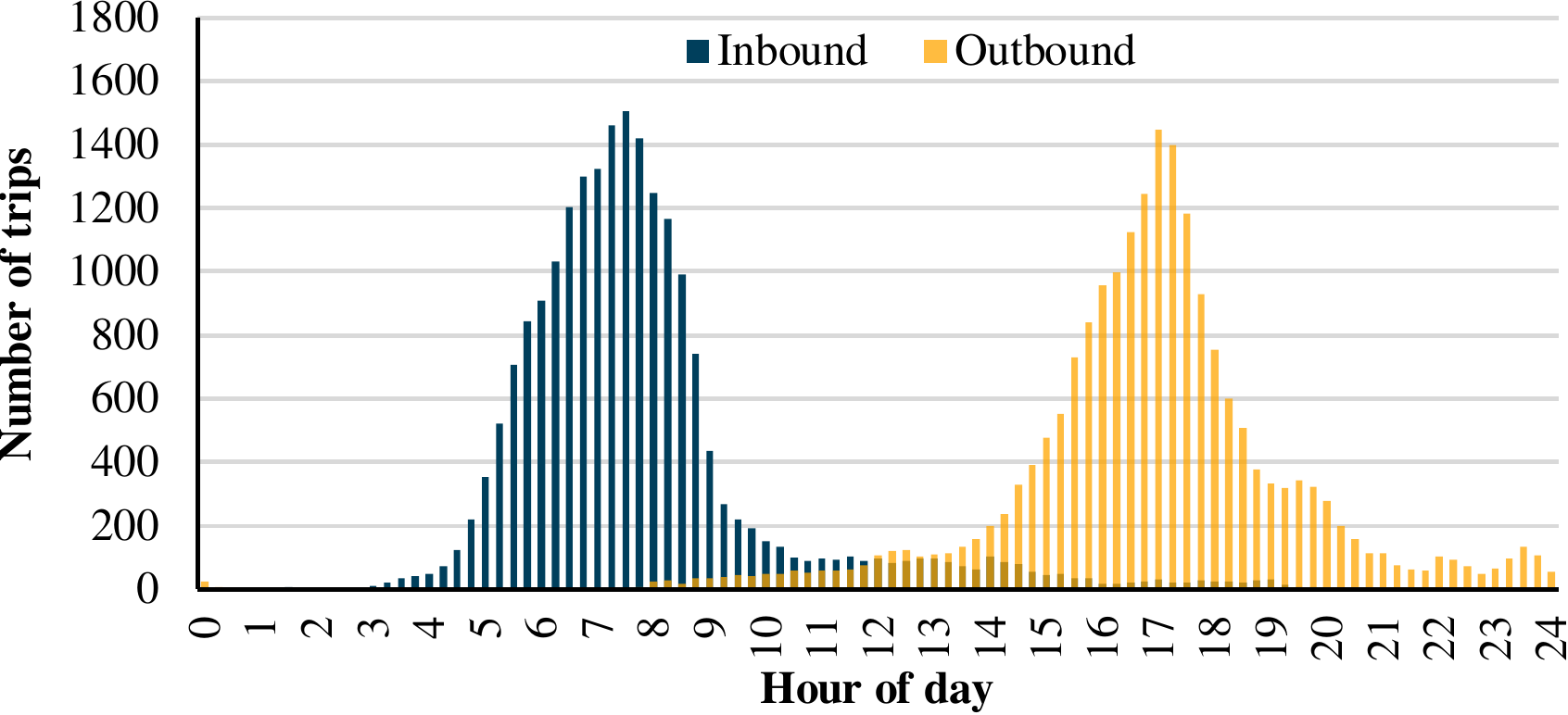} 
		\caption{Thursday} 
		\label{fig:demand_d} 
	\end{subfigure} 
	\caption{Commute Trip Demand Over 15-Minute Intervals on Week 2.}
	\label{fig:demand} 
\end{figure}

Figure \ref{fig:demand} provides an overview of the trip demand from
the dataset and reports the number of ongoing trips for every
15-minute interval throughout the four days considered. The data
exhibits clear and consistent commuting patterns: the inbound demand
peaks between 7--8 am, and the outbound demand peaks at around 5 pm
every day. The highly consistent nature of the trip distributions
highlights the opportunities in optimizing them.

\begin{figure}[!t]
	\centering
	\includegraphics[width=1.0\linewidth]{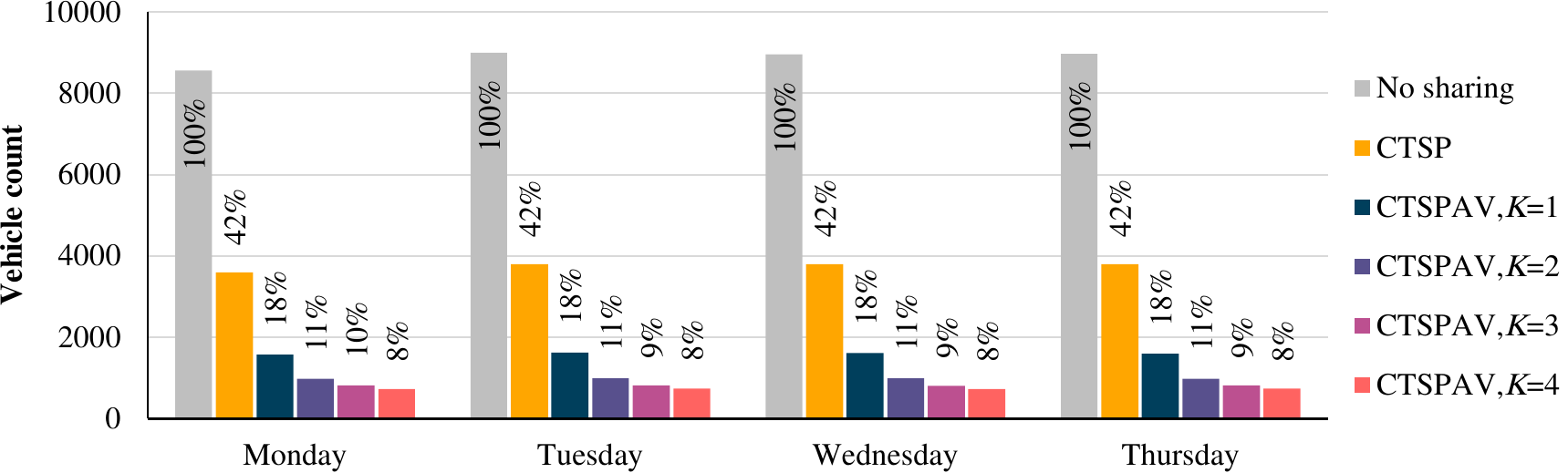}
	\caption{Total Number of Cars Used on Week 2.}
	\label{fig:vehicle_count}
\end{figure}

\subsection{Reductions in Vehicle Counts and Travel Distances}

Figure \ref{fig:vehicle_count} summarizes results of the primary
objective of the CTSPAV for various vehicle capacities
$K\in\{1,2,3,4\}$. It reports the total number of vehicles needed to
cover all trips for each $K$ value by aggregating the final vehicle
count results of every cluster. The number of vehicles utilized under
{\em no-sharing} conditions (i.e., when commuters travel using their
personal vehicles) and under the original CTSP (with $K=4$) (i.e.,
when drivers are selected from the set of commuters) are included for
additional perspectives. The percentages in the figure report each
count as a fraction of the no-sharing count. {\em The figure highlights the
significant capability of the CTSPAV in reducing the number of
vehicles. Indeed, the CTSPAV reduces the vehicle counts by up to 92\%
every day}, and improves upon the original CTSP by an additional
34\%. In fact, the results show that, even without any ride sharing
(i.e., when $K=1$), AVs still reduce the number of vehicles by 82\%
and improve upon the CTSP by an additional 24\%. This reduction in
vehicle count can be translated into a significant reduction in
parking spaces, which can then be utilized for other, more useful,
infrastructures. The difference in vehicle counts between the CTSP and
the CTSPAV is due to autonomy: the vehicles are not associated with
drivers and can travel back and forth between residential
neighborhoods and workplaces. In the CTSP, vehicles only make a single
inbound and outbound trip every day as their routes are restricted to
begin and end at the trip origins and destinations of their drivers.

\begin{figure}[!t]
	\centering
	\includegraphics[width=1.0\linewidth]{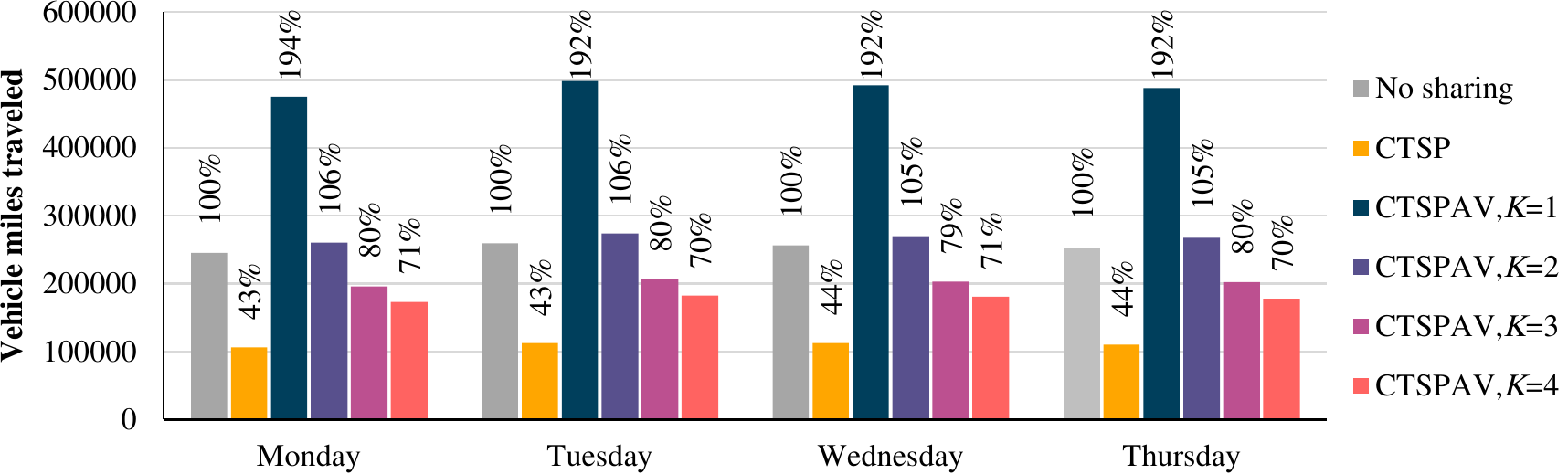}
	\caption{Total Travel Distance on Week 2.}
	\label{fig:vehicle_miles_traveled}
\end{figure}

\begin{figure}[!t]
	\centering
        \includegraphics[width=1.0\linewidth]{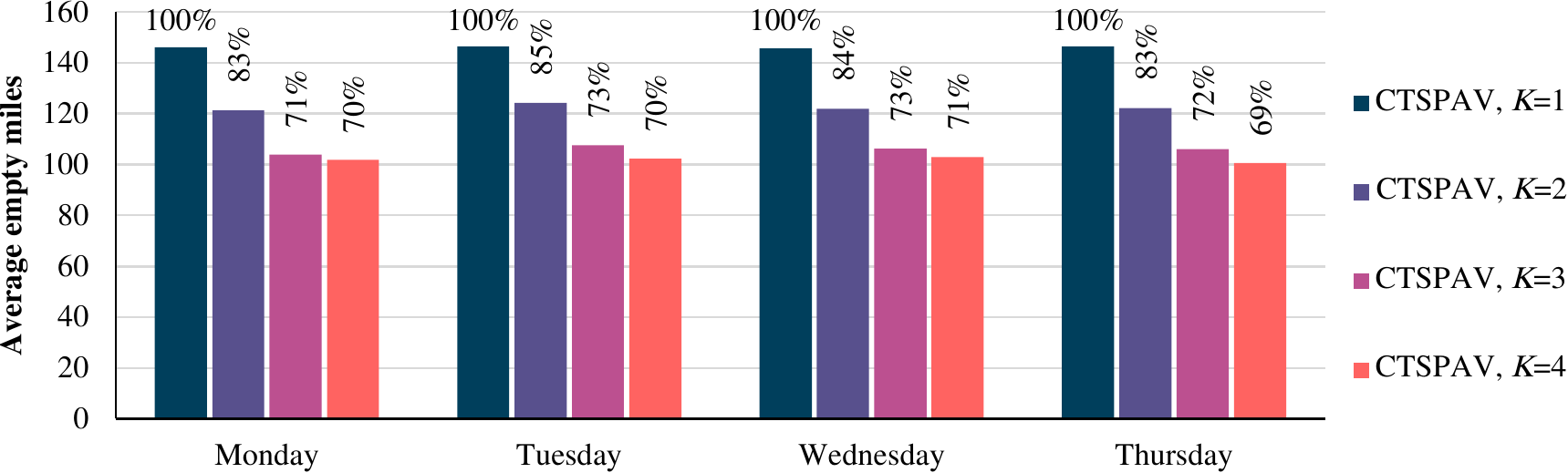}
	\caption{Average Empty Miles Per Vehicle on Week 2.}
	\label{fig:average_empty_miles}
\end{figure}

Figure \ref{fig:vehicle_miles_traveled} summarizes the total travel
distance of the vehicles, which is the secondary objective of the
CTSPAV, under the same configurations. The results are again obtained
by aggregating the results from every cluster and the percentages
represent each quantity as a fraction of the no-sharing total. The
first result that stands out is how many more miles are traveled by
the CTSPAV when $K=1$ (92--94\% more than those under no-sharing
conditions). When $K=1$ for the CTSPAV, the autonomous vehicles need
to perform significantly more back-and-forth traveling between the
neighborhoods and the workplace to cover the same amount of trips,
which consequently leads to their inflated total travel distance. The
results improve significantly when $K$ is increased to 2 as the
vehicles allow for more trip aggregations, yet the traveled miles are
still 5--6\% more than those for private vehicles. Net savings in
travel distance are only realized when $K\geq 3$: beyond this point,
{\em the reduction in travel distance from ride sharing exceeds the
additional empty miles (the miles traveled by an AV with no passengers
onboard) introduced by the back-and-forth traveling of the
AVs.} Nevertheless, the 29--30\% reduction in miles traveled when $K=4$
is still not as significant as that offered by the original CTSP which
is around 56--57\%. Indeed, the CTSP does not introduce any empty
miles and benefits from all the distance savings from ride sharing. On
the other hand, the CTSPAV total will necessarily include some empty
miles from when the vehicles travel without any passengers onboard as
they go from the workplace back to the residential neighborhoods in
the morning (or vice versa in the evening) to pick up more
trips. There is obviously a tradeoff between the reductions in vehicle
counts and travel distances. Figure \ref{fig:average_empty_miles}
provides a closer look at the average empty miles per vehicle for the
various vehicle capacities. The results are quite intuitive: the
average decreases as $K$ increases, since the larger vehicle
capacities allow for more ridesharing and require less back-and-forth
traveling to cover the same amount of trips.

\begin{figure}[!t]
	\centering
	\includegraphics[width=0.48\linewidth]{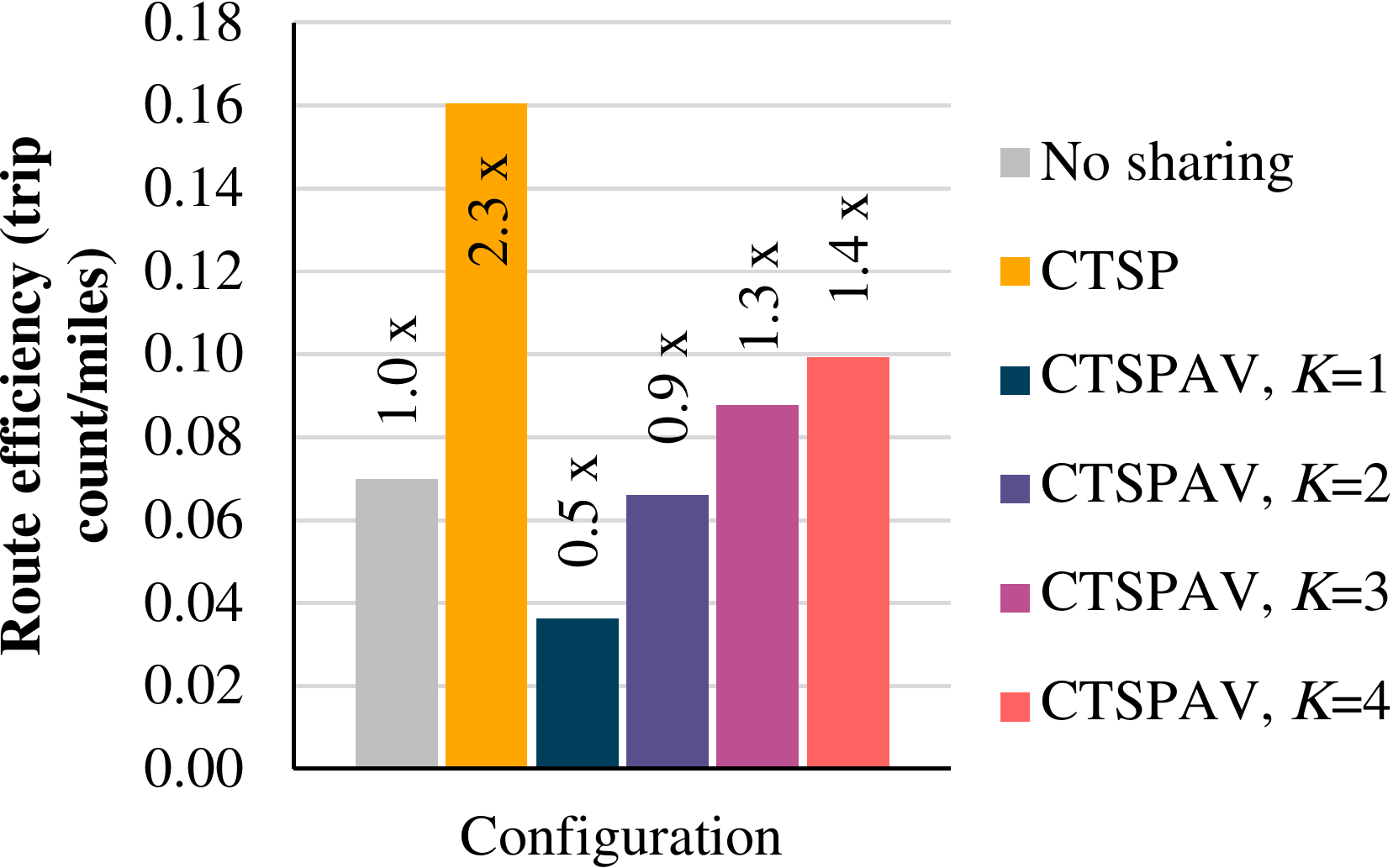}
	\caption{Efficiency of Vehicle Routes}
	\label{fig:route_eff}
\end{figure}

Figure \ref{fig:route_eff} then attempts to quantify the route
efficiency of of the various configurations, i.e., the number of trips
covered per mile traveled. It also includes a multiplicative factor
for each quantity as a multiple of the no-sharing value.
The results indicate that the CTSP produces the most efficient routes,
whereas the CTSPAV, when $K=1$, is the least efficient. The CTSPAV
gains more efficiency (albeit at a decreasing rate) as its vehicle
capacity increases: while its routes are more efficient than those of
the private vehicles when $K=4$, they still cannot outperform those of
the CTSP. There is an intuitive explanation for this observation. The
CTSPAV loses its route efficiency from its empty miles and then has to
recover them by maximizing ridesharing to cover as many trips as
possible. In contrast, the CTSP does not have to contend with any
efficiency losses due to empty miles.

\subsection{Congestion Analysis}

\begin{figure}[!t]
	\begin{subfigure}[b]{0.5\linewidth}
		\centering
		\includegraphics[width=0.98\linewidth]{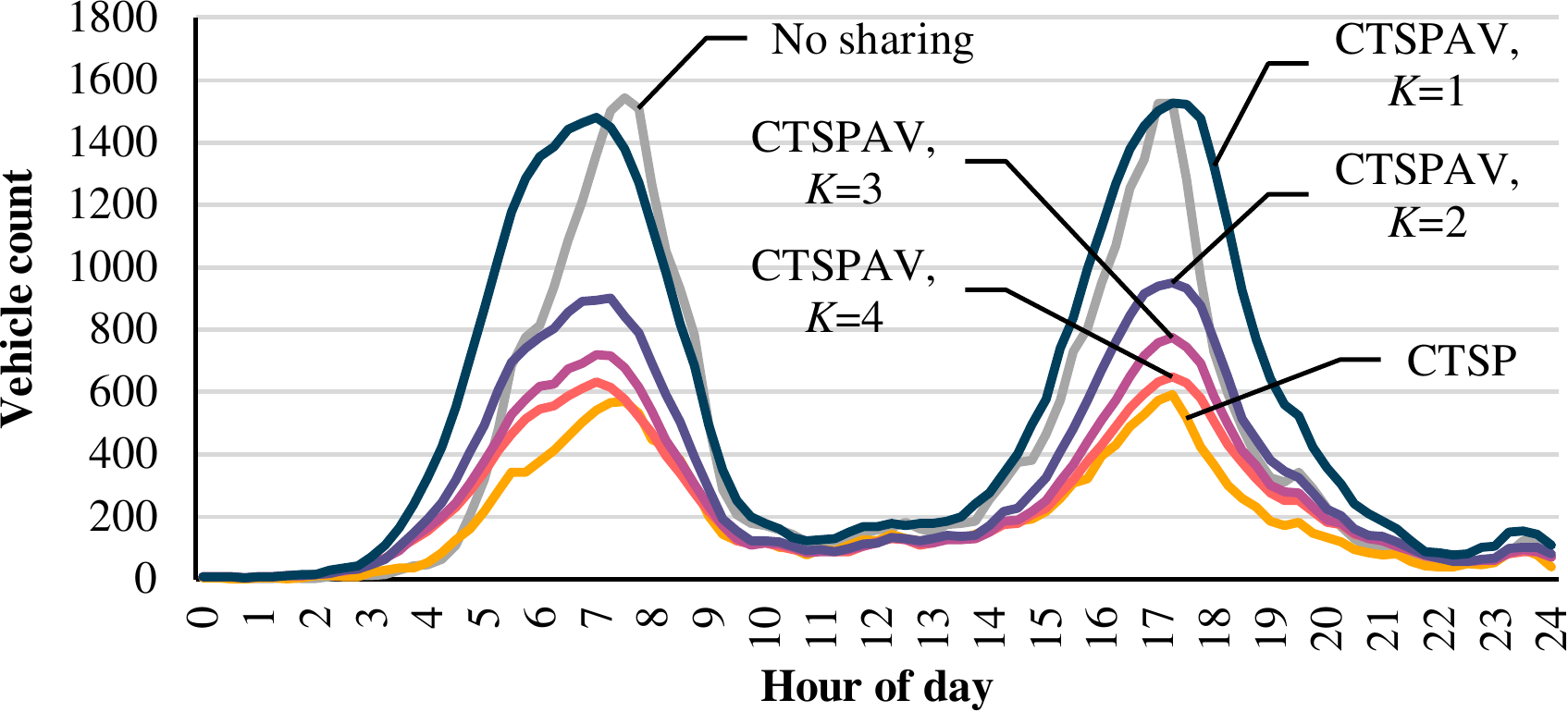} 
		\caption{Monday} 
		\label{fig:vehicle_usage_a} 
		\vspace{1ex}
	\end{subfigure}
	\begin{subfigure}[b]{0.5\linewidth}
		\centering
		\includegraphics[width=0.98\linewidth]{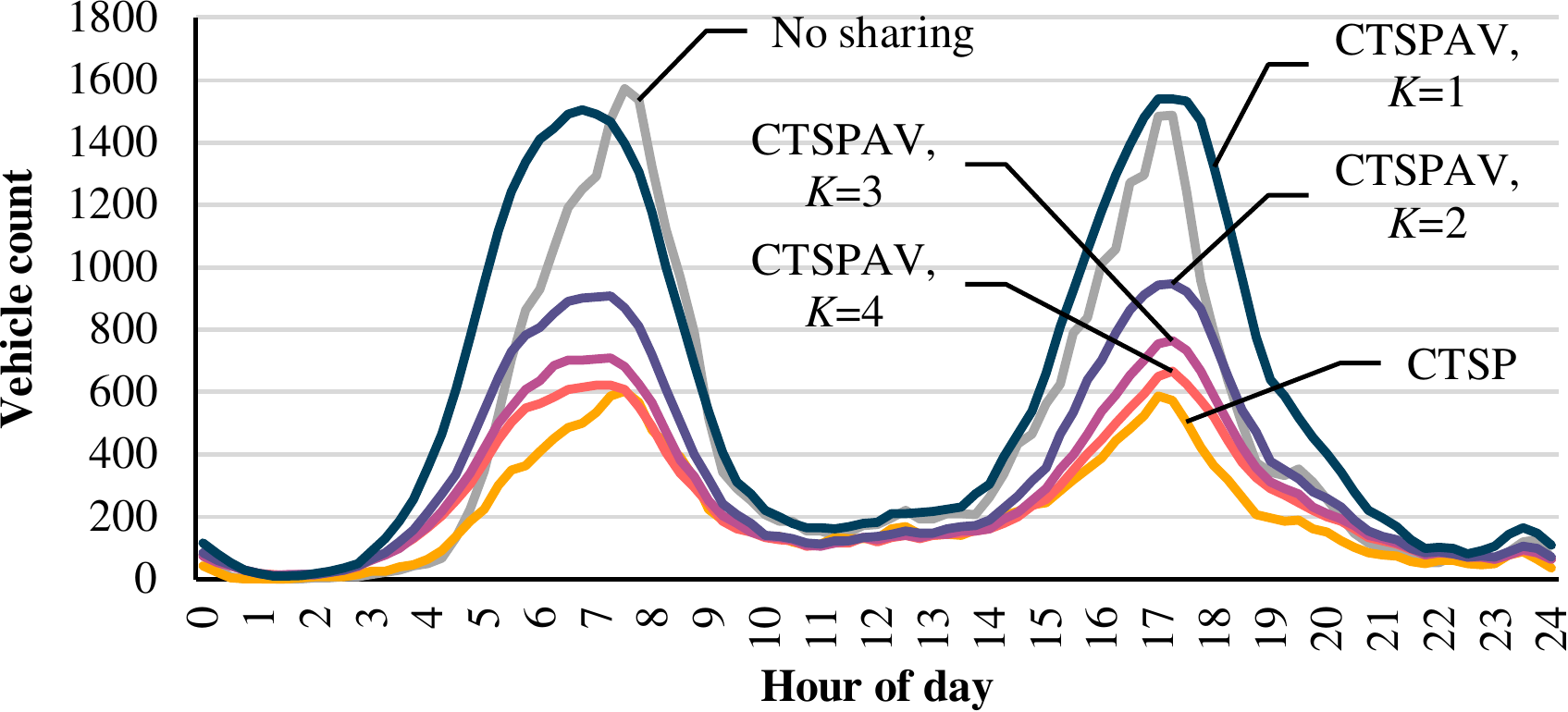} 
		\caption{Tuesday} 
		\label{fig:vehicle_usage_b} 
		\vspace{1ex}
	\end{subfigure} 
	\begin{subfigure}[b]{0.5\linewidth}
		\centering
		\includegraphics[width=0.98\linewidth]{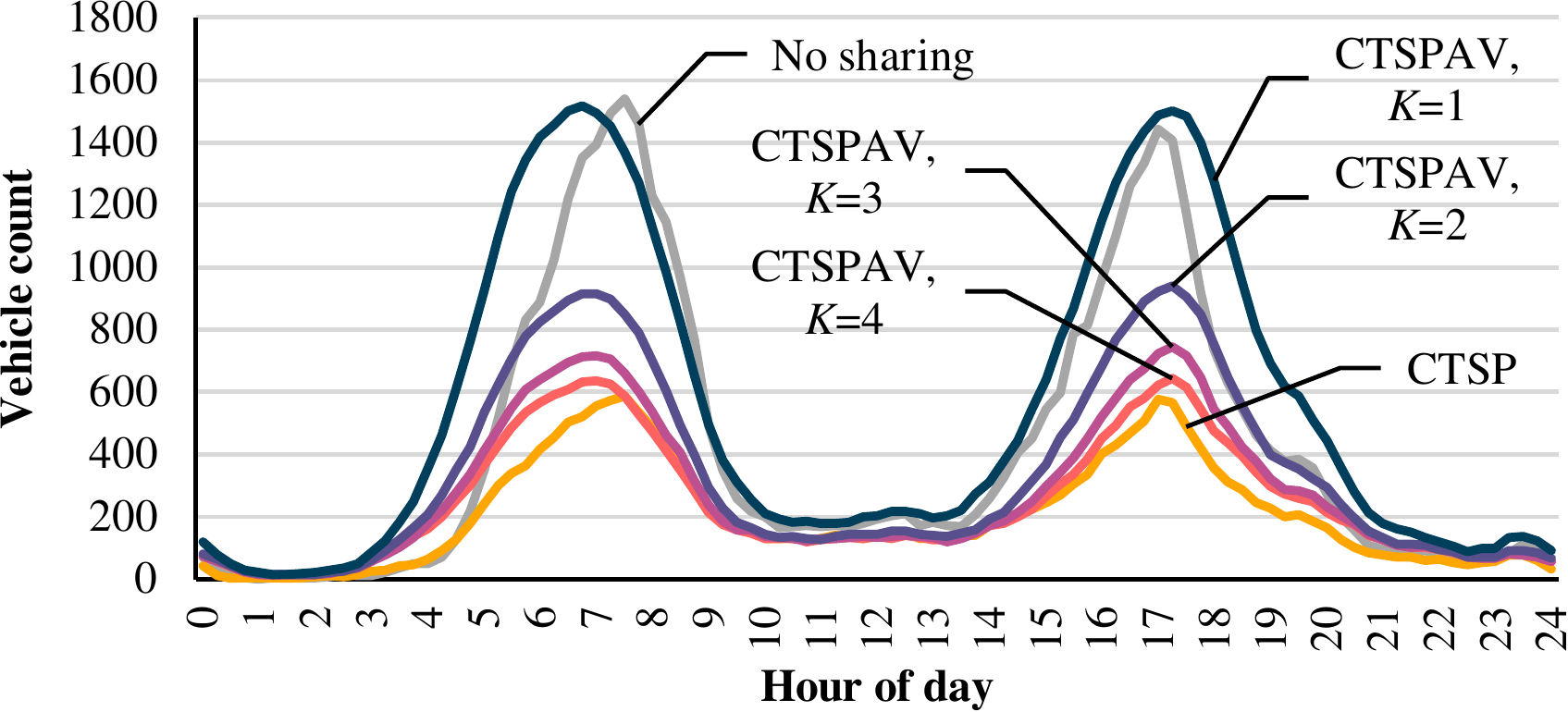} 
		\caption{Wednesday} 
		\label{fig:vehicle_usage_c} 
	\end{subfigure}
	\begin{subfigure}[b]{0.5\linewidth}
		\centering
		\includegraphics[width=0.98\linewidth]{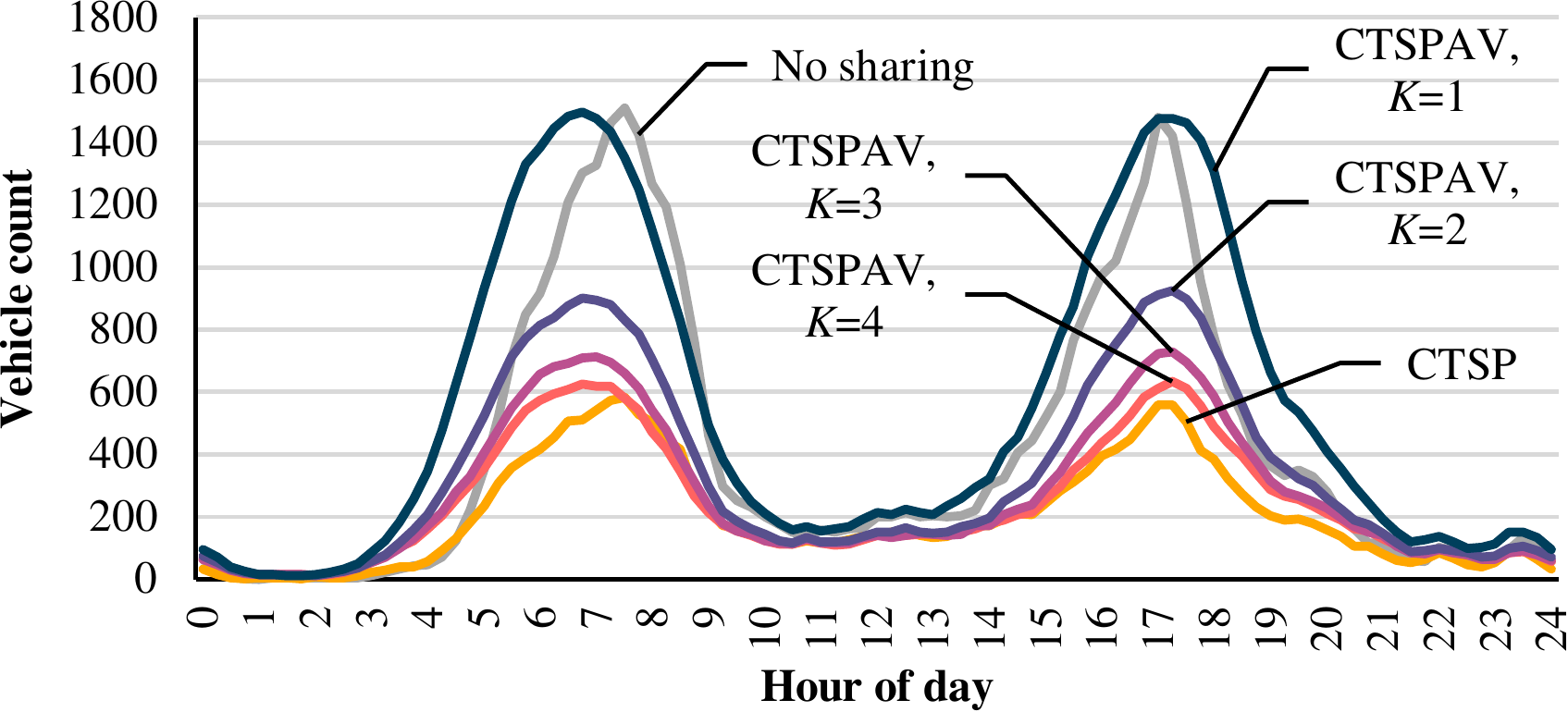} 
		\caption{Thursday} 
		\label{fig:vehicle_usage_d} 
	\end{subfigure} 
	\caption{Number of Vehicles on the Road Over 15-Minute Intervals on Week 2.}
	\label{fig:vehicle_usage} 
\end{figure}

Figure \ref{fig:vehicle_usage} presents results on congestion to
understand the reduction (or increase) in traffic caused by AVs
compared to the no-sharing condition. It tallies the total number of
vehicles used by each configuration over every 15-minute interval
throughout the four days considered. The goal is to investigate,
qualitatively and comparatively, the capability of each configuration
in flattening the traffic curve originally produced by the private
vehicles. The CTSPAV with $K = 1$ appears to aggravate traffic as its
curve is as tall as, and is wider than, that of private vehicles. This
is not surprising. As illustrated earlier, this configuration produces
the largest amount of vehicle miles traveled and also the most empty
miles. The curve is drastically flattened as soon as $K$ increases to
2, and it keeps becoming flatter (at a decreasing rate) as $K$ further
increases. When $K=4$, the CTSPAV produces about a 60\% reduction in
traffic. The traffic curves of the CTSP appear to dominate slightly
those of the CTSPAV with $K=4$ most of the time. This observation is
also in line with the route efficiency calculations.  However,
regardless of their relative performance, Figure
\ref{fig:vehicle_usage} provides evidence that {\em both the CTSP and
CTSPAV have the potential to significantly reduce traffic congestion
and parking utilization.}

\subsection{Analysis of Commuting Properties}

\begin{figure}[!t]
	\begin{subfigure}[b]{0.5\linewidth}
		\centering
		\includegraphics[width=0.98\linewidth]{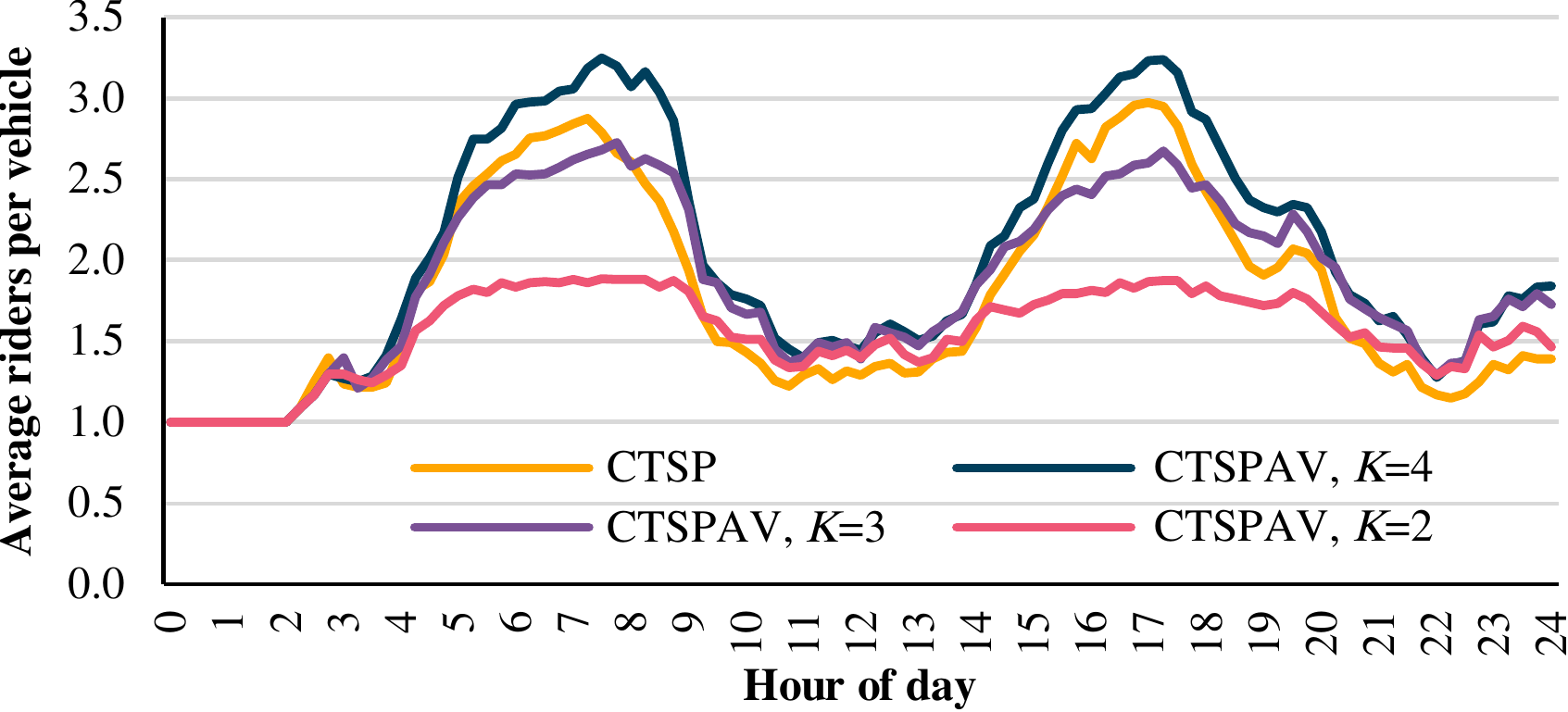} 
		\caption{Monday} 
		\label{fig:sharing_a} 
		\vspace{1ex}
	\end{subfigure}
	\begin{subfigure}[b]{0.5\linewidth}
		\centering
		\includegraphics[width=0.98\linewidth]{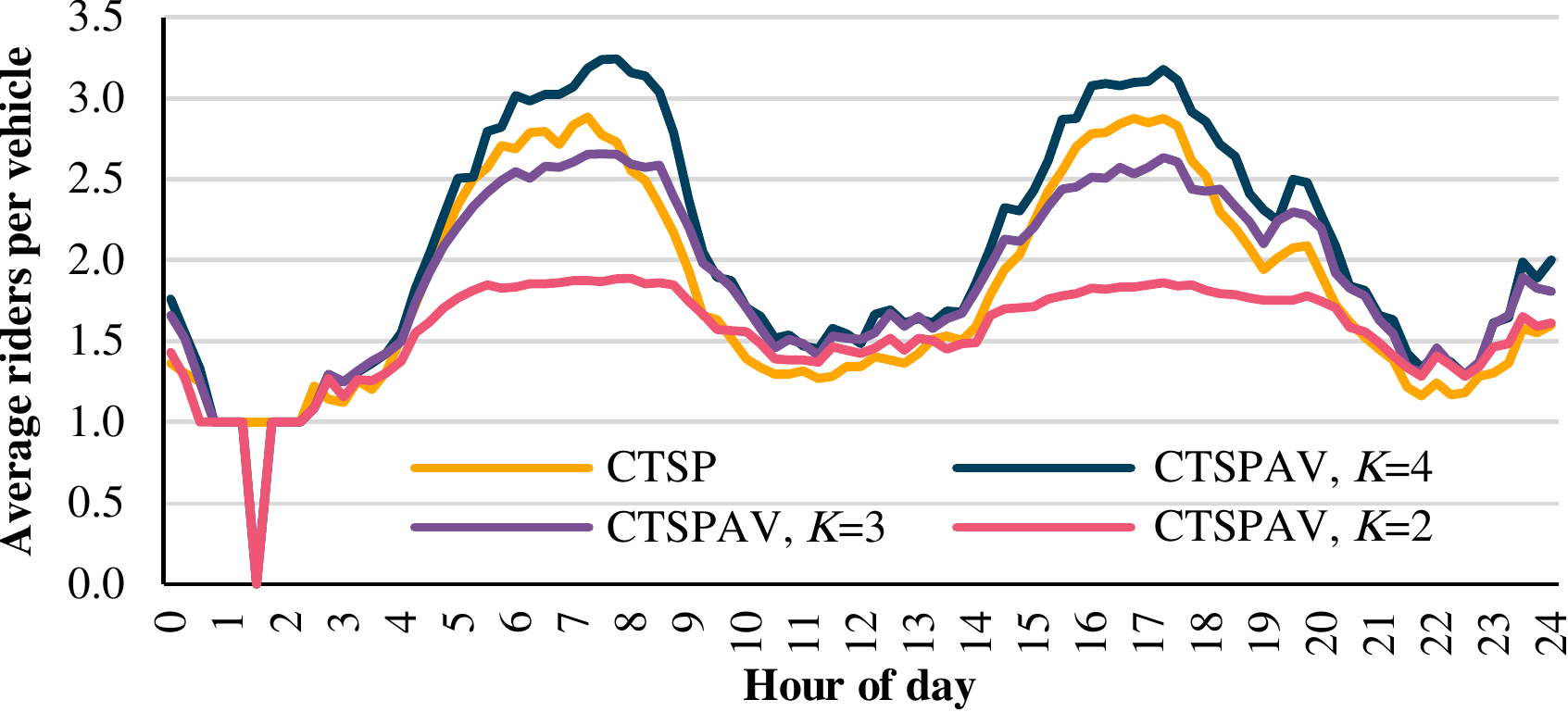} 
		\caption{Tuesday} 
		\label{fig:sharing_b} 
		\vspace{1ex}
	\end{subfigure} 
	\begin{subfigure}[b]{0.5\linewidth}
		\centering
		\includegraphics[width=0.98\linewidth]{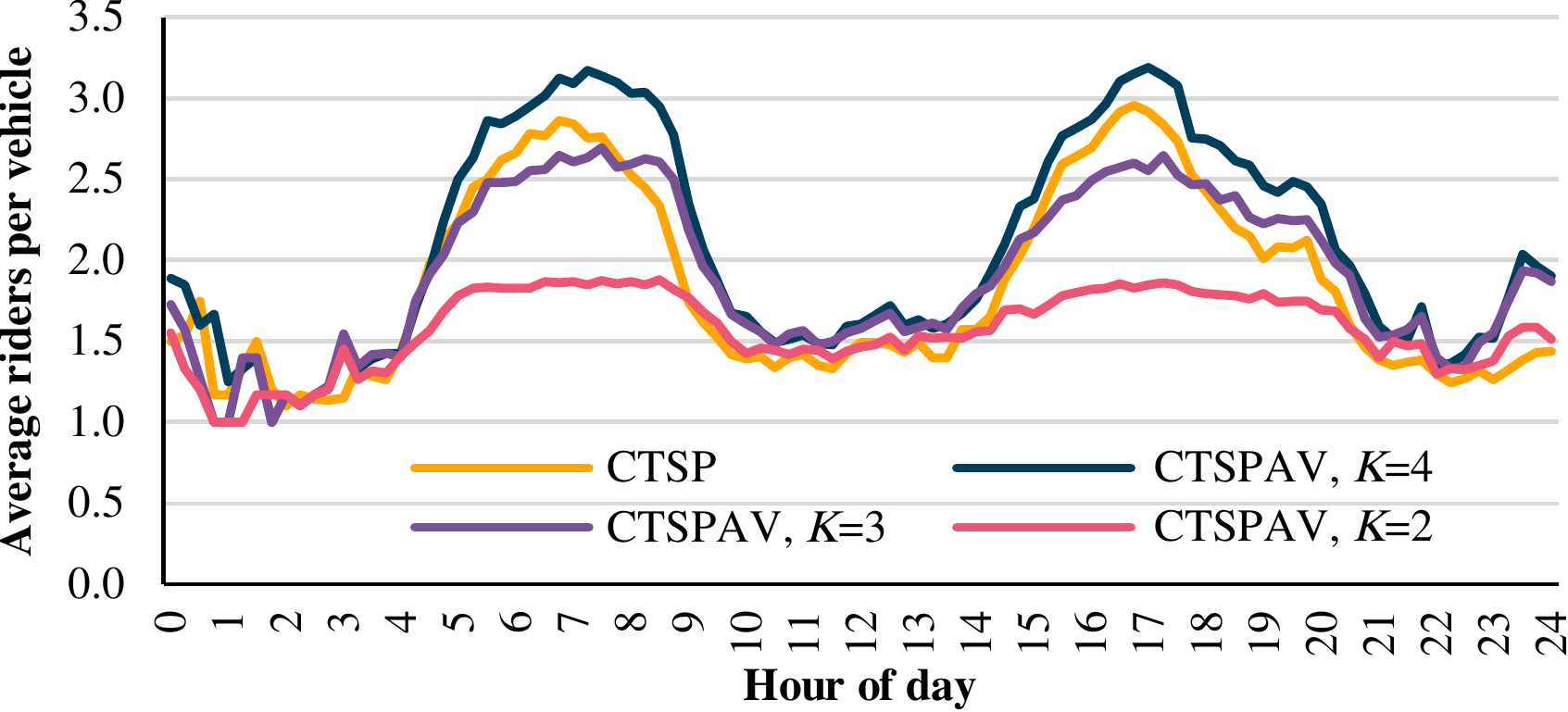} 
		\caption{Wednesday} 
		\label{fig:sharing_c} 
	\end{subfigure}
	\begin{subfigure}[b]{0.5\linewidth}
		\centering
		\includegraphics[width=0.98\linewidth]{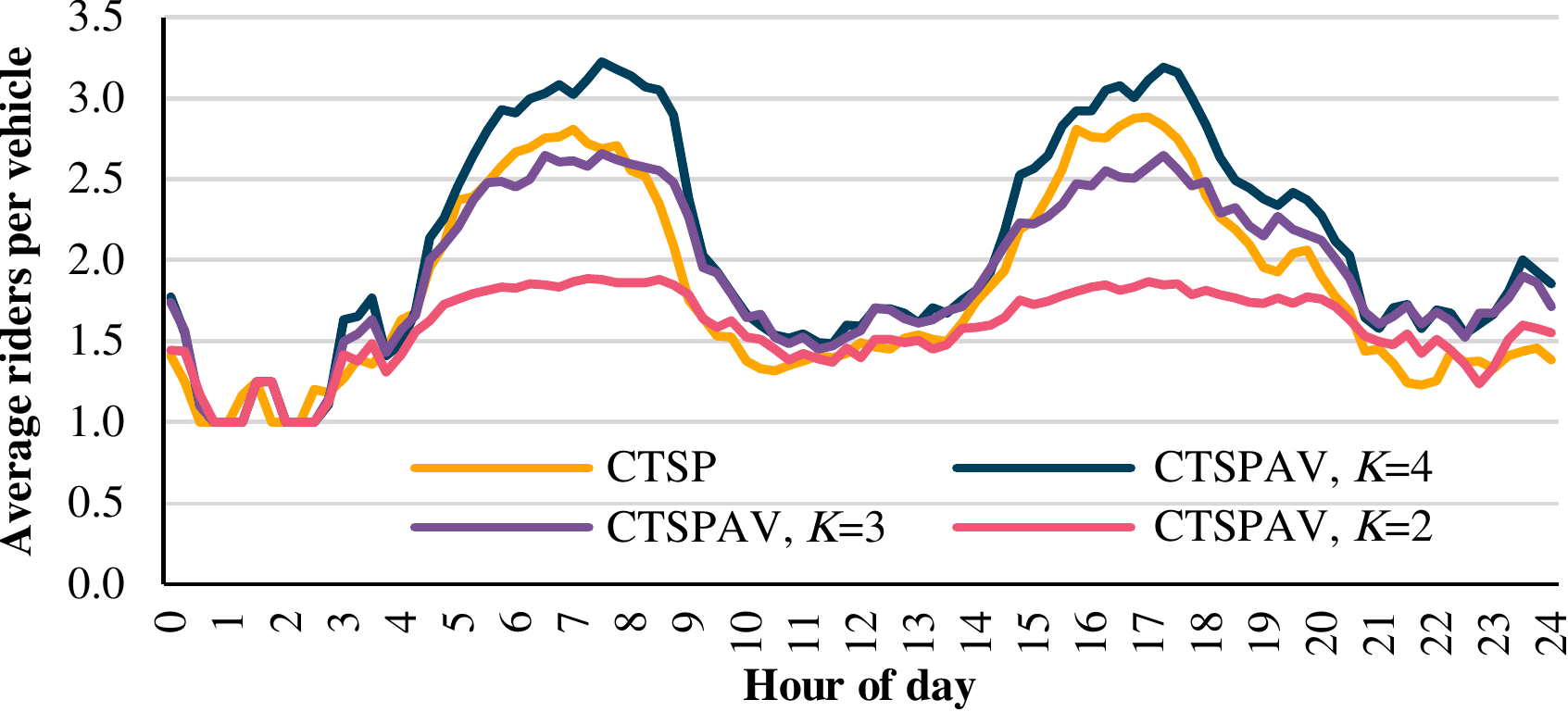} 
		\caption{Thursday} 
		\label{fig:sharing_d} 
	\end{subfigure} 
	\caption{Average Riders Per Vehicle Over 15-Minute Intervals on Week 2}
	\label{fig:sharing} 
\end{figure}

Figure \ref{fig:sharing} aims to quantify the relative amount of ride
sharing taking place throughout each day for the different
configurations. It reports the average number of riders per vehicle
for every 15-minute interval throughout the four days
considered. Results for the private vehicles and for the CTSPAV with
$K=1$ are not included for obvious reasons (they do not allow any
sharing). The amount of ride sharing throughout a typical weekday
mimics the shape of the trip demand: they both peak during the same
periods of the day. This is to be expected as the CSTSP and CTSPAV
maximize ride sharing, which is easier when the trip demand is
higher. The figure also shows that the relative amount of sharing for
the CTSPAV increases with vehicle capacity. {\em Moreover, when $K=4$,
  there is more ride sharing in the CTSPAV than in the CTSP most of
  the time.} This can be attributed to the relative flexibility of the
mini routes of the CTSPAV compared to those of the CTSP. Indeed, a
CTSP route must start and end at the orign and destination of its
driver, which constrains its total duration by the ride-duration
constraints on its driver.  Mini routes of the CTSPAV are not
subjected to these restrictions, allowing for more flexibility in
serving trips. Interestingly, during peaks, the average amount of ride
sharing is between 3.0 and 3.5 due to the spatial and temporal
properties of the commuting trips. This also indicates the types of
autonmous vehicles that will be most useful in the future, at least
for cities like Ann Arbor.

\begin{figure}[!t]
	\centering
	\includegraphics[width=1.0\linewidth]{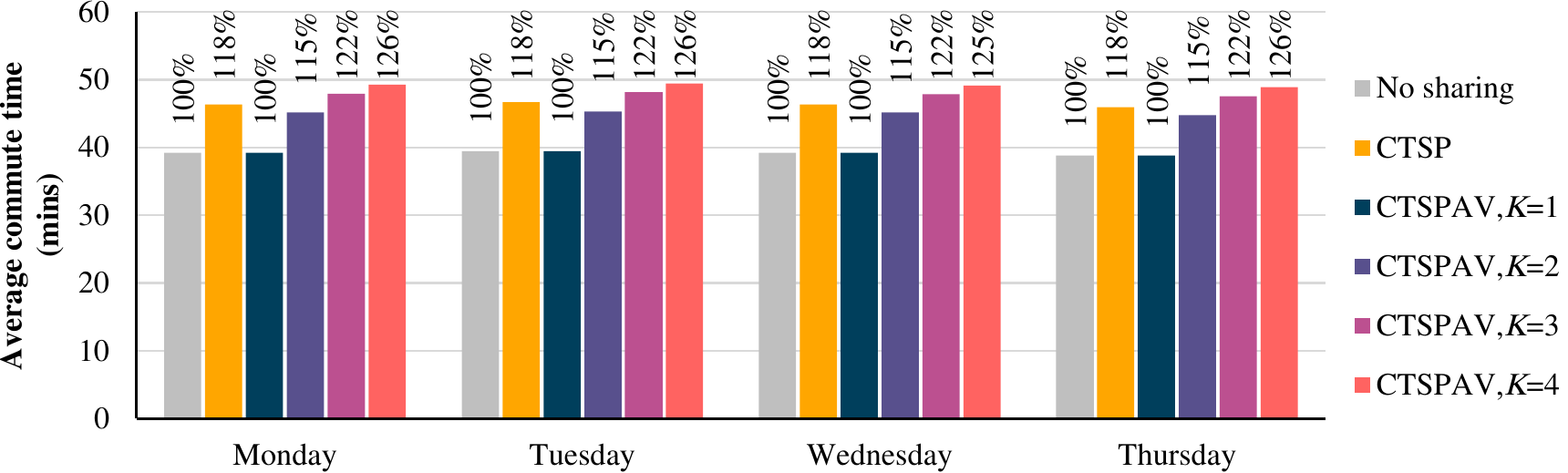}
	\caption{Average Commute Time on Week 2.}
	\label{fig:average_commute_time}
\end{figure}

Figure \ref{fig:average_commute_time} reports the average commute
times, i.e., the average time spent on the vehicle by each rider.  The
percentages of each quantity are calculated relative to the no-sharing
value. The results shed light on another inherent trade-off in
ride-sharing service as the ride duration necessarily
incraeses. During ridesharing, a route may deviate from the optimal
path to pickup or drop off other riders. This, combined with possible
wait times incurred at the pickup locations, contribute to the
increased ride duration. The results reveal an expected trend for the
CTSPAV: the average commute times increase with an increase in vehicle
capacity. However, {\em it is interesting to observe that, although
parameter $R$ was set to 50\% for the case study, the commute times of
the CTSPAV with $K=4$ only increase by an average of 26\%. The CTSPAV
thus guarantees a high quality of service for its riders.}

\section{Conclusion} \label{sec:conclusion}

The purpose of the CTSPAV is to synthesize an optimal routing plan for
serving a large set of commute trips with AVs. Its design was
originally motivated by the desire to address the growing parking and
traffic congestion problems induced by the average of 9,000 daily
commuters traveling to parking lots operated by the University of
Michigan located in downtown Ann Arbor, Michigan. Utilization of AVs
was seen as the key to addressing the shortcomings of the original
CTSP---a conventional car-pooling problem with the same objectives as
the CTSPAV---by obviating any driver-related requirements that could
limit its ridesharing potential. A first attempt at solving the
problem by \cite{hasan2021} investigated two different methods: (1) A
CTSPAV procedure which used column-generation to discover mini
routes---short routes covering only inbound or outbound trips that
have distinct pickup, transit, and drop-off phases---with negative
reduced costs which are chained together to form longer AV routes in
its master problem and (2) A DARP procedure which uses a classical
column-generation approach originally developed for the DARP to solve
the CTSPAV. Both methods utilized identical lexicographic objectives
which sought to first minimize the required vehicle count and then
minimize their total travel distance. To deal with the complexity of
handling the massive volume of trips, the commuters were first
clustered into groups representing artificial neighborhoods, after
which ridesharing within each cluster was optimized exclusively. They
discovered that each method had a trade-off: The CTSPAV procedure
produced strong integer solutions but had weak primal lower
bounds. Conversely, the DARP procedure generated stonger primal lower
bounds especially for the primary objective, but it was slow and
therefore could not obtain strong integer solutions within
time-constrained scenarios.

{\em The trade-offs of the two procedures presented an opportunity for
  exploring a method that could leverage the strengths of both, which
  is the primary methodological contribution of this work.} This paper
thus proposed a branch-and-cut procedure that exploits a dual-modeling
approach for solving the CTSPAV. The core of the procedure is a MIP
formulation of the CTSPAV that chains (exhaustively enumerated) mini
routes to form longer AV routes and is capable of producing
high-quality integer solutions quality. This core is complemented by a
DARP formulation whose relaxation (for minimizing vehicle counts) is
obtained through a column-generation procedure. The DARP formulation
is less effective in finding high-quality integer solution, but its
relaxation produces stronger lower bounds. The overall algorithm
solves the core branch-and-cut procedure and the DARP relaxation in
parallel, transmitting new lower bounds asynchrously from the
relaxation to the branch and cut procedure. Computational evaluations
that use instances derived from the Ann Arbor commute-trip data
demonstrated that this hybrid algorithm consistently outperforms a
similar branch-and-cut procedure that utilizes other well-established
valid inequalities like 2-path cuts and successor and predecessor
inequalities. It also successfully closes the optimality gaps for
several large and medium-sized instances as well as those for all
tight problem instances considered in the evaluation, of which none
could be optimally solved by the CTSPAV procedure of \cite{hasan2021}.

With the availabilty of an exact branch and cut procedure, the paper
then provided a comprehensive analysis of the potential of AVs for
ride-sharing platforms and relieving parking pressure and congestion
in medium-sized cities. In particular, the paper presented results of
a case study which applies the clustering-CTSPAV optimization workflow
on a large-scale dataset of commute trips from the city of Ann Arbor,
Michigan. {\em The analysis revealed several invaluable insights,
  including the CTSPAV capability of reducing daily vehicle counts by
  92\%, further improving upon the already massive 57\% vehicle
  reductions of the original CTSP}. It does so by generating AV routes
that are very long---a stark contrast to the short routes of the
CTSP---allowing each AV to cover significantly more trips every
day. It could also effectively flatten the vehicle usage curve (i.e.,
the number of vehicles used per unit time), suggesting a concomitant
ability to effectively reduce traffic congestion. The CTSPAV also
produced higher averages for trips shared per unit time than the CTSP,
indicating that it is superior at aggregating more trips for
ridesharing. The analysis also revealed some drawbacks, the most
significant being the introduction of empty miles into the daily
travel distance totals. The empty miles degrade the efficiency of the
CTSPAV routes, which measures the average number of trips covered per
distance traveled, making them less efficient than the routes of the
CTSP. Empty miles are unfortunately a by-product that is inherent to
the utilization of AVs, and its introduction is a trade-off that will
need to be carefully weighed against the benefits of AVs by the
ridesharing platform operator. Nonetheless, the results indicate that
the CTSPAV routing plan, even with its empty miles, is still able to
reduce the total miles traveled by private vehicles by 30\% while
producing routes that at 1.4 times more efficient. On the whole, the
case study shows that a CTSPAV-based ridesharing platform could
significantly reduce daily vehicle counts, as well as the number of
vehicles used per unit time. Such a platform would be highly effective
at aggregating trips, making it a very promising solution for reducing
parking space utilization and for mitigating traffic congestion
induced by large-scale commuting.

%
%
%

\ACKNOWLEDGMENT{We would like to thank Stephen Dolen from Logistics,
  Transportation, and Parking of the University of Michigan for his
  assistance in obtaining the dataset used in this research. Part of
  this research was funded by the Rackham Graduate Student Research
  Grant, computational resources and services provided by Advanced
  Research Computing at the University of Michigan, NSF Leap HI
  proposal NSF-1854684, and Department of Energy Research Award
  7F-30154.}



\begin{thebibliography}{78}
\providecommand{\natexlab}[1]{#1}
\providecommand{\url}[1]{\texttt{#1}}
\providecommand{\urlprefix}{URL }

\bibitem[{Agatz et~al.(2012)Agatz, Erera, Savelsbergh, \protect\BIBand{}
  Wang}]{agatz2012}
Agatz N, Erera A, Savelsbergh M, Wang X (2012) Optimization for dynamic
  ride-sharing: A review. \emph{European Journal of Operational Research}
  223(2):295 -- 303, ISSN 0377-2217,
  \urlprefix\url{http://dx.doi.org/https://doi.org/10.1016/j.ejor.2012.05.028}.

\bibitem[{Agatz et~al.(2011)Agatz, Erera, Savelsbergh, \protect\BIBand{}
  Wang}]{agatz2011}
Agatz NA, Erera AL, Savelsbergh MW, Wang X (2011) Dynamic ride-sharing: A
  simulation study in metro atlanta. \emph{Transportation Research Part B:
  Methodological} 45(9):1450 -- 1464, ISSN 0191-2615,
  \urlprefix\url{http://dx.doi.org/https://doi.org/10.1016/j.trb.2011.05.017},
  select Papers from the 19th ISTTT.

\bibitem[{Alazzawi et~al.(2018)Alazzawi, Hummel, Kordt, Sickenberger,
  Wieseotte, \protect\BIBand{} Wohak}]{alazzawi2018}
Alazzawi S, Hummel M, Kordt P, Sickenberger T, Wieseotte C, Wohak O (2018)
  Simulating the impact of shared, autonomous vehicles on urban mobility -- a
  case study of milan. Wie\{\textbackslash{}ss\}ner E, L\textbackslash{}"ucken
  L, Hilbrich R, Fl\textbackslash{}"otter\textbackslash{}"od YP, Erdmann J,
  Bieker-Walz L, Behrisch M, eds., \emph{SUMO 2018- Simulating Autonomous and
  Intermodal Transport Systems}, volume~2 of \emph{EPiC Series in Engineering},
  94--110 (EasyChair), ISSN 2516-2330,
  \urlprefix\url{http://dx.doi.org/10.29007/2n4h}.

\bibitem[{Alonso-Mora et~al.(2017)Alonso-Mora, Samaranayake, Wallar, Frazzoli,
  \protect\BIBand{} Rus}]{alonso-mora2017}
Alonso-Mora J, Samaranayake S, Wallar A, Frazzoli E, Rus D (2017) On-demand
  high-capacity ride-sharing via dynamic trip-vehicle assignment.
  \emph{Proceedings of the National Academy of Sciences} 114(3):462--467, ISSN
  0027-8424, \urlprefix\url{http://dx.doi.org/10.1073/pnas.1611675114}.

\bibitem[{Arthur \protect\BIBand{} Vassilvitskii(2007)}]{arthur2007}
Arthur D, Vassilvitskii S (2007) K-means++: The advantages of careful seeding.
  \emph{Proceedings of the Eighteenth Annual ACM-SIAM Symposium on Discrete
  Algorithms}, 1027–1035, SODA '07 (USA: Society for Industrial and Applied
  Mathematics), ISBN 9780898716245.

\bibitem[{Ascheuer et~al.(2000)Ascheuer, Fischetti, \protect\BIBand{}
  Gr{\"o}tschel}]{ascheuer2000}
Ascheuer N, Fischetti M, Gr{\"o}tschel M (2000) A polyhedral study of the
  asymmetric traveling salesman problem with time windows. \emph{Networks}
  36(2):69--79, ISSN 0028-3045,
  \urlprefix\url{http://dx.doi.org/10.1002/1097-0037(200009)36:2<69::AID-NET1>3.0.CO;2-Q}.

\bibitem[{Ascheuer et~al.(2001)Ascheuer, Fischetti, \protect\BIBand{}
  Gr{\"o}tschel}]{ascheuer2001}
Ascheuer N, Fischetti M, Gr{\"o}tschel M (2001) Solving the asymmetric
  travelling salesman problem with time windows by branch-and-cut.
  \emph{Mathematical Programming} 90(3):475--506, ISSN 1436-4646,
  \urlprefix\url{http://dx.doi.org/10.1007/PL00011432}.

\bibitem[{Balas et~al.(1995)Balas, Fischetti, \protect\BIBand{}
  Pulleyblank}]{balas1995}
Balas E, Fischetti M, Pulleyblank WR (1995) The precedence-constrained
  asymmetric traveling salesman polytope. \emph{Mathematical Programming}
  68(1):241--265, ISSN 1436-4646,
  \urlprefix\url{http://dx.doi.org/10.1007/BF01585767}.

\bibitem[{Baldacci et~al.(2004)Baldacci, Maniezzo, \protect\BIBand{}
  Mingozzi}]{baldacci2004}
Baldacci R, Maniezzo V, Mingozzi A (2004) An exact method for the car pooling
  problem based on lagrangean column generation. \emph{Operations Research}
  52(3):422--439, \urlprefix\url{http://dx.doi.org/10.1287/opre.1030.0106}.

\bibitem[{Bard et~al.(2002)Bard, Kontoravdis, \protect\BIBand{} Yu}]{bard2002}
Bard JF, Kontoravdis G, Yu G (2002) A branch-and-cut procedure for the vehicle
  routing problem with time windows. \emph{Transportation Science}
  36(2):250--269, \urlprefix\url{http://dx.doi.org/10.1287/trsc.36.2.250.565}.

\bibitem[{Beasley \protect\BIBand{} Christofides(1989)}]{beasley1989}
Beasley JE, Christofides N (1989) An algorithm for the resource constrained
  shortest path problem. \emph{Networks} 19(4):379--394,
  \urlprefix\url{http://dx.doi.org/10.1002/net.3230190402}.

\bibitem[{Boland et~al.(2006)Boland, Dethridge, \protect\BIBand{}
  Dumitrescu}]{boland2006}
Boland N, Dethridge J, Dumitrescu I (2006) Accelerated label setting algorithms
  for the elementary resource constrained shortest path problem.
  \emph{Operations Research Letters} 34(1):58 -- 68, ISSN 0167-6377,
  \urlprefix\url{http://dx.doi.org/https://doi.org/10.1016/j.orl.2004.11.011}.

\bibitem[{Bornd{\"o}rfer et~al.(2001)Bornd{\"o}rfer, Gr{\"o}tschel,
  \protect\BIBand{} L{\"o}bel}]{borndrfer2001}
Bornd{\"o}rfer R, Gr{\"o}tschel M, L{\"o}bel A (2001) Scheduling duties by
  adaptive column generation. {ZIB-Report 01-02}. Konrad-Zuse-Zentrum f{\"u}r
  Informationstechnik Berlin.

\bibitem[{Cordeau(2006)}]{cordeau2006}
Cordeau JF (2006) A branch-and-cut algorithm for the dial-a-ride problem.
  \emph{Operations Research} 54(3):573--586,
  \urlprefix\url{http://dx.doi.org/10.1287/opre.1060.0283}.

\bibitem[{Cordeau et~al.(2002)Cordeau, Desaulniers, Desrosiers, Solomon,
  \protect\BIBand{} Soumis}]{cordeau2002}
Cordeau JF, Desaulniers G, Desrosiers J, Solomon MM, Soumis F (2002) {VRP} with
  time windows. Toth P, Vigo D, eds., \emph{The Vehicle Routing Problem},
  chapter~7, 157--193 (Philadelphia, PA, USA: SIAM monographs on discrete
  mathematics and applications),
  \urlprefix\url{http://dx.doi.org/10.1137/1.9780898718515.ch7}.

\bibitem[{Cordeau \protect\BIBand{} Laporte(2003{\natexlab{a}})}]{cordeau2003a}
Cordeau JF, Laporte G (2003{\natexlab{a}}) The dial-a-ride problem (darp):
  Variants, modeling issues and algorithms. \emph{Quarterly Journal of the
  Belgian, French and Italian Operations Research Societies} 1(2):89--101, ISSN
  1619-4500, \urlprefix\url{http://dx.doi.org/10.1007/s10288-002-0009-8}.

\bibitem[{Cordeau \protect\BIBand{} Laporte(2003{\natexlab{b}})}]{cordeau2003}
Cordeau JF, Laporte G (2003{\natexlab{b}}) A tabu search heuristic for the
  static multi-vehicle dial-a-ride problem. \emph{Transportation Research Part
  B: Methodological} 37(6):579 -- 594, ISSN 0191-2615,
  \urlprefix\url{http://dx.doi.org/https://doi.org/10.1016/S0191-2615(02)00045-0}.

\bibitem[{Cordeau \protect\BIBand{} Laporte(2007)}]{cordeau2007}
Cordeau JF, Laporte G (2007) The dial-a-ride problem: models and algorithms.
  \emph{Annals of Operations Research} 153(1):29--46, ISSN 1572-9338,
  \urlprefix\url{http://dx.doi.org/10.1007/s10479-007-0170-8}.

\bibitem[{Dantzig et~al.(1954)Dantzig, Fulkerson, \protect\BIBand{}
  Johnson}]{dantzig1954}
Dantzig G, Fulkerson R, Johnson S (1954) Solution of a large-scale
  traveling-salesman problem. \emph{Journal of the Operations Research Society
  of America} 2(4):393--410,
  \urlprefix\url{http://dx.doi.org/10.1287/opre.2.4.393}.

\bibitem[{Dantzig \protect\BIBand{} Wolfe(1960)}]{dantzig1960}
Dantzig GB, Wolfe P (1960) Decomposition principle for linear programs.
  \emph{Operations Research} 8(1):101--111,
  \urlprefix\url{http://dx.doi.org/10.1287/opre.8.1.101}.

\bibitem[{Desaulniers et~al.(2008)Desaulniers, Lessard, \protect\BIBand{}
  Hadjar}]{desaulniers2008}
Desaulniers G, Lessard F, Hadjar A (2008) Tabu search, partial elementarity,
  and generalized k-path inequalities for the vehicle routing problem with time
  windows. \emph{Transportation Science} 42(3):387--404,
  \urlprefix\url{http://dx.doi.org/10.1287/trsc.1070.0223}.

\bibitem[{Desrochers(1988)}]{desrochers1988b}
Desrochers M (1988) An algorithm for the shortest path problem with resource
  constraints. Technical Report G-88-27, Les Cahiers du GERAD, Montreal
  (Quebec), Canada.

\bibitem[{Desrochers et~al.(1992)Desrochers, Desrosiers, \protect\BIBand{}
  Solomon}]{desrochers1992}
Desrochers M, Desrosiers J, Solomon M (1992) A new optimization algorithm for
  the vehicle routing problem with time windows. \emph{Operations Research}
  40(2):342--354, \urlprefix\url{http://dx.doi.org/10.1287/opre.40.2.342}.

\bibitem[{Desrochers \protect\BIBand{} Laporte(1991)}]{desrochers1991}
Desrochers M, Laporte G (1991) Improvements and extensions to the
  miller-tucker-zemlin subtour elimination constraints. \emph{Operations
  Research Letters} 10(1):27 -- 36, ISSN 0167-6377,
  \urlprefix\url{http://dx.doi.org/https://doi.org/10.1016/0167-6377(91)90083-2}.

\bibitem[{Desrosiers et~al.(1984)Desrosiers, Soumis, \protect\BIBand{}
  Desrochers}]{desrosiers1984}
Desrosiers J, Soumis F, Desrochers M (1984) Routing with time windows by column
  generation. \emph{Networks} 14(4):545--565,
  \urlprefix\url{http://dx.doi.org/10.1002/net.3230140406}.

\bibitem[{Dia \protect\BIBand{} Javanshour(2017)}]{dia2017}
Dia H, Javanshour F (2017) Autonomous shared mobility-on-demand: Melbourne
  pilot simulation study. \emph{Transportation Research Procedia} 22:285 --
  296, ISSN 2352-1465,
  \urlprefix\url{http://dx.doi.org/https://doi.org/10.1016/j.trpro.2017.03.035},
  {19th EURO Working Group on Transportation Meeting, EWGT2016, 5-7 September
  2016, Istanbul, Turkey}.

\bibitem[{Drexl(2013)}]{drexl2013}
Drexl M (2013) A note on the separation of subtour elimination constraints in
  elementary shortest path problems. \emph{European Journal of Operational
  Research} 229(3):595 -- 598, ISSN 0377-2217,
  \urlprefix\url{http://dx.doi.org/https://doi.org/10.1016/j.ejor.2013.03.009}.

\bibitem[{Dror(1994)}]{dror1994}
Dror M (1994) Note on the complexity of the shortest path models for column
  generation in vrptw. \emph{Operations Research} 42(5):977--978,
  \urlprefix\url{http://dx.doi.org/10.1287/opre.42.5.977}.

\bibitem[{Dumas et~al.(1991)Dumas, Desrosiers, \protect\BIBand{}
  Soumis}]{dumas1991}
Dumas Y, Desrosiers J, Soumis F (1991) The pickup and delivery problem with
  time windows. \emph{European Journal of Operational Research} 54(1):7 -- 22,
  ISSN 0377-2217,
  \urlprefix\url{http://dx.doi.org/https://doi.org/10.1016/0377-2217(91)90319-Q}.

\bibitem[{Farhan \protect\BIBand{} Chen(2018)}]{farhan2018}
Farhan J, Chen TD (2018) Impact of ridesharing on operational efficiency of
  shared autonomous electric vehicle fleet. \emph{Transportation Research Part
  C: Emerging Technologies} 93:310 -- 321, ISSN 0968-090X,
  \urlprefix\url{http://dx.doi.org/https://doi.org/10.1016/j.trc.2018.04.022}.

\bibitem[{Farley(1990)}]{farley1990}
Farley AA (1990) A note on bounding a class of linear programming problems,
  including cutting stock problems. \emph{Operations Research} 38(5):922--923,
  \urlprefix\url{http://dx.doi.org/10.1287/opre.38.5.922}.

\bibitem[{Firat \protect\BIBand{} Woeginger(2011)}]{firat2011}
Firat M, Woeginger GJ (2011) Analysis of the dial-a-ride problem of hunsaker
  and savelsbergh. \emph{Operations Research Letters} 39(1):32 -- 35, ISSN
  0167-6377,
  \urlprefix\url{http://dx.doi.org/https://doi.org/10.1016/j.orl.2010.11.004}.

\bibitem[{Fischetti \protect\BIBand{} Toth(1997)}]{fischetti1997}
Fischetti M, Toth P (1997) A polyhedral approach to the asymmetric traveling
  salesman problem. \emph{Management Science} 43(11):1520--1536,
  \urlprefix\url{http://dx.doi.org/10.1287/mnsc.43.11.1520}.

\bibitem[{Friedrich(2015)}]{friedrich2015}
Friedrich B (2015) Verkehrliche wirkung autonomer fahrzeuge. Maurer M, Gerdes
  JC, Lenz B, Winner H, eds., \emph{Autonomes Fahren: Technische, rechtliche
  und gesellschaftliche Aspekte}, 331--350 (Berlin, Heidelberg: Springer Berlin
  Heidelberg), ISBN 978-3-662-45854-9,
  \urlprefix\url{http://dx.doi.org/10.1007/978-3-662-45854-9\_16}.

\bibitem[{Gomory \protect\BIBand{} Hu(1961)}]{gomory1961}
Gomory RE, Hu TC (1961) Multi-terminal network flows. \emph{Journal of the
  Society for Industrial and Applied Mathematics} 9(4):551--570,
  \urlprefix\url{http://dx.doi.org/10.1137/0109047}.

\bibitem[{Gouveia \protect\BIBand{} Pires(1999)}]{gouveia1999}
Gouveia L, Pires JM (1999) The asymmetric travelling salesman problem and a
  reformulation of the miller–tucker–zemlin constraints. \emph{European
  Journal of Operational Research} 112(1):134 -- 146, ISSN 0377-2217,
  \urlprefix\url{http://dx.doi.org/https://doi.org/10.1016/S0377-2217(97)00358-5}.

\bibitem[{Gr{\"o}tschel \protect\BIBand{} Padberg(1985)}]{grotshcel1985}
Gr{\"o}tschel M, Padberg M (1985) Polyhedral theory. Lawler E, Lenstra J,
  Rinnooy~Kan A, Shmoys D, eds., \emph{The Traveling Salesman Problem},
  chapter~8, 251--305, A Wiley-Interscience publication (John Wiley \& Sons,
  Incorporated), ISBN 9780471904137,
  \urlprefix\url{https://books.google.com/books?id=EPFQAAAAMAAJ}.

\bibitem[{Gr{\"o}tschel \protect\BIBand{} Padberg(1975)}]{grotshcel1975}
Gr{\"o}tschel M, Padberg MW (1975) Partial linear characterizations of the
  asymmetric travelling salesman polytope. \emph{Mathematical Programming}
  8(1):378--381, ISSN 1436-4646,
  \urlprefix\url{http://dx.doi.org/10.1007/BF01580454}.

\bibitem[{Gschwind \protect\BIBand{} Irnich(2015)}]{gschwind2015}
Gschwind T, Irnich S (2015) Effective handling of dynamic time windows and its
  application to solving the dial-a-ride problem. \emph{Transportation Science}
  49(2):335--354, \urlprefix\url{http://dx.doi.org/10.1287/trsc.2014.0531}.

\bibitem[{Hasan \protect\BIBand{} Van~Hentenryck(2020)}]{hasan2020b}
Hasan MH, Van~Hentenryck P (2020) The flexible and real-time commute trip
  sharing problems. \emph{Constraints} 25(3):160--179, ISSN 1572-9354,
  \urlprefix\url{http://dx.doi.org/10.1007/s10601-020-09310-5}.

\bibitem[{Hasan \protect\BIBand{} Van~Hentenryck(\noop{3001}in press
  2021)}]{hasan2021}
Hasan MH, Van~Hentenryck P (\noop{3001}in press 2021) The benefits of
  autonomous vehicles for community-based trip sharing. \emph{Transportation
  Research Part C: Emerging Technologies} .

\bibitem[{Hasan et~al.(2018)Hasan, Van~Hentenryck, Budak, Chen,
  \protect\BIBand{} Chaudhry}]{hasan2018}
Hasan MH, Van~Hentenryck P, Budak C, Chen J, Chaudhry C (2018) Community-based
  trip sharing for urban commuting. McIlraith S, Weinberger K, eds.,
  \emph{Proceedings of the Thirty-Second AAAI Conference on Artificial
  Intelligence}, 6589--6597, AAAI-18 (Palo Alto, California, USA: AAAI Press).

\bibitem[{Hasan et~al.(2020)Hasan, Van~Hentenryck, \protect\BIBand{}
  Legrain}]{hasan2020}
Hasan MH, Van~Hentenryck P, Legrain A (2020) The commute trip-sharing problem.
  \emph{Transportation Science} 54(6):1640--1675,
  \urlprefix\url{http://dx.doi.org/10.1287/trsc.2019.0969}.

\bibitem[{Haugland \protect\BIBand{} Ho(2010)}]{haughland2010}
Haugland D, Ho SC (2010) Feasibility testing for dial-a-ride problems. Chen B,
  ed., \emph{Algorithmic Aspects in Information and Management}, 170--179
  (Berlin, Heidelberg: Springer Berlin Heidelberg), ISBN 978-3-642-14355-7.

\bibitem[{Hunsaker \protect\BIBand{} Savelsbergh(2002)}]{hunsaker2002}
Hunsaker B, Savelsbergh M (2002) Efficient feasibility testing for dial-a-ride
  problems. \emph{Operations Research Letters} 30(3):169 -- 173, ISSN
  0167-6377,
  \urlprefix\url{http://dx.doi.org/https://doi.org/10.1016/S0167-6377(02)00120-7}.

\bibitem[{Irnich \protect\BIBand{} Desaulniers(2005)}]{irnich2005}
Irnich S, Desaulniers G (2005) Shortest path problems with resource
  constraints. Desaulniers G, Desrosiers J, Solomon MM, eds., \emph{Column
  Generation}, 33--65 (Boston, MA: Springer US), ISBN 978-0-387-25486-9,
  \urlprefix\url{http://dx.doi.org/10.1007/0-387-25486-2\_2}.

\bibitem[{Jaw et~al.(1986)Jaw, Odoni, Psaraftis, \protect\BIBand{}
  Wilson}]{jaw1986}
Jaw JJ, Odoni AR, Psaraftis HN, Wilson NH (1986) A heuristic algorithm for the
  multi-vehicle advance request dial-a-ride problem with time windows.
  \emph{Transportation Research Part B: Methodological} 20(3):243 -- 257, ISSN
  0191-2615,
  \urlprefix\url{http://dx.doi.org/https://doi.org/10.1016/0191-2615(86)90020-2}.

\bibitem[{Kallehauge et~al.(2007)Kallehauge, Boland, \protect\BIBand{}
  Madsen}]{kallehauge2007}
Kallehauge B, Boland N, Madsen OB (2007) Path inequalities for the vehicle
  routing problem with time windows. \emph{Networks} 49(4):273--293,
  \urlprefix\url{http://dx.doi.org/10.1002/net.20178}.

\bibitem[{Kohl et~al.(1999)Kohl, Desrosiers, Madsen, Solomon, \protect\BIBand{}
  Soumis}]{kohl1999}
Kohl N, Desrosiers J, Madsen OBG, Solomon MM, Soumis F (1999) 2-path cuts for
  the vehicle routing problem with time windows. \emph{Transportation Science}
  33(1):101--116, \urlprefix\url{http://dx.doi.org/10.1287/trsc.33.1.101}.

\bibitem[{Langevin et~al.(1990)Langevin, Soumis, \protect\BIBand{}
  Desrosiers}]{langevin1990}
Langevin A, Soumis F, Desrosiers J (1990) Classification of travelling salesman
  problem formulations. \emph{Operations Research Letters} 9(2):127 -- 132,
  ISSN 0167-6377,
  \urlprefix\url{http://dx.doi.org/https://doi.org/10.1016/0167-6377(90)90052-7}.

\bibitem[{Liberti(2004)}]{Liberti2004}
Liberti L (2004) Reduction constraints for the global optimization of nlps.
  \emph{International Transactions in Operational Research} 11(1):33--41,
  \urlprefix\url{http://dx.doi.org/https://doi.org/10.1111/j.1475-3995.2004.00438.x}.

\bibitem[{{Lloyd}(1982)}]{lloyd1982}
{Lloyd} S (1982) Least squares quantization in pcm. \emph{IEEE Transactions on
  Information Theory} 28(2):129--137, ISSN 1557-9654,
  \urlprefix\url{http://dx.doi.org/10.1109/TIT.1982.1056489}.

\bibitem[{Ma et~al.(2017)Ma, Li, Zhou, \protect\BIBand{} Hao}]{ma2017}
Ma J, Li X, Zhou F, Hao W (2017) Designing optimal autonomous vehicle sharing
  and reservation systems: A linear programming approach. \emph{Transportation
  Research Part C: Emerging Technologies} 84:124 -- 141, ISSN 0968-090X,
  \urlprefix\url{http://dx.doi.org/https://doi.org/10.1016/j.trc.2017.08.022}.

\bibitem[{Martinez \protect\BIBand{} Viegas(2017)}]{martinez2017}
Martinez LM, Viegas JM (2017) Assessing the impacts of deploying a shared
  self-driving urban mobility system: An agent-based model applied to the city
  of lisbon, portugal. \emph{International Journal of Transportation Science
  and Technology} 6(1):13 -- 27, ISSN 2046-0430,
  \urlprefix\url{http://dx.doi.org/https://doi.org/10.1016/j.ijtst.2017.05.005},
  connected and Automated Vehicles: Effects on Traffic, Mobility and Urban
  Design.

\bibitem[{{Mena-Oreja} et~al.(2018){Mena-Oreja}, {Gozalvez}, \protect\BIBand{}
  {Sepulcre}}]{menaoreja2018}
{Mena-Oreja} J, {Gozalvez} J, {Sepulcre} M (2018) Effect of the configuration
  of platooning maneuvers on the traffic flow under mixed traffic scenarios.
  \emph{2018 IEEE Vehicular Networking Conference (VNC)}, 1--4, ISSN 2157-9865,
  \urlprefix\url{http://dx.doi.org/10.1109/VNC.2018.8628381}.

\bibitem[{Miller et~al.(1960)Miller, Tucker, \protect\BIBand{}
  Zemlin}]{miller1960}
Miller CE, Tucker AW, Zemlin RA (1960) Integer programming formulation of
  traveling salesman problems. \emph{J. ACM} 7(4):326–329, ISSN 0004-5411,
  \urlprefix\url{http://dx.doi.org/10.1145/321043.321046}.

\bibitem[{Mourad et~al.(2019)Mourad, Puchinger, \protect\BIBand{}
  Chu}]{mourad2019}
Mourad A, Puchinger J, Chu C (2019) A survey of models and algorithms for
  optimizing shared mobility. \emph{Transportation Research Part B:
  Methodological} 123:323 -- 346, ISSN 0191-2615,
  \urlprefix\url{http://dx.doi.org/https://doi.org/10.1016/j.trb.2019.02.003}.

\bibitem[{Naddef \protect\BIBand{} Rinaldi(2001)}]{naddef2001}
Naddef D, Rinaldi G (2001) Branch-and-cut algorithms for the capacitated vrp.
  \emph{The Vehicle Routing Problem}, 53–84 (USA: Society for Industrial and
  Applied Mathematics), ISBN 0898714982.

\bibitem[{Narayanan et~al.(2020)Narayanan, Chaniotakis, \protect\BIBand{}
  Antoniou}]{narayanan2020}
Narayanan S, Chaniotakis E, Antoniou C (2020) Shared autonomous vehicle
  services: A comprehensive review. \emph{Transportation Research Part C:
  Emerging Technologies} 111:255 -- 293, ISSN 0968-090X,
  \urlprefix\url{http://dx.doi.org/https://doi.org/10.1016/j.trc.2019.12.008}.

\bibitem[{{NYC Taxi \& Limousine Commission}(2020)}]{nyctriprecord}
{NYC Taxi \& Limousine Commission} (2020) {TLC} trip record data.
  \url{https://www1.nyc.gov/site/tlc/about/tlc-trip-record-data.page},
  accessed: 2020-11-20.

\bibitem[{Olia et~al.(2018)Olia, Razavi, Abdulhai, \protect\BIBand{}
  Abdelgawad}]{olia2018}
Olia A, Razavi S, Abdulhai B, Abdelgawad H (2018) Traffic capacity implications
  of automated vehicles mixed with regular vehicles. \emph{Journal of
  Intelligent Transportation Systems} 22(3):244--262,
  \urlprefix\url{http://dx.doi.org/10.1080/15472450.2017.1404680}.

\bibitem[{Padberg \protect\BIBand{} Rinaldi(1990)}]{padberg1990}
Padberg M, Rinaldi G (1990) An efficient algorithm for the minimum capacity cut
  problem. \emph{Mathematical Programming} 47(1):19--36, ISSN 1436-4646,
  \urlprefix\url{http://dx.doi.org/10.1007/BF01580850}.

\bibitem[{Padberg \protect\BIBand{} Rinaldi(1991)}]{padberg1991}
Padberg M, Rinaldi G (1991) A branch-and-cut algorithm for the resolution of
  large-scale symmetric traveling salesman problems. \emph{SIAM Review}
  33(1):60--100, \urlprefix\url{http://dx.doi.org/10.1137/1033004}.

\bibitem[{Ropke \protect\BIBand{} Cordeau(2006)}]{ropke2006}
Ropke S, Cordeau JF (2006) \emph{Heuristic and exact algorithms for vehicle
  routing problems}. Ph.D. thesis, University of Copenhagen,
  branch-and-cut-and-price for the pickup and delivery problem with time
  windows.

\bibitem[{Ropke \protect\BIBand{} Cordeau(2009)}]{ropke2009}
Ropke S, Cordeau JF (2009) Branch and cut and price for the pickup and delivery
  problem with time windows. \emph{Transportation Science} 43(3):267--286,
  \urlprefix\url{http://dx.doi.org/10.1287/trsc.1090.0272}.

\bibitem[{Rousseau et~al.(2007)Rousseau, Gendreau, \protect\BIBand{}
  Feillet}]{rousseau2007}
Rousseau LM, Gendreau M, Feillet D (2007) Interior point stabilization for
  column generation. \emph{Operations Research Letters} 35(5):660 -- 668, ISSN
  0167-6377,
  \urlprefix\url{http://dx.doi.org/https://doi.org/10.1016/j.orl.2006.11.004}.

\bibitem[{Rousseau et~al.(2004)Rousseau, Gendreau, Pesant, \protect\BIBand{}
  Focacci}]{rousseau2004}
Rousseau LM, Gendreau M, Pesant G, Focacci F (2004) Solving vrptws with
  constraint programming based column generation. \emph{Annals of Operations
  Research} 130(1):199--216, ISSN 1572-9338,
  \urlprefix\url{http://dx.doi.org/10.1023/B:ANOR.0000032576.73681.29}.

\bibitem[{Ruiz \protect\BIBand{} Grossmann(2011)}]{Ruiz2011}
Ruiz JP, Grossmann IE (2011) Using redundancy to strengthen the relaxation for
  the global optimization of minlp problems. \emph{Computers \& Chemical
  Engineering} 35(12):2729 -- 2740, ISSN 0098-1354,
  \urlprefix\url{http://dx.doi.org/https://doi.org/10.1016/j.compchemeng.2011.01.035}.

\bibitem[{Ruland \protect\BIBand{} Rodin(1997)}]{ruland1997}
Ruland K, Rodin E (1997) The pickup and delivery problem: Faces and
  branch-and-cut algorithm. \emph{Computers \& Mathematics with Applications}
  33(12):1 -- 13, ISSN 0898-1221,
  \urlprefix\url{http://dx.doi.org/https://doi.org/10.1016/S0898-1221(97)00090-4}.

\bibitem[{{Salazar} et~al.(2018){Salazar}, {Rossi}, {Schiffer}, {Onder},
  \protect\BIBand{} {Pavone}}]{salazar2018}
{Salazar} M, {Rossi} F, {Schiffer} M, {Onder} CH, {Pavone} M (2018) On the
  interaction between autonomous mobility-on-demand and public transportation
  systems. \emph{2018 21st International Conference on Intelligent
  Transportation Systems (ITSC)}, 2262--2269, ISSN 2153-0017,
  \urlprefix\url{http://dx.doi.org/10.1109/ITSC.2018.8569381}.

\bibitem[{Santi et~al.(2014)Santi, Resta, Szell, Sobolevsky, Strogatz,
  \protect\BIBand{} Ratti}]{santi2014}
Santi P, Resta G, Szell M, Sobolevsky S, Strogatz SH, Ratti C (2014)
  Quantifying the benefits of vehicle pooling with shareability networks.
  \emph{Proceedings of the National Academy of Sciences} 111(37):13290--13294,
  ISSN 0027-8424, \urlprefix\url{http://dx.doi.org/10.1073/pnas.1403657111}.

\bibitem[{Savelsbergh(1985)}]{savelsbergh1985}
Savelsbergh MWP (1985) Local search in routing problems with time windows.
  \emph{Annals of Operations Research} 4(1):285--305, ISSN 1572-9338,
  \urlprefix\url{http://dx.doi.org/10.1007/BF02022044}.

\bibitem[{Talebpour \protect\BIBand{} Mahmassani(2016)}]{talebpour2016}
Talebpour A, Mahmassani HS (2016) Influence of connected and autonomous
  vehicles on traffic flow stability and throughput. \emph{Transportation
  Research Part C: Emerging Technologies} 71:143 -- 163, ISSN 0968-090X,
  \urlprefix\url{http://dx.doi.org/https://doi.org/10.1016/j.trc.2016.07.007}.

\bibitem[{Tang et~al.(2010)Tang, Kong, Lau, \protect\BIBand{} Ip}]{tang2010}
Tang J, Kong Y, Lau H, Ip AW (2010) A note on “efficient feasibility testing
  for dial-a-ride problems”. \emph{Operations Research Letters} 38(5):405 --
  407, ISSN 0167-6377,
  \urlprefix\url{http://dx.doi.org/https://doi.org/10.1016/j.orl.2010.05.002}.

\bibitem[{Tarjan(1972)}]{tarjan1972}
Tarjan R (1972) Depth-first search and linear graph algorithms. \emph{SIAM
  Journal on Computing} 1(2):146--160,
  \urlprefix\url{http://dx.doi.org/10.1137/0201010}.

\bibitem[{{Tientrakool} et~al.(2011){Tientrakool}, {Ho}, \protect\BIBand{}
  {Maxemchuk}}]{tientrakool2011}
{Tientrakool} P, {Ho} Y, {Maxemchuk} NF (2011) Highway capacity benefits from
  using vehicle-to-vehicle communication and sensors for collision avoidance.
  \emph{2011 IEEE Vehicular Technology Conference (VTC Fall)}, 1--5, ISSN
  1090-3038, \urlprefix\url{http://dx.doi.org/10.1109/VETECF.2011.6093130}.

\bibitem[{Zhang \protect\BIBand{} Guhathakurta(2017)}]{zhang2017}
Zhang W, Guhathakurta S (2017) Parking spaces in the age of shared autonomous
  vehicles: How much parking will we need and where? \emph{Transportation
  Research Record} 2651(1):80--91,
  \urlprefix\url{http://dx.doi.org/10.3141/2651-09}.

\bibitem[{Zhang et~al.(2015)Zhang, Guhathakurta, Fang, \protect\BIBand{}
  Zhang}]{zhang2015}
Zhang W, Guhathakurta S, Fang J, Zhang G (2015) Exploring the impact of shared
  autonomous vehicles on urban parking demand: An agent-based simulation
  approach. \emph{Sustainable Cities and Society} 19:34 -- 45, ISSN 2210-6707,
  \urlprefix\url{http://dx.doi.org/https://doi.org/10.1016/j.scs.2015.07.006}.

\end{thebibliography}



\bibliographystyle{informs2014}

\newcommand{\noop}[1]{}

\newpage

\APPENDIX{Filtering of Graph $\mathcal{G}$}

Graph $\mathcal{G}$ can be made more compact by only retaining edges
that satisfy a priori route-feasibility constraints. This is done by
pre-processing time-window, pairing, precedence, and ride-duration
limit constraints on $\mathcal{A}$ to identify and eliminate edges
that are infeasible, i.e., those that cannot belong to any feasible AV
route. In this work, the set of infeasible edges is identified using a
combination of rules proposed by \cite{dumas1991} and
\cite{cordeau2006}. These rules are presented in the Appendix.

\begin{enumerate}[(a)]
	\item Direct trips to and from the depot:
	\begin{itemize}
		\item $\{(v_s,v_t), (v_t,v_s)\}$
		\item $\{(i,v_s), (i,v_t), (v_t,i) : i\in\mathcal{P}\}$
		\item $\{(v_s,i), (i,v_s), (v_t,i) : i\in\mathcal{D}\}$
	\end{itemize}
	\item Precedence of pickup and drop-off nodes of inbound and outbound trips of each commuter (constraints \eqref{eqn:precedence}): $\{(i,2n+i), (i,3n+i), (n+i,i), (n+i,3n+i), (2n+i,i), (2n+i,n+i), (3n+i,i), (3n+i,n+i), (3n+i,2n+i) : i\in\mathcal{P}^+\}$
	\item Precedence of pickup and drop-off nodes of inbound and outbound mini routes:
	\begin{itemize}
		\item $\{(i,j):i\in\mathcal{P}^+ \wedge j\in\mathcal{P}^-\cup\mathcal{D}^-\}$
		\item $\{(i,j):i\in\mathcal{D}^+ \wedge j\in\mathcal{D}^-\}$
		\item $\{(i,j):i\in\mathcal{P}^- \wedge j\in\mathcal{P}^+\cup\mathcal{D}^+\}$
		\item $\{(i,j):i\in\mathcal{D}^- \wedge j\in\mathcal{D}^+\}$
	\end{itemize}
	\item Time windows along each edge: $\{(i,j):(i,j)\in\mathcal{A} \setminus \{\delta^+(v_s)\cup\delta^-(v_t)\} \wedge a_i + s_i + \tau_{(i,j)} > b_j\}$
	\item Ride-duration limit of each commuter: $\{(i,j),(j,n+i):i\in\mathcal{P} \wedge j\in\mathcal{P}\cup\mathcal{D} \wedge i\neq j \wedge \tau_{(i,j)} + s_j + \tau_{(j,n+i)} > L_i\}$
	\item Time windows and ride-duration limits of pairs of trips:
	\begin{itemize}
		\item $\{(i,n+j):i,j\in\mathcal{P} \wedge i\neq j \wedge \neg feasible(j\rightarrow i\rightarrow n+j\rightarrow n+i)\}$
		\item $\{(n+i,j):i,j\in\mathcal{P} \wedge i\neq j \wedge \neg feasible(i\rightarrow n+i\rightarrow j\rightarrow n+j)\}$
		\item $\{(i,j):i,j\in\mathcal{P} \wedge i\neq j \wedge \neg feasible(i\rightarrow j\rightarrow n+i\rightarrow n+j) \wedge \neg feasible(i\rightarrow j\rightarrow n+j\rightarrow n+i)\}$
		\item $\{(n+i,n+j):i,j\in\mathcal{P} \wedge i\neq j \wedge \neg feasible(i\rightarrow j\rightarrow n+i\rightarrow n+j) \wedge \neg feasible(j\rightarrow i\rightarrow n+i\rightarrow n+j)\}$
	\end{itemize}
\end{enumerate}

\begin{figure}[!t]
	\centering
	\includegraphics[width=0.5\linewidth]{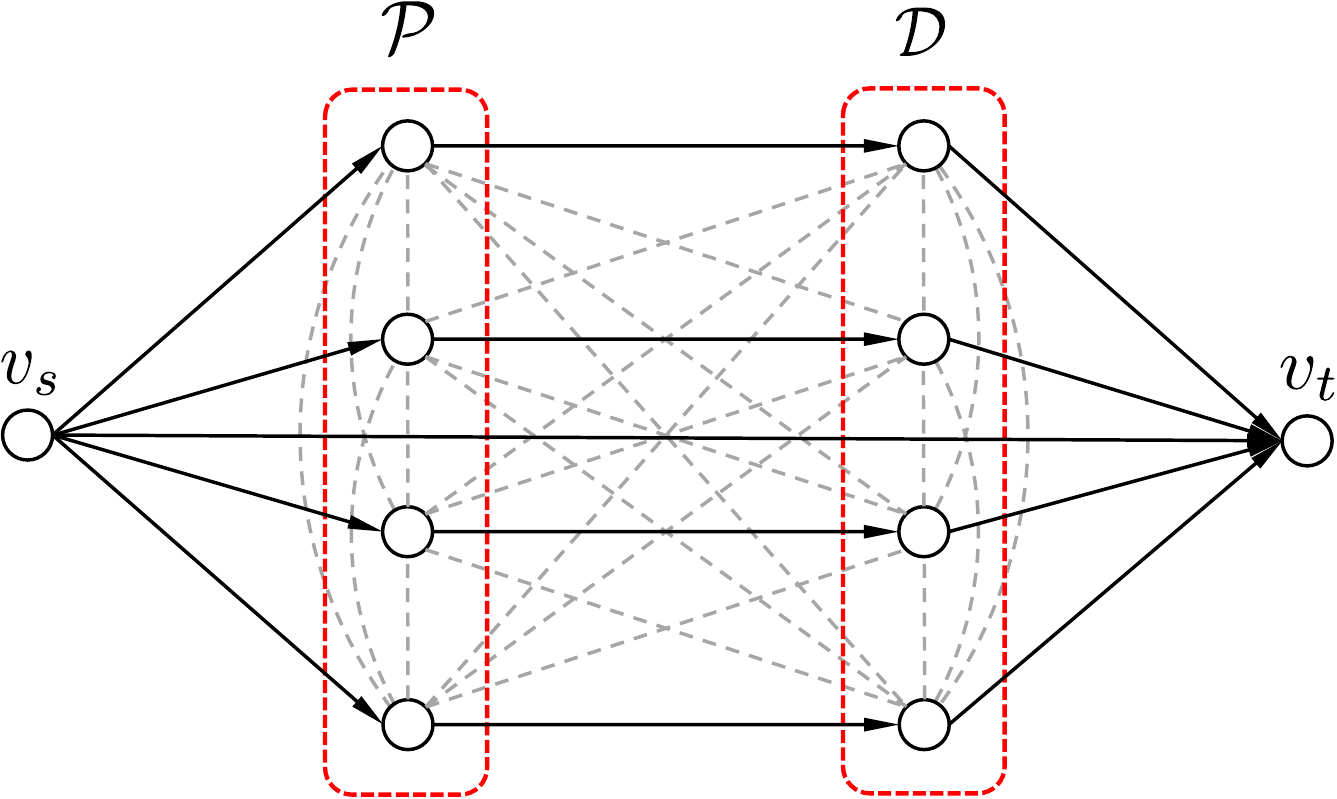}
	\caption{Graph $\mathcal{G}$ (Each Dotted Line Represents a Pair of Bidirectional Edges).}
	\label{fig:graph1}
\end{figure}

\noindent
Note that the sets of edges in (f) utilize the $feasible$ function to
determine if a partial route satisfies time-window and ride-duration
limit constraints. For instance, the first condition indicates that edge
$(i,n+j)$ is infeasible if the partial route $j\rightarrow i\rightarrow
n+j\rightarrow n+i$ is infeasible. Figure \ref{fig:graph1} illustrates
an example of graph $\mathcal{G}$ resulting from this pre-processing step.

\newpage

\APPENDIX{Computational Results}

Table \ref{tab:my-table} summarizes the results of
CTSPAV\textsubscript{Hybrid} for every large problem instance. Its
first column shows the name of every instance. The next three columns
display properties that characterize the size of each instance. They
list the node count of graph $\mathcal{G}$, $|\mathcal{N}|$, the edge
count of the graph (after the pre-processing step), $|\mathcal{A}|$,
and finally the number of mini routes generated by the MREA,
$|\mathrm{\Omega}|$, for every instance. The next column shows the
wall time spent to the enumerate the mini routes. The remaining
columns summarize the results of CTSPAV\textsubscript{Hybrid}. The
first two show the vehicle count and total travel distance from its
best incumbent solution. The next two display the absolute gap for the
vehicle count and the optimality gap for the objective value of the
best incumbent solution.  The following column shows the number of
tree nodes explored in the solution process. The last two columns
display the (total) wall time spent to solve the MIP and that spent to
close the vehicle count gap. For the very last column, values are only
listed for instances whereby the vehicle count gap could be closed
within the 2-hour time limit. It is left blank otherwise. Tables
\ref{tab:my-table3} and \ref{tab:my-table5} provide the same set of
information for CTSPAV\textsubscript{Hybrid} for every medium and
tight problem instance respectively. On the other hand, Tables
\ref{tab:my-table2}, \ref{tab:my-table4}, and \ref{tab:my-table6} show
the results of CTSPAV\textsubscript{SEC} and
CTSPAV\textsubscript{Base} for all large, medium, and tight problem
instances respectively.

\begin{table}[!t]
	\centering
	\caption{Results of CTSPAV\textsubscript{Hybrid} for the Large Problem Instances}
	\label{tab:my-table}
	\setlength{\tabcolsep}{3pt}
	\resizebox{\textwidth}{!}{%
		\begin{tabular}{cccccccccccc}
			\hline\up
			\multirow{2}{*}{\begin{tabular}[c]{@{}c@{}}Instance\\ name\end{tabular}} & \multirow{2}{*}{\begin{tabular}[c]{@{}c@{}}Node\\ count\end{tabular}} & \multirow{2}{*}{\begin{tabular}[c]{@{}c@{}}Edge\\ count\end{tabular}} & \multirow{2}{*}{\begin{tabular}[c]{@{}c@{}}Mini\\ route\\ count\end{tabular}} & \multirow{2}{*}{\begin{tabular}[c]{@{}c@{}}Route\\ enumeration\\ time (s)\end{tabular}} & \multirow{2}{*}{\begin{tabular}[c]{@{}c@{}}Vehicle\\ count\end{tabular}} & \multirow{2}{*}{\begin{tabular}[c]{@{}c@{}}Total\\ distance\\ (m)\end{tabular}} & \multirow{2}{*}{\begin{tabular}[c]{@{}c@{}}Vehicle\\ count\\ gap\end{tabular}} & \multirow{2}{*}{\begin{tabular}[c]{@{}c@{}}Optimality\\ gap (\%)\end{tabular}} & \multirow{2}{*}{\begin{tabular}[c]{@{}c@{}}Nodes\\ explored\end{tabular}} & \multicolumn{2}{c}{Wall time (s)\down}\\
			\cline{11-12}\rule{0pt}{20pt} 
			&  &  &  &  &  &  &  &  &  & MIP & \begin{tabular}[c]{@{}c@{}}Optimal\\ count\down\end{tabular}\\
			\hline\up
			L0 & 402 & 23983 & 3730 & 22 & 3 & 642049 & 0 & 0.0 & 156016 & 5360 & 1284 \\
			L1 & 402 & 22621 & 1093 & 21 & 3 & 463065 & 1 & 33.3 & 524584 & 7200 & - \\
			L2 & 402 & 26781 & 51175 & 24 & 4 & 817348 & 2 & 49.9 & 6424 & 7200 & - \\
			L3 & 402 & 26496 & 63597 & 24 & 4 & 841180 & 2 & 49.9 & 7430 & 7202 & - \\
			L4 & 402 & 25309 & 49147 & 23 & 4 & 813018 & 1 & 24.9 & 11734 & 7201 & - \\
			L5 & 402 & 22425 & 1605 & 20 & 3 & 512675 & 1 & 33.3 & 189596 & 7200 & - \\
			L6 & 402 & 26420 & 20060 & 23 & 4 & 955285 & 2 & 49.9 & 7935 & 7201 & - \\
			L7 & 402 & 24699 & 21403 & 23 & 4 & 888490 & 1 & 24.9 & 22067 & 7201 & - \\
			L8 & 402 & 25710 & 14818 & 23 & 4 & 844674 & 1 & 24.9 & 23822 & 7200 & - \\
			L9 & 402 & 27315 & 191067 & 25 & 5 & 737361 & 3 & 59.9 & 1511 & 7200 & - \\
			L10 & 402 & 24386 & 5807 & 25 & 3 & 555102 & 1 & 33.3 & 30016 & 7201 & - \\
			L11 & 402 & 25639 & 18237 & 23 & 3 & 570036 & 1 & 33.3 & 13176 & 7201 & - \\
			L12 & 402 & 23748 & 3631 & 21 & 3 & 581863 & 1 & 33.3 & 125059 & 7200 & - \\
			L13 & 402 & 24581 & 6835 & 24 & 3 & 624843 & 1 & 33.3 & 23394 & 7202 & - \\
			L14 & 402 & 26287 & 72200 & 23 & 4 & 949361 & 2 & 49.9 & 5138 & 7201 & - \\
			L15 & 402 & 24898 & 114817 & 38 & 4 & 1108007 & 2 & 49.9 & 7258 & 7200 & - \\
			L16 & 402 & 24203 & 9231 & 22 & 4 & 847394 & 1 & 24.9 & 75500 & 7200 & - \\
			L17 & 402 & 23734 & 6404 & 22 & 4 & 863265 & 0 & 0.0 & 22485 & 7200 & 5883 \\
			L18 & 402 & 24712 & 4417 & 33 & 4 & 914762 & 1 & 24.9 & 33188 & 7201 & - \\
			L19 & 402 & 25513 & 35873 & 24 & 3 & 698599 & 1 & 33.3 & 11984 & 7201 & - \\
			L20 & 402 & 25528 & 58833 & 23 & 3 & 779684 & 1 & 33.3 & 8639 & 7200 & - \\
			L21 & 402 & 22832 & 4870 & 21 & 2 & 457911 & 0 & 0.0 & 166142 & 7200 & 2217\down\\ 
			\hline
		\end{tabular}%
	}
\end{table}

\begin{table}[!th]
	\centering
	\caption{Results of CTSPAV\textsubscript{SEC} and CTSPAV\textsubscript{Base} for the Large Problem Instances}
	\label{tab:my-table2}
	\setlength{\tabcolsep}{2pt}
	\resizebox{\textwidth}{!}{%
		\begin{tabular}{cccccccc|ccccccc}
			\hline\up
			\multirow{4}{*}{\begin{tabular}[c]{@{}c@{}}Instance\\ name\end{tabular}} & \multicolumn{14}{c}{CTSPAV variant\down} \\
			\cline{2-15} \rule{0pt}{12pt}
			& \multicolumn{7}{c|}{SEC} & \multicolumn{7}{c}{Base\down} \\
			\cline{2-15} \rule{0pt}{14pt}
			& \multirow{2}{*}{\begin{tabular}[c]{@{}c@{}}Vehicle\\ count\end{tabular}} & \multirow{2}{*}{\begin{tabular}[c]{@{}c@{}}Total\\ distance\\ (m)\end{tabular}} & \multirow{2}{*}{\begin{tabular}[c]{@{}c@{}}Vehicle\\ count\\ gap\end{tabular}} & \multirow{2}{*}{\begin{tabular}[c]{@{}c@{}}Optimality\\ gap (\%)\end{tabular}} & \multirow{2}{*}{\begin{tabular}[c]{@{}c@{}}Nodes\\ explored\end{tabular}} & \multicolumn{2}{c|}{Wall time (s)} & \multirow{2}{*}{\begin{tabular}[c]{@{}c@{}}Vehicle\\ count\end{tabular}} & \multirow{2}{*}{\begin{tabular}[c]{@{}c@{}}Total\\ distance\\ (m)\end{tabular}} & \multirow{2}{*}{\begin{tabular}[c]{@{}c@{}}Vehicle\\ count\\ gap\end{tabular}} & \multirow{2}{*}{\begin{tabular}[c]{@{}c@{}}Optimality\\ gap (\%)\end{tabular}} & \multirow{2}{*}{\begin{tabular}[c]{@{}c@{}}Nodes\\ explored\end{tabular}} & \multicolumn{2}{c}{Wall time (s)\down} \\
			\cline{7-8} \cline{14-15} \rule{0pt}{19pt}
			&  &  &  &  &  & MIP & \begin{tabular}[c]{@{}c@{}}Optimal\\ count\end{tabular} &  &  &  &  &  & MIP & \begin{tabular}[c]{@{}c@{}}Optimal\\ count\down\end{tabular} \\
			\hline\up
			L0  & 3 & 646884  & 1 & 33.3 & 43103  & 7200 & - & 3 & 652906  & 2 & 66.5 & 24638  & 7201 & - \\
			L1  & 3 & 463065  & 1 & 33.3 & 135613 & 7228 & - & 3 & 463065  & 2 & 66.6 & 408157 & 7202 & - \\
			L2  & 4 & 821989  & 2 & 49.9 & 6369   & 7218 & - & 4 & 824321  & 3 & 74.8 & 5229   & 7201 & - \\
			L3  & 4 & 849844  & 2 & 49.9 & 4713   & 7215 & - & 4 & 843208  & 3 & 74.8 & 5291   & 7200 & - \\
			L4  & 5 & 820800  & 3 & 59.9 & 10005  & 7202 & - & 5 & 831319  & 3 & 59.9 & 20952  & 7201 & - \\
			L5  & 3 & 512838  & 1 & 33.3 & 73463  & 7202 & - & 3 & 512675  & 2 & 66.5 & 195089 & 7201 & - \\
			L6  & 4 & 971911  & 2 & 49.9 & 9541   & 7207 & - & 4 & 967746  & 3 & 74.8 & 11540  & 7204 & - \\
			L7  & 4 & 891808  & 2 & 49.9 & 7244   & 7206 & - & 4 & 893550  & 3 & 74.8 & 15275  & 7201 & - \\
			L8  & 4 & 845333  & 2 & 49.9 & 8301   & 7201 & - & 4 & 845100  & 3 & 74.8 & 16814  & 7200 & - \\
			L9  & 5 & 730915  & 3 & 59.9 & 2023   & 7200 & - & 5 & 720023  & 4 & 79.9 & 1906   & 7200 & - \\
			L10 & 3 & 555102  & 1 & 33.3 & 21162  & 7200 & - & 3 & 555102  & 2 & 66.5 & 21223  & 7201 & - \\
			L11 & 3 & 573246  & 1 & 33.3 & 3428   & 7203 & - & 3 & 574227  & 2 & 66.5 & 21195  & 7200 & - \\
			L12 & 3 & 581863  & 1 & 33.3 & 34193  & 7200 & - & 3 & 581863  & 2 & 66.5 & 57588  & 7202 & - \\
			L13 & 3 & 626100  & 1 & 33.3 & 15871  & 7221 & - & 3 & 625042  & 2 & 66.5 & 36251  & 7201 & - \\
			L14 & 4 & 949659  & 2 & 49.9 & 5431   & 7213 & - & 4 & 932389  & 3 & 74.8 & 4986   & 7200 & - \\
			L15 & 4 & 1108620 & 2 & 49.9 & 4435   & 7202 & - & 4 & 1116187 & 3 & 74.8 & 2732   & 7201 & - \\
			L16 & 4 & 857161  & 2 & 49.9 & 12595  & 7203 & - & 4 & 846684  & 3 & 74.8 & 21489  & 7200 & - \\
			L17 & 4 & 867674  & 2 & 49.9 & 21259  & 7201 & - & 4 & 865011  & 2 & 49.9 & 21691  & 7200 & - \\
			L18 & 4 & 917395  & 2 & 49.9 & 18251  & 7201 & - & 4 & 914762  & 2 & 49.9 & 20825  & 7200 & - \\
			L19 & 4 & 697540  & 2 & 49.9 & 4925   & 7298 & - & 4 & 706887  & 3 & 74.9 & 15757  & 7200 & - \\
			L20 & 3 & 772418  & 1 & 33.3 & 6277   & 7318 & - & 3 & 778248  & 2 & 66.5 & 7573   & 7200 & - \\
			L21 & 3 & 447435  & 2 & 66.6 & 1632   & 7259 & - & 2 & 458460  & 1 & 49.9 & 86453  & 7205 & -\down\\ 
			\hline
		\end{tabular}%
	}
\end{table}

\begin{table}[!th]
	\centering
	\caption{Results of CTSPAV\textsubscript{Hybrid} for the Medium Problem Instances}
	\label{tab:my-table3}
	\setlength{\tabcolsep}{3pt}
	\resizebox{\textwidth}{!}{%
		\begin{tabular}{cccccccccccc}
			\hline\up
			\multirow{2}{*}{\begin{tabular}[c]{@{}c@{}}Instance\\ name\end{tabular}} & \multirow{2}{*}{\begin{tabular}[c]{@{}c@{}}Node\\ count\end{tabular}} & \multirow{2}{*}{\begin{tabular}[c]{@{}c@{}}Edge\\ count\end{tabular}} & \multirow{2}{*}{\begin{tabular}[c]{@{}c@{}}Mini\\ route\\ count\end{tabular}} & \multirow{2}{*}{\begin{tabular}[c]{@{}c@{}}Route\\ enumeration\\ time (s)\end{tabular}} & \multirow{2}{*}{\begin{tabular}[c]{@{}c@{}}Vehicle\\ count\end{tabular}} & \multirow{2}{*}{\begin{tabular}[c]{@{}c@{}}Total\\ distance\\ (m)\end{tabular}} & \multirow{2}{*}{\begin{tabular}[c]{@{}c@{}}Vehicle\\ count\\ gap\end{tabular}} & \multirow{2}{*}{\begin{tabular}[c]{@{}c@{}}Optimality\\ gap (\%)\end{tabular}} & \multirow{2}{*}{\begin{tabular}[c]{@{}c@{}}Nodes\\ explored\end{tabular}} & \multicolumn{2}{c}{Wall time (s)\down}\\
			\cline{11-12}\rule{0pt}{20pt} 
			&  &  &  &  &  &  &  &  &  & MIP & \begin{tabular}[c]{@{}c@{}}Optimal\\ count\down\end{tabular}\\
			\hline\up
			M0 & 302 & 14024 & 3233 & 7 & 2 & 481141 & 0 & 0.0 & 109840 & 7200 & 445 \\
			M1 & 262 & 11267 & 8986 & 6 & 3 & 605515 & 1 & 33.3 & 39142 & 7200 & - \\
			M2 & 302 & 13973 & 31559 & 7 & 3 & 847030 & 0 & 0.0 & 27300 & 7200 & 4567 \\
			M3 & 302 & 15253 & 30739 & 10 & 3 & 668490 & 1 & 33.3 & 18968 & 7201 & - \\
			M4 & 302 & 14426 & 28359 & 9 & 3 & 535195 & 1 & 33.3 & 19036 & 7201 & - \\
			M5 & 302 & 12739 & 503 & 6 & 2 & 333048 & 0 & 0.0 & 1348803 & 3409 & 340 \\
			M6 & 302 & 15515 & 47521 & 8 & 3 & 657988 & 1 & 33.3 & 12023 & 7200 & - \\
			M7 & 302 & 14485 & 3485 & 7 & 3 & 595519 & 1 & 33.3 & 123341 & 7200 & - \\
			M8 & 302 & 15404 & 10828 & 8 & 3 & 689147 & 1 & 33.3 & 21890 & 7201 & - \\
			M9 & 302 & 15882 & 55026 & 9 & 3 & 489997 & 1 & 33.3 & 14828 & 7201 & - \\
			M10 & 302 & 14898 & 119198 & 10 & 3 & 719639 & 1 & 33.3 & 18473 & 7200 & - \\
			M11 & 302 & 13800 & 5845 & 10 & 2 & 602968 & 0 & 0.0 & 205444 & 7200 & 1814 \\
			M12 & 302 & 13542 & 1884 & 7 & 2 & 417175 & 0 & 0.0 & 61043 & 1007 & 122 \\
			M13 & 302 & 14564 & 28922 & 9 & 3 & 652724 & 1 & 33.3 & 18510 & 7200 & - \\
			M14 & 302 & 13902 & 3207 & 7 & 2 & 401064 & 0 & 0.0 & 51325 & 2406 & 270 \\
			M15 & 302 & 14801 & 14693 & 7 & 3 & 627967 & 0 & 0.0 & 39332 & 7200 & 7030 \\
			M16 & 254 & 10233 & 3968 & 4 & 3 & 599126 & 0 & 0.0 & 30465 & 2949 & 2787 \\
			M17 & 302 & 13224 & 1380 & 7 & 2 & 490178 & 0 & 0.0 & 14669 & 134 & 73 \\
			M18 & 290 & 11758 & 749 & 5 & 2 & 347259 & 0 & 0.0 & 30780 & 418 & 416 \\
			M19 & 302 & 13043 & 3174 & 7 & 2 & 339073 & 0 & 0.0 & 278853 & 6566 & 6004 \\
			M20 & 302 & 14184 & 4380 & 7 & 3 & 551547 & 1 & 33.3 & 81164 & 7200 & - \\
			M21 & 258 & 10135 & 1696 & 6 & 3 & 620764 & 0 & 0.0 & 273752 & 4256 & 4116 \\
			M22 & 302 & 14856 & 19435 & 8 & 3 & 683612 & 1 & 33.3 & 18247 & 7200 & - \\
			M23 & 302 & 14230 & 12339 & 7 & 3 & 556522 & 1 & 33.3 & 31373 & 7200 & - \\
			M24 & 302 & 14694 & 23970 & 7 & 3 & 588191 & 1 & 33.3 & 18586 & 7200 & - \\
			M25 & 286 & 13139 & 19056 & 6 & 3 & 596412 & 1 & 33.3 & 24223 & 7201 & - \\
			M26 & 302 & 13505 & 1547 & 11 & 3 & 445952 & 0 & 0.0 & 55454 & 1576 & 1311 \\
			M27 & 262 & 10980 & 4981 & 4 & 3 & 712881 & 0 & 0.0 & 34804 & 1648 & 1422 \\
			M28 & 302 & 13883 & 3565 & 11 & 2 & 394323 & 0 & 0.0 & 1737 & 183 & 104 \\
			M29 & 302 & 15142 & 38021 & 10 & 3 & 729149 & 1 & 33.3 & 18677 & 7200 & -\down\\ 
			\hline
		\end{tabular}%
	}
\end{table}

\begin{table}[!th]
	\centering
	\caption{Results of CTSPAV\textsubscript{SEC} and CTSPAV\textsubscript{Base} for the Medium Problem Instances}
	\label{tab:my-table4}
	\setlength{\tabcolsep}{2pt}
	\resizebox{\textwidth}{!}{%
		\begin{tabular}{cccccccc|ccccccc}
			\hline\up
			\multirow{4}{*}{\begin{tabular}[c]{@{}c@{}}Instance\\ name\end{tabular}} & \multicolumn{14}{c}{CTSPAV variant\down} \\
			\cline{2-15} \rule{0pt}{12pt}
			& \multicolumn{7}{c|}{SEC} & \multicolumn{7}{c}{Base\down} \\
			\cline{2-15} \rule{0pt}{14pt}
			& \multirow{2}{*}{\begin{tabular}[c]{@{}c@{}}Vehicle\\ count\end{tabular}} & \multirow{2}{*}{\begin{tabular}[c]{@{}c@{}}Total\\ distance\\ (m)\end{tabular}} & \multirow{2}{*}{\begin{tabular}[c]{@{}c@{}}Vehicle\\ count\\ gap\end{tabular}} & \multirow{2}{*}{\begin{tabular}[c]{@{}c@{}}Optimality\\ gap (\%)\end{tabular}} & \multirow{2}{*}{\begin{tabular}[c]{@{}c@{}}Nodes\\ explored\end{tabular}} & \multicolumn{2}{c|}{Wall time (s)} & \multirow{2}{*}{\begin{tabular}[c]{@{}c@{}}Vehicle\\ count\end{tabular}} & \multirow{2}{*}{\begin{tabular}[c]{@{}c@{}}Total\\ distance\\ (m)\end{tabular}} & \multirow{2}{*}{\begin{tabular}[c]{@{}c@{}}Vehicle\\ count\\ gap\end{tabular}} & \multirow{2}{*}{\begin{tabular}[c]{@{}c@{}}Optimality\\ gap (\%)\end{tabular}} & \multirow{2}{*}{\begin{tabular}[c]{@{}c@{}}Nodes\\ explored\end{tabular}} & \multicolumn{2}{c}{Wall time (s)\down} \\
			\cline{7-8} \cline{14-15} \rule{0pt}{19pt}
			&  &  &  &  &  & MIP & \begin{tabular}[c]{@{}c@{}}Optimal\\ count\end{tabular} &  &  &  &  &  & MIP & \begin{tabular}[c]{@{}c@{}}Optimal\\ count\down\end{tabular} \\
			\hline\up
			M0  & 2 & 480223 & 0 & 0.0  & 229608 & 7200 & 2756 & 2 & 480225 & 1 & 49.9 & 143747  & 7201 & - \\
			M1  & 3 & 603771 & 1 & 33.3 & 21394  & 7203 & -    & 3 & 604103 & 2 & 66.5 & 20833   & 7200 & - \\
			M2  & 3 & 846579 & 1 & 33.2 & 16221  & 7207 & -    & 3 & 846597 & 1 & 33.2 & 21540   & 7201 & - \\
			M3  & 3 & 668248 & 1 & 33.3 & 15377  & 7205 & -    & 3 & 682726 & 2 & 66.5 & 21125   & 7200 & - \\
			M4  & 3 & 535334 & 1 & 33.3 & 7076   & 7298 & -    & 3 & 535195 & 2 & 66.5 & 21423   & 7200 & - \\
			M5  & 2 & 333048 & 0 & 0.0  & 14122  & 335  & 95   & 2 & 333366 & 1 & 49.9 & 1119738 & 7200 & - \\
			M6  & 3 & 656983 & 1 & 33.3 & 6422   & 7201 & -    & 4 & 655969 & 3 & 74.9 & 20905   & 7200 & - \\
			M7  & 3 & 595519 & 1 & 33.3 & 45152  & 7204 & -    & 3 & 595519 & 2 & 66.5 & 65095   & 7201 & - \\
			M8  & 3 & 679167 & 1 & 33.3 & 21493  & 7201 & -    & 3 & 687498 & 2 & 66.5 & 21259   & 7200 & - \\
			M9  & 3 & 489461 & 1 & 33.3 & 5734   & 7238 & -    & 3 & 497878 & 2 & 66.6 & 21032   & 7200 & - \\
			M10 & 3 & 719788 & 1 & 33.3 & 8781   & 7215 & -    & 3 & 722278 & 2 & 66.5 & 3697    & 7201 & - \\
			M11 & 2 & 601111 & 0 & 0.0  & 35354  & 7202 & 1730 & 2 & 601041 & 1 & 49.8 & 29943   & 7200 & - \\
			M12 & 2 & 417175 & 0 & 0.0  & 50401  & 1911 & 195  & 2 & 417185 & 1 & 49.9 & 251358  & 7200 & - \\
			M13 & 3 & 655996 & 1 & 33.3 & 7966   & 7212 & -    & 3 & 653183 & 2 & 66.5 & 21185   & 7202 & - \\
			M14 & 2 & 401064 & 0 & 0.0  & 26183  & 5314 & 983  & 2 & 401064 & 1 & 49.9 & 92983   & 7200 & - \\
			M15 & 4 & 622760 & 2 & 49.9 & 20584  & 7203 & -    & 4 & 622717 & 2 & 49.9 & 23019   & 7200 & - \\
			M16 & 3 & 599126 & 1 & 33.3 & 80256  & 7205 & -    & 3 & 599442 & 2 & 66.5 & 32695   & 7200 & - \\
			M17 & 2 & 490178 & 0 & 0.0  & 4120   & 141  & 58   & 2 & 490178 & 1 & 49.9 & 272323  & 7200 & - \\
			M18 & 2 & 347259 & 0 & 0.0  & 436    & 156  & 151  & 2 & 347259 & 1 & 49.9 & 1064235 & 7201 & - \\
			M19 & 2 & 339073 & 0 & 0.0  & 4695   & 1645 & 637  & 2 & 339073 & 1 & 49.9 & 192573  & 7200 & - \\
			M20 & 3 & 551547 & 1 & 33.3 & 41920  & 7203 & -    & 3 & 551547 & 2 & 66.5 & 39175   & 7200 & - \\
			M21 & 3 & 620783 & 1 & 33.3 & 211984 & 7200 & -    & 3 & 620764 & 2 & 66.5 & 319796  & 7200 & - \\
			M22 & 3 & 685043 & 1 & 33.3 & 15662  & 7205 & -    & 3 & 683300 & 1 & 33.3 & 24972   & 7200 & - \\
			M23 & 3 & 556571 & 1 & 33.3 & 21292  & 7200 & -    & 3 & 555996 & 2 & 66.5 & 20915   & 7200 & - \\
			M24 & 3 & 588191 & 1 & 33.3 & 17174  & 7223 & -    & 3 & 587860 & 2 & 66.5 & 21374   & 7200 & - \\
			M25 & 3 & 597367 & 1 & 33.3 & 20807  & 7200 & -    & 3 & 596653 & 2 & 66.5 & 21527   & 7201 & - \\
			M26 & 3 & 445952 & 1 & 33.3 & 114827 & 7202 & -    & 3 & 445952 & 2 & 66.6 & 189359  & 7200 & - \\
			M27 & 3 & 712881 & 1 & 33.3 & 25439  & 7202 & -    & 3 & 712881 & 2 & 66.5 & 21500   & 7200 & - \\
			M28 & 2 & 394323 & 0 & 0.0  & 1830   & 431  & 241  & 2 & 394323 & 1 & 49.9 & 139970  & 7200 & - \\
			M29 & 3 & 731148 & 1 & 33.3 & 10958  & 7204 & -    & 3 & 729946 & 2 & 66.5 & 19975   & 7203 & -\down\\ 
			\hline
		\end{tabular}%
	}
\end{table}

\begin{table}[!th]
	\centering
	\caption{Results of CTSPAV\textsubscript{Hybrid} for the Tight Problem Instances}
	\label{tab:my-table5}
	\setlength{\tabcolsep}{3pt}
	\resizebox{\textwidth}{!}{%
		\begin{tabular}{cccccccccccc}
			\hline\up
			\multirow{2}{*}{\begin{tabular}[c]{@{}c@{}}Instance\\ name\end{tabular}} & \multirow{2}{*}{\begin{tabular}[c]{@{}c@{}}Node\\ count\end{tabular}} & \multirow{2}{*}{\begin{tabular}[c]{@{}c@{}}Edge\\ count\end{tabular}} & \multirow{2}{*}{\begin{tabular}[c]{@{}c@{}}Mini\\ route\\ count\end{tabular}} & \multirow{2}{*}{\begin{tabular}[c]{@{}c@{}}Route\\ enumeration\\ time (s)\end{tabular}} & \multirow{2}{*}{\begin{tabular}[c]{@{}c@{}}Vehicle\\ count\end{tabular}} & \multirow{2}{*}{\begin{tabular}[c]{@{}c@{}}Total\\ distance\\ (m)\end{tabular}} & \multirow{2}{*}{\begin{tabular}[c]{@{}c@{}}Vehicle\\ count\\ gap\end{tabular}} & \multirow{2}{*}{\begin{tabular}[c]{@{}c@{}}Optimality\\ gap (\%)\end{tabular}} & \multirow{2}{*}{\begin{tabular}[c]{@{}c@{}}Nodes\\ explored\end{tabular}} & \multicolumn{2}{c}{Wall time (s)\down}\\
			\cline{11-12}\rule{0pt}{20pt} 
			&  &  &  &  &  &  &  &  &  & MIP & \begin{tabular}[c]{@{}c@{}}Optimal\\ count\down\end{tabular}\\
			\hline\up
			S0 & 402 & 20870 & 374 & 19 & 5 & 961566 & 0 & 0.0 & 144186 & 544 & 129 \\
			S1 & 402 & 20847 & 267 & 18 & 3 & 619257 & 0 & 0.0 & 19909 & 143 & 124 \\
			S2 & 402 & 21424 & 971 & 20 & 5 & 1246019 & 0 & 0.0 & 27515 & 459 & 333 \\
			S3 & 402 & 21472 & 1268 & 21 & 5 & 1192722 & 0 & 0.0 & 19049 & 830 & 721 \\
			S4 & 402 & 21352 & 1204 & 20 & 5 & 1187914 & 0 & 0.0 & 957524 & 5084 & 238 \\
			S5 & 402 & 20918 & 304 & 17 & 3 & 676142 & 0 & 0.0 & 1887 & 28 & 24 \\
			S6 & 402 & 21050 & 707 & 20 & 6 & 1503404 & 0 & 0.0 & 14494 & 224 & 187 \\
			S7 & 402 & 21022 & 687 & 20 & 5 & 1345009 & 0 & 0.0 & 121198 & 1524 & 1180 \\
			S8 & 402 & 20896 & 581 & 31 & 5 & 1310231 & 0 & 0.0 & 2705 & 37 & 32 \\
			S9 & 402 & 21876 & 1666 & 30 & 6 & 1094536 & 0 & 0.0 & 14475 & 384 & 262 \\
			S10 & 402 & 21044 & 430 & 29 & 4 & 805606 & 0 & 0.0 & 17905 & 228 & 228 \\
			S11 & 402 & 21614 & 835 & 29 & 4 & 819652 & 0 & 0.0 & 11194 & 211 & 188 \\
			S12 & 402 & 20946 & 393 & 32 & 4 & 837723 & 0 & 0.0 & 448878 & 1504 & 86 \\
			S13 & 402 & 21137 & 504 & 20 & 4 & 914708 & 0 & 0.0 & 136179 & 1149 & 667 \\
			S14 & 402 & 21438 & 1056 & 32 & 5 & 1450697 & 0 & 0.0 & 10064 & 71 & 17 \\
			S15 & 402 & 21156 & 2825 & 31 & 5 & 1613836 & 0 & 0.0 & 2646 & 20 & 8 \\
			S16 & 402 & 21005 & 528 & 32 & 5 & 1220586 & 0 & 0.0 & 8396 & 147 & 136 \\
			S17 & 402 & 20844 & 499 & 30 & 5 & 1252397 & 0 & 0.0 & 9523 & 68 & 34 \\
			S18 & 402 & 20713 & 392 & 31 & 6 & 1452716 & 0 & 0.0 & 18044 & 201 & 200 \\
			S19 & 402 & 21377 & 1267 & 31 & 4 & 1030225 & 0 & 0.0 & 513121 & 4069 & 1218 \\
			S20 & 402 & 21542 & 1541 & 33 & 4 & 1144849 & 0 & 0.0 & 8369 & 222 & 211 \\
			S21 & 402 & 20959 & 313 & 19 & 3 & 580008 & 0 & 0.0 & 204235 & 2267 & 2131\down\\ 
			\hline
		\end{tabular}%
	}
\end{table}

\begin{table}[!th]
	\centering
	\caption{Results of CTSPAV\textsubscript{SEC} and CTSPAV\textsubscript{Base} for the Tight Problem Instances}
	\label{tab:my-table6}
	\setlength{\tabcolsep}{2pt}
	\resizebox{\textwidth}{!}{%
		\begin{tabular}{cccccccc|ccccccc}
			\hline\up
			\multirow{4}{*}{\begin{tabular}[c]{@{}c@{}}Instance\\ name\end{tabular}} & \multicolumn{14}{c}{CTSPAV variant\down} \\
			\cline{2-15} \rule{0pt}{12pt}
			& \multicolumn{7}{c|}{SEC} & \multicolumn{7}{c}{Base\down} \\
			\cline{2-15} \rule{0pt}{14pt}
			& \multirow{2}{*}{\begin{tabular}[c]{@{}c@{}}Vehicle\\ count\end{tabular}} & \multirow{2}{*}{\begin{tabular}[c]{@{}c@{}}Total\\ distance\\ (m)\end{tabular}} & \multirow{2}{*}{\begin{tabular}[c]{@{}c@{}}Vehicle\\ count\\ gap\end{tabular}} & \multirow{2}{*}{\begin{tabular}[c]{@{}c@{}}Optimality\\ gap (\%)\end{tabular}} & \multirow{2}{*}{\begin{tabular}[c]{@{}c@{}}Nodes\\ explored\end{tabular}} & \multicolumn{2}{c|}{Wall time (s)} & \multirow{2}{*}{\begin{tabular}[c]{@{}c@{}}Vehicle\\ count\end{tabular}} & \multirow{2}{*}{\begin{tabular}[c]{@{}c@{}}Total\\ distance\\ (m)\end{tabular}} & \multirow{2}{*}{\begin{tabular}[c]{@{}c@{}}Vehicle\\ count\\ gap\end{tabular}} & \multirow{2}{*}{\begin{tabular}[c]{@{}c@{}}Optimality\\ gap (\%)\end{tabular}} & \multirow{2}{*}{\begin{tabular}[c]{@{}c@{}}Nodes\\ explored\end{tabular}} & \multicolumn{2}{c}{Wall time (s)\down} \\
			\cline{7-8} \cline{14-15} \rule{0pt}{19pt}
			&  &  &  &  &  & MIP & \begin{tabular}[c]{@{}c@{}}Optimal\\ count\end{tabular} &  &  &  &  &  & MIP & \begin{tabular}[c]{@{}c@{}}Optimal\\ count\down\end{tabular} \\
			\hline\up
			S0  & 5 & 961566  & 0 & 0.0  & 95266  & 388  & 82   & 5 & 961566  & 0 & 0.0  & 151291 & 533  & 90   \\
			S1  & 3 & 619257  & 0 & 0.0  & 9643   & 97   & 89   & 3 & 619257  & 0 & 0.0  & 24230  & 326  & 323  \\
			S2  & 5 & 1246019 & 0 & 0.0  & 13952  & 277  & 203  & 5 & 1246019 & 0 & 0.0  & 21299  & 917  & 902  \\
			S3  & 5 & 1192722 & 1 & 20.0 & 178946 & 7201 & -    & 5 & 1192722 & 1 & 19.9 & 241540 & 7201 & -    \\
			S4  & 5 & 1187914 & 0 & 0.0  & 400941 & 2668 & 187  & 5 & 1187914 & 0 & 0.0  & 17315  & 406  & 225  \\
			S5  & 3 & 676142  & 0 & 0.0  & 3023   & 13   & 5    & 3 & 676142  & 0 & 0.0  & 4393   & 14   & 6    \\
			S6  & 6 & 1503404 & 0 & 0.0  & 14190  & 284  & 284  & 6 & 1503404 & 0 & 0.0  & 73967  & 1653 & 1567 \\
			S7  & 5 & 1345009 & 0 & 0.0  & 216353 & 3000 & 2352 & 5 & 1345009 & 0 & 0.0  & 243780 & 2824 & 1953 \\
			S8  & 5 & 1310231 & 0 & 0.0  & 1459   & 22   & 21   & 5 & 1310231 & 0 & 0.0  & 3948   & 46   & 44   \\
			S9  & 6 & 1094536 & 1 & 16.6 & 152214 & 7202 & -    & 6 & 1094536 & 1 & 16.6 & 193079 & 7201 & -    \\
			S10 & 4 & 805606  & 0 & 0.0  & 16966  & 236  & 226  & 4 & 805606  & 0 & 0.0  & 9222   & 84   & 80   \\
			S11 & 4 & 819652  & 0 & 0.0  & 9997   & 168  & 150  & 4 & 819652  & 0 & 0.0  & 12619  & 210  & 197  \\
			S12 & 4 & 837723  & 0 & 0.0  & 161991 & 665  & 99   & 4 & 837723  & 0 & 0.0  & 155274 & 554  & 157  \\
			S13 & 4 & 914708  & 0 & 0.0  & 94311  & 1553 & 1301 & 4 & 914708  & 0 & 0.0  & 304917 & 5579 & 5231 \\
			S14 & 5 & 1450697 & 0 & 0.0  & 1250   & 44   & 31   & 5 & 1450697 & 0 & 0.0  & 14449  & 39   & 7    \\
			S15 & 5 & 1613836 & 0 & 0.0  & 3338   & 24   & 12   & 5 & 1613836 & 0 & 0.0  & 1600   & 19   & 11   \\
			S16 & 5 & 1220586 & 0 & 0.0  & 3471   & 78   & 72   & 5 & 1220586 & 0 & 0.0  & 1348   & 54   & 52   \\
			S17 & 5 & 1252397 & 0 & 0.0  & 7338   & 48   & 32   & 5 & 1252397 & 0 & 0.0  & 10428  & 76   & 55   \\
			S18 & 6 & 1452716 & 0 & 0.0  & 26594  & 278  & 268  & 6 & 1452716 & 0 & 0.0  & 22340  & 375  & 374  \\
			S19 & 4 & 1030225 & 0 & 0.0  & 553629 & 4371 & 1324 & 4 & 1030225 & 0 & 0.0  & 173744 & 1736 & 326  \\
			S20 & 4 & 1144849 & 0 & 0.0  & 9894   & 199  & 188  & 4 & 1144849 & 0 & 0.0  & 22515  & 504  & 500  \\
			S21 & 3 & 580008  & 0 & 0.0  & 138098 & 2078 & 2042 & 3 & 580008  & 1 & 33.3 & 886186 & 7201 & -\down\\ 
			\hline
		\end{tabular}%
	}
                \end{table}

Tables \ref{tab:my-table9}, \ref{tab:my-table10}, and
\ref{tab:my-table11} list the heuristic results for every large,
medium, and tight instance respectively. Their first columns show the
instance names, followed by three columns that show the number of
columns (mini routes) generated, the final vehicle count, and the
total travel distance for every instance. The following two columns
display the absolute gap of its vehicle count results and the
optimality gap of its best incumbent solution. Since the heuristic
does not utilize all feasible mini routes, it has to use the optimal
LP-relaxation solution of RMP\textsubscript{CTSPAV} to derive primal
lower bounds for these gap calculations. The final three columns show
the percentage difference between the column count, the vehicle count,
and the total distance of the heuristic relative to those of
CTSPAV\textsubscript{Hybrid}.

                \begin{table}[!th]
	\centering
	\caption{Results of CTSPAV Column-Generation Heuristic by \cite{hasan2021} for Large Problem Instances}
	\label{tab:my-table9}
		\begin{tabular}{crcrccrrr}
			\hline\up
			\multirow{2}{*}{\begin{tabular}[c]{@{}c@{}}Instance\\ name\end{tabular}} & \multirow{2}{*}{\begin{tabular}[c]{@{}c@{}}Column\\ count\end{tabular}} & \multicolumn{1}{c}{\multirow{2}{*}{\begin{tabular}[c]{@{}c@{}}Vehicle\\ count\end{tabular}}} & \multicolumn{1}{c}{\multirow{2}{*}{\begin{tabular}[c]{@{}c@{}}Total\\ distance\\ (m)\end{tabular}}} & \multirow{2}{*}{\begin{tabular}[c]{@{}c@{}}Vehicle\\ count\\ gap\end{tabular}} & \multirow{2}{*}{\begin{tabular}[c]{@{}c@{}}Optimality\\ gap (\%)\end{tabular}} & \multicolumn{3}{c}{Percentage difference\down} \\
			\cline{7-9} \rule{0pt}{20pt}
			&  & \multicolumn{1}{c}{} & \multicolumn{1}{c}{} &  &  & \multicolumn{1}{c}{\begin{tabular}[c]{@{}c@{}}Column\\ count\end{tabular}} & \multicolumn{1}{c}{\begin{tabular}[c]{@{}c@{}}Vehicle\\ count\end{tabular}} & \multicolumn{1}{c}{\begin{tabular}[c]{@{}c@{}}Total\\ distance\down\end{tabular}} \\
			\hline\up
			L0 & 2231 & 3 & 647661 & 2 & 66.5 & -40\% & 0\% & -0.75\% \\
			L1 & 901 & 3 & 463065 & 2 & 66.6 & -18\% & 0\% & 0.00\% \\
			L2 & 8713 & 4 & 817348 & 3 & 74.9 & -83\% & 0\% & 0.05\% \\
			L3 & 9347 & 4 & 841180 & 3 & 74.8 & -85\% & 0\% & -0.30\% \\
			L4 & 7253 & 4 & 813018 & 3 & 74.8 & -85\% & 0\% & -0.14\% \\
			L5 & 960 & 3 & 512675 & 2 & 66.6 & -40\% & 0\% & 0.00\% \\
			L6 & 6330 & 4 & 955285 & 3 & 74.8 & -68\% & 0\% & 1.19\% \\
			L7 & 5087 & 4 & 888490 & 3 & 74.8 & -76\% & 0\% & 0.35\% \\
			L8 & 4902 & 4 & 844674 & 3 & 74.8 & -67\% & 0\% & 0.00\% \\
			L9 & 13892 & 5 & 737361 & 4 & 79.9 & -93\% & 0\% & -0.22\% \\
			L10 & 2884 & 3 & 555102 & 2 & 66.5 & -50\% & 0\% & 0.05\% \\
			L11 & 5659 & 3 & 570036 & 2 & 66.5 & -69\% & 0\% & 0.88\% \\
			L12 & 2116 & 3 & 581863 & 2 & 66.5 & -42\% & 0\% & 0.00\% \\
			L13 & 3106 & 3 & 624843 & 2 & 66.5 & -55\% & 0\% & 0.01\% \\
			L14 & 9539 & 4 & 949361 & 3 & 74.8 & -87\% & 0\% & -0.66\% \\
			L15 & 8161 & 4 & 1108007 & 3 & 74.8 & -93\% & 0\% & 0.80\% \\
			L16 & 3513 & 4 & 847394 & 3 & 74.8 & -62\% & 0\% & 0.37\% \\
			L17 & 2886 & 4 & 862155 & 3 & 74.8 & -55\% & 0\% & 0.11\% \\
			L18 & 2912 & 4 & 914762 & 3 & 74.8 & -34\% & 0\% & 0.49\% \\
			L19 & 6278 & 3 & 698599 & 2 & 74.9 & -82\% & 33\% & 0.27\% \\
			L20 & 8291 & 3 & 779684 & 2 & 66.5 & -86\% & 0\% & -1.40\% \\
			L21 & 1397 & 2 & 457911 & 1 & 49.9 & -71\% & 0\% & -0.01\%\down\\
			\hline
		\end{tabular}
\end{table}

\begin{table}[!th]
	\centering
	\caption{Results of CTSPAV Column-Generation Heuristic by \cite{hasan2021} for Medium Problem Instances}
	\label{tab:my-table10}
		\begin{tabular}{crcrccrrr}
			\hline\up
			\multirow{2}{*}{\begin{tabular}[c]{@{}c@{}}Instance\\ name\end{tabular}} & \multirow{2}{*}{\begin{tabular}[c]{@{}c@{}}Column\\ count\end{tabular}} & \multicolumn{1}{c}{\multirow{2}{*}{\begin{tabular}[c]{@{}c@{}}Vehicle\\ count\end{tabular}}} & \multicolumn{1}{c}{\multirow{2}{*}{\begin{tabular}[c]{@{}c@{}}Total\\ distance\\ (m)\end{tabular}}} & \multirow{2}{*}{\begin{tabular}[c]{@{}c@{}}Vehicle\\ count\\ gap\end{tabular}} & \multirow{2}{*}{\begin{tabular}[c]{@{}c@{}}Optimality\\ gap (\%)\end{tabular}} & \multicolumn{3}{c}{Percentage difference\down} \\
			\cline{7-9} \rule{0pt}{20pt}
			&  & \multicolumn{1}{c}{} & \multicolumn{1}{c}{} &  &  & \multicolumn{1}{c}{\begin{tabular}[c]{@{}c@{}}Column\\ count\end{tabular}} & \multicolumn{1}{c}{\begin{tabular}[c]{@{}c@{}}Vehicle\\ count\end{tabular}} & \multicolumn{1}{c}{\begin{tabular}[c]{@{}c@{}}Total\\ distance\down\end{tabular}} \\
			\hline\up
			M0 & 1664 & 2 & 481141 & 1 & 49.9 & -49\% & 0\% & -0.19\% \\
			M1 & 2643 & 3 & 605515 & 2 & 66.5 & -71\% & 0\% & -0.17\% \\
			M2 & 4349 & 3 & 846579 & 2 & 66.5 & -86\% & 0\% & 0.75\% \\
			M3 & 5461 & 3 & 668490 & 2 & 66.5 & -82\% & 0\% & 1.07\% \\
			M4 & 3556 & 3 & 535195 & 2 & 66.6 & -87\% & 0\% & 0.04\% \\
			M5 & 464 & 2 & 333048 & 1 & 49.9 & -8\% & 0\% & 0.00\% \\
			M6 & 6217 & 3 & 657988 & 2 & 66.5 & -87\% & 0\% & 0.36\% \\
			M7 & 2081 & 3 & 595519 & 2 & 66.5 & -40\% & 0\% & 0.00\% \\
			M8 & 3728 & 3 & 689147 & 2 & 66.5 & -66\% & 0\% & -0.21\% \\
			M9 & 6545 & 3 & 489997 & 2 & 66.6 & -88\% & 0\% & 0.01\% \\
			M10 & 6938 & 3 & 719639 & 2 & 66.5 & -94\% & 0\% & 0.03\% \\
			M11 & 2142 & 2 & 602968 & 1 & 49.9 & -63\% & 0\% & -0.40\% \\
			M12 & 1198 & 2 & 417175 & 1 & 49.9 & -36\% & 0\% & 0.00\% \\
			M13 & 4821 & 3 & 652724 & 2 & 66.5 & -83\% & 0\% & 0.17\% \\
			M14 & 1712 & 2 & 401064 & 1 & 49.9 & -47\% & 0\% & 0.08\% \\
			M15 & 4122 & 3 & 627967 & 2 & 74.9 & -72\% & 33\% & -1.07\% \\
			M16 & 1849 & 3 & 599126 & 2 & 66.5 & -53\% & 0\% & 0.07\% \\
			M17 & 964 & 2 & 490178 & 1 & 49.9 & -30\% & 0\% & 0.00\% \\
			M18 & 528 & 2 & 347259 & 1 & 49.9 & -30\% & 0\% & 0.00\% \\
			M19 & 914 & 2 & 339073 & 1 & 49.9 & -71\% & 0\% & 0.00\% \\
			M20 & 2153 & 3 & 551547 & 2 & 66.5 & -51\% & 0\% & 0.00\% \\
			M21 & 1172 & 3 & 620764 & 2 & 66.5 & -31\% & 0\% & 0.00\% \\
			M22 & 4527 & 3 & 683612 & 2 & 66.5 & -77\% & 0\% & 0.16\% \\
			M23 & 3416 & 3 & 556522 & 2 & 66.5 & -72\% & 0\% & -0.09\% \\
			M24 & 4949 & 3 & 588191 & 2 & 66.5 & -79\% & 0\% & 0.04\% \\
			M25 & 3969 & 3 & 596412 & 2 & 66.5 & -79\% & 0\% & 0.16\% \\
			M26 & 1043 & 3 & 445952 & 2 & 66.6 & -33\% & 0\% & 0.09\% \\
			M27 & 2336 & 3 & 712881 & 2 & 66.5 & -53\% & 0\% & 0.02\% \\
			M28 & 1810 & 2 & 394323 & 1 & 49.9 & -49\% & 0\% & 0.00\% \\
			M29 & 6028 & 3 & 729149 & 2 & 66.5 & -84\% & 0\% & -0.38\%\down\\
			\hline
		\end{tabular}
\end{table}

\begin{table}[!th]
	\centering
	\caption{Results of CTSPAV Column-Generation Heuristic by \cite{hasan2021} for Tight Problem Instances}
	\label{tab:my-table11}
		\begin{tabular}{crcrccrrr}
			\hline\up
			\multirow{2}{*}{\begin{tabular}[c]{@{}c@{}}Instance\\ name\end{tabular}} & \multirow{2}{*}{\begin{tabular}[c]{@{}c@{}}Column\\ count\end{tabular}} & \multicolumn{1}{c}{\multirow{2}{*}{\begin{tabular}[c]{@{}c@{}}Vehicle\\ count\end{tabular}}} & \multicolumn{1}{c}{\multirow{2}{*}{\begin{tabular}[c]{@{}c@{}}Total\\ distance\\ (m)\end{tabular}}} & \multirow{2}{*}{\begin{tabular}[c]{@{}c@{}}Vehicle\\ count\\ gap\end{tabular}} & \multirow{2}{*}{\begin{tabular}[c]{@{}c@{}}Optimality\\ gap (\%)\end{tabular}} & \multicolumn{3}{c}{Percentage difference\down} \\
			\cline{7-9} \rule{0pt}{20pt}
			&  & \multicolumn{1}{c}{} & \multicolumn{1}{c}{} &  &  & \multicolumn{1}{c}{\begin{tabular}[c]{@{}c@{}}Column\\ count\end{tabular}} & \multicolumn{1}{c}{\begin{tabular}[c]{@{}c@{}}Vehicle\\ count\end{tabular}} & \multicolumn{1}{c}{\begin{tabular}[c]{@{}c@{}}Total\\ distance\down\end{tabular}} \\
			\hline\up
			S0 & 359 & 5 & 961566 & 2 & 39.9 & -4\% & 0\% & 0.00\% \\
			S1 & 267 & 3 & 619257 & 1.97 & 65.4 & 0\% & 0\% & 0.00\% \\
			S2 & 813 & 5 & 1246019 & 2 & 39.9 & -16\% & 0\% & 0.00\% \\
			S3 & 969 & 5 & 1192722 & 2.92 & 58.3 & -24\% & 0\% & 0.00\% \\
			S4 & 840 & 5 & 1187914 & 1 & 20.0 & -30\% & 0\% & 0.05\% \\
			S5 & 291 & 3 & 676142 & 1 & 33.3 & -4\% & 0\% & 0.00\% \\
			S6 & 633 & 6 & 1503404 & 2 & 33.3 & -10\% & 0\% & 0.00\% \\
			S7 & 575 & 5 & 1345009 & 2 & 49.9 & -16\% & 20\% & -1.19\% \\
			S8 & 502 & 5 & 1310231 & 1 & 19.9 & -14\% & 0\% & 0.00\% \\
			S9 & 1273 & 6 & 1094536 & 3 & 49.9 & -24\% & 0\% & 0.13\% \\
			S10 & 421 & 4 & 805606 & 2 & 49.9 & -2\% & 0\% & 0.00\% \\
			S11 & 744 & 4 & 819652 & 1 & 25.0 & -11\% & 0\% & 0.01\% \\
			S12 & 377 & 4 & 837723 & 2 & 49.9 & -4\% & 0\% & 0.00\% \\
			S13 & 461 & 4 & 914708 & 2 & 49.9 & -9\% & 0\% & 0.12\% \\
			S14 & 876 & 5 & 1450697 & 1.62 & 32.2 & -17\% & 0\% & 0.01\% \\
			S15 & 1072 & 5 & 1613836 & 1 & 19.9 & -62\% & 0\% & 0.01\% \\
			S16 & 502 & 5 & 1220586 & 2 & 39.9 & -5\% & 0\% & 0.00\% \\
			S17 & 453 & 5 & 1252397 & 2 & 39.9 & -9\% & 0\% & 0.00\% \\
			S18 & 383 & 6 & 1452716 & 2 & 33.3 & -2\% & 0\% & 0.00\% \\
			S19 & 762 & 4 & 1030225 & 1.50 & 37.4 & -40\% & 0\% & 0.15\% \\
			S20 & 934 & 4 & 1144849 & 1 & 24.9 & -39\% & 0\% & 0.00\% \\
			S21 & 305 & 3 & 580008 & 2 & 66.5 & -3\% & 0\% & 0.00\%\down\\
			\hline
		\end{tabular}
\end{table}

\end{document}